\documentclass[leqno]{amsart}
\usepackage[foot]{amsaddr}
\usepackage[T1]{fontenc}
\usepackage[utf8]{inputenc} %
\usepackage{lmodern}
\usepackage{csquotes}

\usepackage[a4paper,margin=1in]{geometry}
\usepackage{amsmath,amssymb,amsthm,mathtools}
\usepackage{bm,bbm}
\usepackage{booktabs,array,longtable}
\usepackage{enumitem}
\usepackage{xcolor}
\usepackage{float}
\usepackage{graphicx}
\usepackage{microtype}
\usepackage{caption}
\usepackage[ruled,vlined,linesnumbered]{algorithm2e}
\SetNlSty{}{}{}
\usepackage[hidelinks,hypertexnames=false]{hyperref}
\usepackage[nameinlink,capitalize]{cleveref}
\usepackage{tikz}
\usetikzlibrary{arrows.meta,positioning}
\newtheorem{theorem}{Theorem}[section]
\newtheorem{proposition}[theorem]{Proposition}

\newtheorem{assumption}[theorem]{Assumption}
\newtheorem{definition}[theorem]{Definition}
\theoremstyle{remark}
\newtheorem{remark}[theorem]{Remark}

\definecolor{SaifColor}{RGB}{34,139,34}
\definecolor{ErikColor}{RGB}{0,130,160}
\definecolor{RaulColor}{RGB}{180,45,55}
\definecolor{ChristianColor}{RGB}{100,65,190}

\newcommand{\R}{\mathbb R}
\newcommand{\E}{\mathbb E}
\newcommand{\Pbb}{\mathbb P}
\newcommand{\mcal}{\mathcal}
\newcommand{\F}{\mathbb F}
\newcommand{\Yhist}{\mathbf Y}
\newcommand{\yhist}{\mathbf y}
\newcommand{\cA}{\mcal A}
\newcommand{\cB}{\mcal B}
\newcommand{\cP}{\mcal P}
\newcommand{\cV}{\mcal V}
\newcommand{\beliefstate}{\mathsf b}

\newcommand{\cG}{\mcal G}
\newcommand{\cJ}{\mcal J}
\newcommand{\dd}{\mathrm d}
\newcommand{\mb}[1]{\mathbb{#1}}
\newcommand{\mf}[1]{\mathbf{#1}}
\newcommand{\ip}[2]{\left\langle #1,#2\right\rangle}
\newcommand{\tH}{\widetilde{\mcal H}}
\newcommand{\Uhat}{\widehat U}

\title[Hamiltonian Potential Optimization for POSOC]{Hamiltonian Potential Optimization for Partially Observed Stochastic Optimal Control}

\author[C. Bayer, S. Ben Naamia, E. von Schwerin, R. Tempone]{Christian Bayer$^{1}$, Saifeddine ben Naamia$^{\#,2}$,  Erik von Schwerin$^{3}$ \& Ra\'{u}l Tempone$^{3}$}
\address{$^{1}$Weierstrass Institut (WIAS), Berlin, Germany,}
\address{$^2$Department of Mathematics, RWTH Aachen University, Geb\"{a}ude-1953 1.OG, Pontdriesch 14-16, 161, 52062 Aachen, Germany}
\address{$^3$King Abdullah University of Science \& Technology (KAUST), Computer, Electrical and Mathematical Sciences \& Engineering Division (CEMSE), Thuwal 23955-6900, Saudi Arabia}
\email{$^\#$Bennaamia@ssd.rwth-aachen.de}


\date{\today}

\newcommand{\AnsParams}{348}

\newcommand{\AnsFeatPoly}{2, 8, 21, 45}
\newcommand{\AnsFeatTrig}{2, 6, 12, 20}
\newcommand{\AnsTimeCoef}{3}
\newcommand{\AnsJ}{1.4313}
\newcommand{\AnsJse}{0.0149}
\newcommand{\AnsD}{1.4207}
\newcommand{\AnsGap}{0.0106}
\newcommand{\AnsGapPct}{0.7}

\newcommand{\AnsSteps}{100}
\newcommand{\AnsDualLR}{0.01}

\newcommand{\AnsRandScale}{0.1}
\newcommand{\AnsEvalPaths}{5\,000}
\newcommand{\AnsTrainIters}{200}
\newcommand{\AnsTrainLR}{0.01}
\newcommand{\AnsHJBPaired}{1.4293}

\newcommand{\AnsPairedDiff}{$+0.0020$}
\newcommand{\AnsPairedDiffse}{0.0027}

\newcommand{\AnsDualSE}{0.0008}

\newcommand{\AnsDualMCPaths}{5\,000}
\newcommand{\AnsBetaOneAns}{0.760}
\newcommand{\AnsBetaOneHJB}{0.783}

\begin{document}

\begin{abstract}
We develop a Hamiltonian potential optimization method for partially observed stochastic optimal control with continuous-time state dynamics and discrete-time noisy observations, including controlled sensing. Rather than solving the infinite-dimensional belief-space Hamilton–Jacobi–Bellman equation, we optimize over parametrized potentials depending on the state variable and the available observation history. The associated Hamiltonian minimizers induce admissible continuous controls and, when observations are controlled, sensing decisions. An extended Hamiltonian incorporating a cost multiplier provides a common representation of the induced belief dynamics and running and observation costs. This construction yields a finite-dimensional optimization problem, with gradients computed by a forward simulation and a backward discrete adjoint sweep. The numerical implementation uses Gaussian beliefs represented by conditional means and covariances. For nonsmooth Hamiltonians arising from switching controls, we introduce smooth Hamiltonian approximations to enable gradient-based optimization. To assess the resulting policies, we develop belief-space information relaxation lower bounds. A separate optimization of the dual test functional provides a lower bound which, paired with the admissible policy cost, yields an a posteriori primal–dual diagnostic. Numerical experiments investigate scalar and multidimensional linear–quadratic–Gaussian problems, controlled sensing with and without control bounds, a nonquadratic obstacle penalty, and Hamiltonian smoothing for a discontinuous control selector.
\end{abstract}

\maketitle
\section*{Acknowledgments}
C. Bayer acknowledges support from DFG CRC/TRR 388 \enquote{Rough Analysis, Stochastic Dynamics and Related Fields} Project B03 and DFG individual grant number 497300407 \enquote{Recursive and sparse approximation in reinforcement learning with applications}.
S. Ben naamia acknowledges support from  the Deutsche Forschungsgemeinschaft (DFG, German Research Foundation) -- 333849990/GRK2379 (IRTG Hierarchical and Hybrid Approaches in Modern Inverse Problems).
This work was supported by the King Abdullah University of Science and Technology (KAUST) Office of Sponsored Research (OSR) under Award  No.~OSR-2019-CRG8-4033.

\section{Introduction}
\label{sec:introduction}

Partially observed stochastic optimal control (POSOC) studies decision-making in stochastic systems where the controller must choose actions based on noisy or incomplete observations, received continuously or at discrete times, rather than direct knowledge of the full state.  This setting appears naturally in controlled sensing, target tracking, robotics under noisy localization, medical and biological decision-making, inventory systems with latent demand, and finance with hidden factors.  The central difficulty is that the controller is affected not only by stochastic state dynamics, but also by the stochastic evolution of its own information.

A classical way to restore the Markov property is to replace the hidden state by the conditional distribution of the state given the observation history.  This conditional distribution, often called the belief state, is a sufficient statistic  for control under standard Markov assumptions.  The partially observed problem is then transformed into a fully observed stochastic control problem on a space of probability measures~\cite{Bensoussan1992,Krishnamurthy2016}.  This reformulation is conceptually powerful, but it also exposes the main computational obstacle: the Hamilton--Jacobi--Bellman equation is posed on an infinite-dimensional belief space~\cite{fabbri2017stochastic}.  Even when the physical state is low dimensional, direct dynamic programming on the belief space is usually impractical unless the model has additional structure.
Neural approximation methods for infinite-dimensional Hamilton--Jacobi--Bellman equations
provide another computational approach, including the Deep Hilbert--Galerkin
methods in~\cite{CohenDeFeoHebnerSirignano2026}.

We consider a hybrid setting with continuous-time hidden dynamics and discrete-time observations, as in the companion Hamilton–Jacobi–Bellman (HJB) work~\cite{bayer2025cont}, where the value is characterized by continuous-time belief-space HJB equations coupled through Bayesian update steps at the observation times.  In Gaussian regimes, this formulation can sometimes be reduced to finite-dimensional variables such as conditional means and covariances.  Nevertheless, the dynamic-programming viewpoint remains demanding, because it requires solving a value equation, or a finite-dimensional reduction of it, on the relevant belief state space.

A complementary analytical route is provided by maximum-principle methods.  Classical partially observed maximum principles are developed in~\cite{Bensoussan1983,Haussmann1987,BarasElliottKohlmann1989,Tang1998}.  More recently, the companion work~\cite{bayer2025pontryagin} derived a Pontryagin maximum principle (PMP) directly on the belief space for continuous-time dynamics with discrete observations. It also proposed a particle-based fixed-point iteration, combining policy extraction, forward belief simulation and regression updates to approximate the resulting optimality system. That result gives a forward–backward system of necessary optimality conditions and states the relationship between the belief-space adjoint, the derivative of the value functional, and the Hamiltonians, under suitable regularity assumptions.  Thus the HJB formulation gives a value-based characterization, while the PMP formulation gives a variational characterization through adjoint equations.

A related setting is adaptive control with continuously observed state
dynamics whose drift depends on an unobserved random parameter
\cite{CohenKnochenhauerMerkel2023}. Although the observation structure
differs from the discrete-time observations considered here, controls
influence both state evolution and learning about the unknown drift.
This creates a similar coupling between control and information acquisition.

Related computational approaches include mean-field formulations of
memory-limited partially observable control and forward--backward
sweeps for the associated HJB--Fokker--Planck system
\cite{TottoriKobayashi2022,TottoriKobayashi2023}, as well as approaches
combining data assimilation with mean-field Pontryagin formulations
for partially observed digital twins~\cite{OpperReich2026,Reich2026Online}.

Here, we develop a computational formulation based on Hamiltonians and parametrized potentials. Instead of solving the belief-space HJB equation directly, and instead of solving a full forward--backward PMP system, we introduce a pathwise potential
\begin{equation}
	U_t(x,\yhist^{[n]}),
\end{equation}
where \(x\) is the hidden state and \(\yhist^{[n]}=(y_1,\ldots,y_n)\) is the observation history available between the \(n\)-th and \((n+1)\)-st observation times.  A candidate potential is inserted into continuous and discrete reduced Hamiltonians.  The corresponding pointwise reductions induce the continuous control between observations and, when present, the sensing or observation-design variable before the next observation is received.  Thus every admissible potential defines a partially observed feedback policy and an on-policy belief flow.

We therefore optimize over a parametrized family of potentials rather than directly over arbitrary adapted policies,
\begin{equation}
	\min_{\theta}\; \cJ_\Delta(\widehat U^\theta),
\end{equation}
where \(\widehat U^\theta\) is a potential parametrized by the finite-dimensional vector \(\theta\), and \(\cJ_\Delta\) is its induced cost in the discretized simulator.  The simulator propagates the belief, applies Bayesian or Kalman updates at observation times, evaluates the reduced Hamiltonian formulas, and accumulates the running, observation, and terminal costs.  In this way, the method keeps the Hamiltonian structure of the analytical theory but turns the computation into a finite-dimensional training problem.

We call a Hamiltonian reduced once the minimization over the corresponding control variable has been carried out.  In the linear–quadratic–Gaussian (LQG) examples, this gives explicit formulas in terms of Gaussian means and variances.  In the controlled-sensing regime, the pre-posterior optimization produces scalar formulas for the sensing variable.  In the bounded-control regime, the reduction produces a Hamiltonian with an absolute-value structure, motivating the Hamiltonian smoothing developed below.  The potential is the object being trained. The physical controls are recovered from the reduced Hamiltonian formulas when the trained potential is evaluated.

The main advantage of this viewpoint is that it separates the modeling,
analytical, and computational tasks. The underlying formulation is not tied
to a particular representation of the belief. In the present paper,
however, the numerical scope is deliberately restricted to Gaussian belief
representations. The belief is represented by its conditional mean and
covariance. All numerical examples use linear-Gaussian state and observation
models, so Gaussian filtering conditional on the full observation history is
exact. The obstacle example introduces a nonquadratic running cost, which
does not change this filtering property. Time discretization remains a
numerical approximation. This structure keeps the belief recursion finite
dimensional and permits analytical evaluation of the Gaussian moments and
probabilities used below.

Particle and ensemble representations could in principle be combined with
the same Hamiltonian-potential formulation. Their treatment would require
additional choices concerning filtering, resampling, quadrature, and
gradient estimation. These questions are not addressed here and are left
for separate work.

A second goal of the manuscript is to provide diagnostics for the learned potential.  A trained potential produces an admissible policy and therefore a primal upper bound on the optimal cost.  To assess how informative this upper bound is, we also develop belief-space information relaxation dual bounds, building on~\cite{Rogers2007,HaughKogan2004,BrownSmithSun2010,BankDeFeo2026}.  The resulting primal--dual gap is an a posteriori diagnostic for the induced policy.
Weak duality holds for every admissible test functional, so the dual test functional can be
chosen independently of the primal potential and optimized separately.
The exact dual functional provides a theoretical lower bound. Its numerical
evaluation also depends on inner minimization and integration or sampling
errors, so the reported numerical values are diagnostics unless those
errors are controlled in the required direction.

We also address a practical difficulty in potential optimization: Hamiltonian minimization can produce a nonsmooth reduced Hamiltonian and a discontinuous control selector. In the affine-in-control switching example, the minimizing control switches between endpoints of the control set, making gradient-based optimization difficult. We introduce smooth Hamiltonian approximations to obtain differentiable control selectors. Numerical experiments demonstrate improved training stability and assess the resulting policies using the original control cost and the associated primal--dual diagnostic.

The main contributions are:
\begin{itemize}[leftmargin=2em]
	\item A Hamiltonian potential optimization method for partially observed control, with induced admissible policies and discrete adjoint gradients for differentiable belief updates.
	\item Belief-space information relaxation bounds for policy assessment, including an exact quadratic dual test family for the scalar controlled-sensing LQG setting studied here.
	\item Smooth Hamiltonian approximations for gradient-based training in problems with discontinuous switching controls.
	\item Numerical validation on fixed-noise and controlled-sensing LQG problems, including bounded controls and dimensions up to $d=100$, a nonquadratic obstacle problem, and a switching-control problem.
\end{itemize}

The paper is organized as follows.
Section~\ref{sec:problem} states the partially observed control problem and
the belief dynamics.
Section~\ref{sec:method} introduces the reduced Hamiltonian potential
formulation and explains the relation between the potential-induced cost and
the original admissible control cost.
Section~\ref{sec:duality} develops the belief-space information relaxation
lower bound and the associated primal--dual certificate.
Section~\ref{sec:training-workflow} turns the reduced Hamiltonian formulation
into a differentiable finite-dimensional training problem.
Section~\ref{sec:numerical-method} gives the numerical discretization and the
discrete adjoint gradient formula, together with smoothing of nonsmooth reduced
Hamiltonians.
Section~\ref{sec:numerics} contains the numerical experiments.
The appendix contains derivations, potential parametrizations (Appendix~\ref{sec:parametrization}), scalar Gaussian formulas (Appendix~\ref{sec:gaussian-formulas}), and implementation details.

\section{Partially observed control problem}
\label{sec:problem}

We consider a finite time horizon $[0,T]$ and a probability
space $(\Omega,\mathcal F,\mathbb P)$ supporting a
$d_w$-dimensional Brownian motion $W$ and an initial random
variable $X_0\sim\mu_0$, independent of $W$.
The state process $X_t\in\mb R^{d_x}$ evolves according to
the controlled stochastic differential equation
\begin{equation}
	\dd X_t
	= b(t,X_t,\alpha_t)\,\dd t
	+ \sigma(t,X_t,\alpha_t)\,\dd W_t,
	\qquad X_0\sim\mu_0,
	\label{eq:state-sde}
\end{equation}
where
\[
b:[0,T]\times\mb R^{d_x}\times\mb R^{d_\alpha}
\to\mb R^{d_x},
\qquad
\sigma:[0,T]\times\mb R^{d_x}\times\mb R^{d_\alpha}
\to\mb R^{d_x\times d_w}
\]
are the drift and diffusion coefficients, respectively,
and $\alpha_t\in\mb R^{d_\alpha}$ is the continuous
control.

The controller receives observations
$Y_{t_n}\in\mb R^{d_y}$ at deterministic times
\[
0=t_0<t_1<\cdots<t_{N_o}<t_{N_o+1}=T.
\]
At each observation time $t_n$, a sensing or design
control $\beta_{t_n}\in\mb R^{d_\beta}$ is selected
before the observation is revealed. The observation
model is given by
\begin{equation}
	Y_{t_n}
	=h_n(t_n,X_{t_n},\beta_{t_n},\xi_n),
	\qquad n=1,\ldots,N_o,
	\label{eq:observation-model}
\end{equation}
where $\xi_n\in\mb R^{d_\xi}$ represents measurement
noise and
\[
h_n:[0,T]\times\mb R^{d_x}
\times\mb R^{d_\beta}\times\mb R^{d_\xi}
\to\mb R^{d_y}
\]
is a measurable observation function.

Conditional on the state $X_{t_n}=x$, the previously
observed history, and the chosen sensing control
$\beta_{t_n}=\beta$, we assume that $Y_{t_n}$ admits
a conditional density $\pi_n$, so that
\begin{equation}
	\mathbb P\left(
	Y_{t_n}\in\dd y
	\,\middle|\,
	X_{t_n}=x,\,
	\mf Y^{[n-1]}=\mf y^{[n-1]},\,
	\beta_{t_n}=\beta
	\right)
	=
	\pi_n(y\mid x,\mf y^{[n-1]},\beta)\,\dd y.
	\label{eq:likelihood}
\end{equation}
Here, $\mf y^{[n-1]}$ denotes a realization of the
past observation history, defined below. The explicit
dependence on the observation history allows for
history-dependent likelihoods.

\subsection{Information filtrations}

Let
\[
\mf Y^{[n]}:=\{Y_{t_1},\ldots,Y_{t_n}\},
\]
and, for $t\in[0,T]$, define
\[
\nu(t)
:=\max\{k\in\{1,\ldots,N_o\}:t_k\le t\},
\]
with the conventions $\max\emptyset=0$ and
$\mf Y^{[0]}:=\emptyset$.

We introduce the observation and state filtrations
through the $\sigma$-algebras
\begin{align}
	\mb F_t^Y
	&:=\sigma\bigl(\mf Y^{[\nu(t)]}\bigr),
	&& t\in[0,T],
	\label{eq:filtration-Y}\\
	\mb F_t^X
	&:=\sigma(X_s:0\le s\le t),
	&& t\in[0,T].
	\label{eq:filtration-X}
\end{align}
In particular, $\mb F_t^Y$ is piecewise constant
between consecutive observation times and changes
only when a new observation becomes available.

\subsection{Controls and admissibility}

We adopt a closed-loop (feedback) formulation in which
controls are selected using the information available
to the controller. We distinguish between continuous
controls $\alpha_t\in\mb R^{d_\alpha}$ and discrete
sensing controls $\beta_{t_n}\in\mb R^{d_\beta}$.
The former are progressively measurable with respect
to the relevant information filtration, while the
latter are selected immediately before the
corresponding observations.

Let $\mcal A$ and $\mcal B$ denote the prescribed
constraints on the continuous and sensing controls,
respectively. Throughout, admissibility includes these constraints,
well-posedness of the controlled state equation, and integrability of
the running, sensing and terminal cost terms introduced below.

\begin{definition}[Admissible control classes]
	\label{def:admissible}
	For $t\in[0,T]$, we distinguish the following
	two information structures:
	\begin{align*}
		\mcal U^{\mathrm{PO}}[t,T]
		:=\Big\{
		&u=(\alpha,\beta)\in
		\mcal A\times\mcal B:
		\ \alpha\text{ is }\mb F^Y
		\text{-progressively measurable on }[t,T],\\
		&\beta_{t_n}\text{ is }
		\mb F^Y_{t_n^-}\text{-measurable
			for every }t_n\in[t,T]
		\Big\},
		\\[0.5em]
		\mcal U^{\mathrm{FO}}[t,T]
		:=\Big\{
		&u=(\alpha,\beta)\in
		\mcal A\times\mcal B:
		\ \alpha\text{ is }\mb F^X
		\text{-progressively measurable on }[t,T],\\
		&\beta_{t_n}\text{ is }
		\mb F^X_{t_n^-}\text{-measurable
			for every }t_n\in[t,T]
		\Big\}.
	\end{align*}
	Here,
	\[
	\mb F_{t_n^-}^Y
	:=\sigma\bigl(\mf Y^{[n-1]}\bigr),
	\qquad
	\mb F_{t_n^-}^X
	:=\sigma(X_s:0\le s<t_n).
	\]
	The \emph{partially observed} class
	$\mcal U^{\mathrm{PO}}$ represents the admissible
	policies for the original control problem.
	The \emph{fully observed} class
	$\mcal U^{\mathrm{FO}}$ consists of state-feedback
	policies and will be used to construct a lower
	envelope for the partially observed value function.
\end{definition}

\subsection{Performance criterion}

For an admissible policy $u=(\alpha,\beta)$, we consider
the expected total cost
\begin{equation}
	J(\alpha,\beta)
	:=
	\E\left[
	\int_0^T f_t(X_t,\alpha_t)\,\dd t
	+\sum_{n=1}^{N_o}
	c_{t_n}(X_{t_n},\beta_{t_n})
	+g(X_T)
	\right],
	\label{eq:original-cost}
\end{equation}
where $f_t:\mb R^{d_x}\times\mb R^{d_\alpha}
\to\mb R$ is the running cost,
$c_{t_n}:\mb R^{d_x}\times\mb R^{d_\beta}
\to\mb R$ is the sensing cost at time $t_n$,
and $g:\mb R^{d_x}\to\mb R$ is the terminal cost.

The objective is to minimize the expected total cost
over all admissible partially observed policies:
\begin{equation}
	V_0^{\mathrm{PO}}
	:=
	\inf_{(\alpha,\beta)\in
		\mcal U^{\mathrm{PO}}[0,T]}
	J(\alpha,\beta).
	\label{eq:optimal-value}
\end{equation}

\subsection{Belief dynamics}

The belief process is the conditional law
\[
  \mu_t = \mcal L(X_t\mid \F^Y_t).
\]
Here $\F_t^Y$ is the observation filtration defined in~\eqref{eq:filtration-Y}. Between observations, $\mu_t$ is propagated by prediction. If $\varphi$ is a smooth test function and $\cG_\alpha$ denotes the generator of~\eqref{eq:state-sde}, then
\begin{equation}
  \frac{\dd}{\dd t}\ip{\varphi}{\mu_t}
  = \ip{\cG_{\alpha_t}\varphi}{\mu_t},
  \qquad t\in(t_n,t_{n+1}).
  \label{eq:weak-filter-prediction}
\end{equation}
At an observation time, the forecast belief $\mu_{t_n^-}$ is updated by Bayes' rule:
\begin{equation}
  \mu_{t_n}(\dd x)
  = \frac{\pi_n(Y_{t_n}\mid x,\Yhist^{[n-1]},\beta_{t_n})\,\mu_{t_n^-}(\dd x)}
  {\int_{\R^{d_x}} \pi_n(Y_{t_n}\mid x',\Yhist^{[n-1]},\beta_{t_n})\,\mu_{t_n^-}(\dd x')}.
  \label{eq:bayes-update}
\end{equation}
Conditioning the cost on the observation filtration gives the belief-space form
\begin{equation}
  J(\alpha,\beta)
  = \E\bigg[
    \int_0^T \ip{f_t(\cdot,\alpha_t)}{\mu_t}\,\dd t
    +\sum_{n=1}^{N_o}\ip{c_{t_n}(\cdot,\beta_{t_n})}{\mu_{t_n^-}}
    +\ip{g}{\mu_T}
  \bigg].
  \label{eq:belief-cost}
\end{equation}

\begin{assumption}[Standing regularity]
The coefficients in~\eqref{eq:state-sde} are regular enough for the state equation to be well posed under admissible policies. The likelihoods in~\eqref{eq:likelihood} are measurable, nonnegative, normalized, and chosen so that~\eqref{eq:bayes-update} is well defined. The partially observed admissible class is nonempty. The costs and state moments satisfy the integrability conditions needed for the conditioning and expectation identities below. Whenever derivatives of Hamiltonians are used, we assume either smooth reduced Hamiltonians or a smooth approximation of the type described in~\cref{subsec:smoothing}.
\end{assumption}

\section{Reduced Hamiltonian potential formulation}
\label{sec:method}

This section introduces the reduced Hamiltonians used by the algorithm. The guiding convention is that Hamiltonians are evaluated after the pointwise minimization over the relevant continuous or sensing variable has been carried out. Thus the training algorithm works directly with reduced Hamiltonian formulas and their induced numerical maps.
\subsection{Continuous reduced Hamiltonian}

Let $U:\R^{d_x}\to\R$ be a smooth test potential and let $\gamma\ge 0$ mark the running-cost contribution. The augmented continuous belief Hamiltonian is
\begin{equation}
  \tH^c_t(\mu,U,\gamma)
  := \inf_{a\in\cA}
  \ip{\cG_aU+\gamma f_t(\cdot,a)}{\mu}.
  \label{eq:continuous-aug-H}
\end{equation}
In every concrete experiment below, this infimum is either computed analytically or replaced by a smooth reduced expression. The cost slice is obtained by differentiating with respect to $\gamma$ on the value slice $\gamma=1$:
\begin{equation}
  \mathfrak f_t(\mu,U)
  :=\partial_\gamma \tH^c_t(\mu,U,1).
  \label{eq:continuous-cost-slice}
\end{equation}
This is the running cost of the Hamiltonian-induced dynamics.

\subsection{Discrete reduced Hamiltonian}

At an observation time $t_{n+1}$, the sensing decision is chosen before the next observation is known. Given a forecast belief $\mu^-$ and a post-observation potential $U^+(x,(\yhist^{[n]},y))$, define
\begin{align}
  \tH^d_{n+1}(\mu^-,U^+,\gamma)
  := \inf_{\beta\in\cB}\bigg\{&
  \gamma\ip{c_{t_{n+1}}(\cdot,\beta)}{\mu^-}
  +\int_{\R^{d_y}}
  \int_{\R^{d_x}}U^+(x,(\yhist^{[n]},y))
  \pi_{n+1}(y\mid x,\yhist^{[n]},\beta)\,\mu^-(\dd x)\,\dd y
  \bigg\}.
  \label{eq:discrete-aug-H}
\end{align}
The discrete cost slice is
\begin{equation}
  \mathfrak c_{n+1}(\mu^-,U^+)
  :=\partial_\gamma \tH^d_{n+1}(\mu^-,U^+,1).
  \label{eq:discrete-cost-slice}
\end{equation}
For fixed-noise observation models, such as regime A below, $\mathfrak c_n\equiv 0$. For controlled sensing models, such as regime B, this term contains the sensing cost.

\begin{proposition}[Cost slices from augmented Hamiltonians]
	\label{prop:cost-slice-main}
	For each augmented Hamiltonian, suppose that the action set is independent
	of $\gamma$, the objective and its $\gamma$-derivative are jointly
	continuous in the action and $\gamma$, and minimizers exist in a common
	compact subset of the action set for $\gamma$ near $1$.
	If the minimizer at $\gamma=1$ is unique, then
	$\partial_\gamma\tH^c_t(\mu,U,1)$ equals the running-cost contribution
	at that minimizer. The analogous derivative of $\tH^d_{n+1}$ equals
	the sensing-cost contribution.
	The same conclusions hold with multiple minimizers if all minimizers
	at $\gamma=1$ have the same corresponding cost contribution.
	No differentiability of the minimizing selector is required.
	The proof is given in \cref{app:cost-slice-proof}.
\end{proposition}

\subsection{Potential-induced simulation}

Let
\[
  \Uhat_t(x,\yhist^{[n]}),\qquad t\in[t_n,t_{n+1}),
\]
be a candidate pathwise potential. Through the reduced Hamiltonians, the potential determines:
\begin{enumerate}[leftmargin=2em]
  \item a continuous prediction map for the belief between observations.
  \item a pre-posterior sensing map at observation times when sensing is controlled.
  \item a Bayes or Kalman update map once the realized observation is sampled or read.
  \item continuous and discrete cost slices.
\end{enumerate}
The resulting on-policy cost can be written abstractly as
\begin{equation}
  \cJ(U)
  = \E\bigg[
  \int_0^T \mathfrak f_t(\mu^U_t,U_t)\,\dd t
  +\sum_{n=1}^{N_o}\mathfrak c_n(\mu^U_{t_n^-},U_{t_n})
  +\ip{g}{\mu^U_T}
  \bigg].
  \label{eq:on-policy-H-cost}
\end{equation}
This expression is not meant to introduce a second cost functional unrelated to the original partially observed control problem. Its role is to write the original cost after the reduced Hamiltonian has selected, from the potential, the continuous and discrete decisions used by the closed-loop simulator. This relation is made explicit in the next subsection.

\subsection{Relation with the original admissible control cost}
\label{subsec:cJ-versus-J}

A candidate potential does not automatically generate an admissible
partially observed policy. In particular, the minimizing branches of the
reduced Hamiltonians must exist, admit measurable selections, and respect
the information available when each decision is made.

Let $\cV$ be a candidate class of sufficiently regular potentials. For
$U\in\cV$, suppose that the continuous Hamiltonian minimization induces a
control $\alpha^U$ and, when sensing is controlled, that the discrete
pre-posterior Hamiltonian minimization induces controls
$\beta_{t_n}^U$, $n=1,\ldots,N_o$.

We call $U$ an admissible potential if its induced controls satisfy:
\begin{enumerate}
	\item $\alpha^U$ is $\cA$-valued and progressively measurable with
	respect to the observation filtration
	$(\mcal F_t^Y)_{t\in[0,T]}$.
	\item for every $n=1,\ldots,N_o$, $\beta_{t_n}^U$ is $\cB$-valued
	and measurable with respect to $\mcal F_{t_n^-}^Y$.
\end{enumerate}
We additionally require that the induced closed-loop state--belief
system is well posed and that $J(\alpha^U,\beta^U)$ is finite.
We denote the corresponding class by
\begin{equation}
	\cV_{\mathrm{adm}}
	:=
	\left\{
	U\in\cV:
	(\alpha^U,\beta^U)
	\text{ satisfies the preceding admissibility conditions}
	\right\}.
	\label{eq:admissible-potential-class}
\end{equation}
Thus admissibility is not a consequence of regularity of $U$ alone. It is
a property of the Hamiltonian selectors and of the closed-loop system
induced by $U$.

For $U\in\cV_{\mathrm{adm}}$, the superscript $U$ only indicates that the
controls are generated by the potential. The actual numerical algorithm
evaluates the reduced Hamiltonians and their induced closed-loop maps
directly.

The cost slices in \eqref{eq:on-policy-H-cost} are chosen to reproduce the
original running and observation costs under this induced closed loop.
More precisely, whenever the envelope identification in
Proposition~\ref{prop:cost-slice-main} applies to the selected minimizing branches,
\begin{equation}
	\mathfrak f_t(\mu_t^U,U_t)
	=
	\ip{f_t(\cdot,\alpha_t^U)}{\mu_t^U},
	\label{eq:cJ-J-running-identification}
\end{equation}
and, at an observation time,
\begin{equation}
	\mathfrak c_n(\mu_{t_n^-}^U,U_{t_n})
	=
	\ip{c_{t_n}(\cdot,\beta_{t_n}^U)}{\mu_{t_n^-}^U}.
	\label{eq:cJ-J-discrete-identification}
\end{equation}
Consequently,
\begin{equation}
	\boxed{
		\mcal J(U)
		=
		J(\alpha^U,\beta^U)
	},
	\qquad U\in\cV_{\mathrm{adm}}.
	\label{eq:cJ-equals-J-induced}
\end{equation}
The proof is given in \cref{app:cJ-versus-J-proof}. In regimes without
controlled sensing, the notation reduces to
$\cJ(U)=J(\alpha^U)$, and the discrete cost slice is absent or identically
zero.

Identity \eqref{eq:cJ-equals-J-induced} gives the potential objective its
control interpretation. We are not minimizing a surrogate cost unrelated
to the original problem. We are evaluating the original POSOC cost along
the admissible policies generated by the chosen potential class.

Once the measurable minimizing selections have been fixed, define the
induced policy class by
\begin{equation}
	\mcal U_{\cV}
	:=
	\left\{
	(\alpha^U,\beta^U):
	U\in\cV_{\mathrm{adm}}
	\right\}.
	\label{eq:potential-induced-policy-class}
\end{equation}
By the definition of $\cV_{\mathrm{adm}}$,
\begin{equation}
	\mcal U_{\cV}\subseteq\mcal U,
\end{equation}
where $\mcal U$ is the full class of admissible partially observed
policies. It follows that
\begin{equation}
	\inf_{U\in\cV_{\mathrm{adm}}}\cJ(U)
	=
	\inf_{(\alpha,\beta)\in\mcal U_{\cV}}
	J(\alpha,\beta)
	\geq
	\inf_{(\alpha,\beta)\in\mcal U}
	J(\alpha,\beta).
	\label{eq:restricted-inf-relation}
\end{equation}
The Hamiltonian potential optimization problem is therefore to find
\begin{equation}
	\overline U
	\in
	\operatorname*{argmin}_{U\in\cV_{\mathrm{adm}}}
	\cJ(U).
	\label{eq:continuous-potential-problem}
\end{equation}
In computation, one considers a finite-dimensional family
$\{\widehat U^\theta:\theta\in\Theta\}$ and restricts the optimization to
\begin{equation}
	\Theta_{\mathrm{adm}}
	:=
	\left\{
	\theta\in\Theta:
	\widehat U^\theta\in\cV_{\mathrm{adm}}
	\right\}.
	\label{eq:admissible-parameter-set}
\end{equation}
In the numerical examples below, the control formulas,
action constraints, and information structure are chosen so that
admissibility can be checked directly from the induced closed-loop maps.

\begin{remark}[Restricted policy interpretation]
	Problem~\eqref{eq:continuous-potential-problem} is a restricted policy
	optimization problem induced by admissible potentials. It is an exact
	reformulation of the unrestricted POSOC problem only if the induced policy
	class attains, or approximates arbitrarily closely, the unrestricted
	optimal value. The numerical experiments therefore report the original
	cost $J$ under the trained admissible policy, together with the associated
	duality diagnostic.
\end{remark}

\section{Duality and information relaxation on belief space}
\label{sec:duality}

An admissible potential $U\in\cV_{\mathrm{adm}}$ generates a partially
observed policy and therefore an upper bound on the optimal cost:
\begin{equation}
	V_0(\mu_0)\leq J(\alpha^U,\beta^U)=\cJ(U).
	\label{eq:duality-primal-upper-bound}
\end{equation}
This section constructs a complementary lower bound by information
relaxation. We allow a hypothetical controller to know future observation
randomness, while subtracting a penalty whose expectation is zero for
every original admissible policy. The relaxed controller minimizes over
a larger class, so its penalized value provides a lower bound.
The construction follows~\cite{Rogers2007,BrownSmithSun2010} and the
belief-state approach of~\cite{HaughRuizLacedelli2020}.

The primal potential $U$ and the dual test functional $\Phi$ have different
roles: $U$ generates the policy being evaluated, whereas $\Phi$ determines
the penalty. They may be parametrized and optimized independently.

\subsection{Dual test functionals and observation-update penalty}
\label{subsec:general-dual-test-functionals}

Let
\begin{equation}
	\Phi_n:[t_n,t_{n+1}]\times\cP(\R^{d_x})
	\times(\R^{d_y})^n\longrightarrow\R,
	\qquad n=0,\ldots,N_o,
	\label{eq:duality-test-family}
\end{equation}
be measurable test functionals with the integrability needed below.
The history $\yhist^{[n]}$ is fixed between observations. Under a Gaussian
belief, one may use a function of the conditional mean and covariance:
\begin{equation}
	\Phi_n(t,\mu;\yhist^{[n]})
	=\phi_n\bigl(t,m(\mu),\Sigma(\mu);\yhist^{[n]}\bigr).
	\label{eq:dual-gaussian-functional}
\end{equation}

\subsubsection{Predictive observations and Bayesian updates}

For a forecast belief $\mu^-$ before $t_n$, past history
$\yhist^{[n-1]}$, and sensing action $\beta\in\cB$, define
\begin{equation}
	L_n(y;\yhist^{[n-1]},\beta,\mu^-)
	:=\int_{\R^{d_x}}\pi_n(y\mid x,\yhist^{[n-1]},\beta)
	\,\mu^-(\dd x).
	\label{eq:dual-predictive-density}
\end{equation}
The Bayes operator is
\begin{equation}
	K_{\beta,y}(\mu^-)(\dd x)
	:=\frac{\pi_n(y\mid x,\yhist^{[n-1]},\beta)}
	{L_n(y;\yhist^{[n-1]},\beta,\mu^-)}\,\mu^-(\dd x).
	\label{eq:dual-bayes-operator}
\end{equation}
Its dependence on the date and past history is suppressed. The operator
is used where the predictive density is positive. Elsewhere a measurable
version may be assigned without affecting the predictive integrals.

The pre-posterior average is the expected test value after the next
observation, computed before that observation is received:
\begin{equation}
	\begin{aligned}
		(\mcal E_n^\beta\Phi_n)(\mu^-;\yhist^{[n-1]})
		:=\int_{\R^{d_y}}&\Phi_n\bigl(t_n,K_{\beta,y}(\mu^-);
		(\yhist^{[n-1]},y)\bigr)\\
		&\times L_n(y;\yhist^{[n-1]},\beta,\mu^-)\,\dd y.
	\end{aligned}
	\label{eq:dual-preposterior-average}
\end{equation}

\subsubsection{A centered observation-update penalty}

For an original admissible policy $u=(\alpha,\beta)$, denote its belief
by $\mu^u$ and its observation history by $\Yhist^{u,[n]}$.
The realized test value minus its pre-posterior average is
\begin{equation}
	\begin{aligned}
		\Delta M_n^{u,\Phi}:={}&
		\Phi_n\bigl(t_n,\mu_{t_n}^u;\Yhist^{u,[n]}\bigr)\\
		&-(\mcal E_n^{\beta_{t_n}}\Phi_n)
		\bigl(\mu_{t_n^-}^u;\Yhist^{u,[n-1]}\bigr).
	\end{aligned}
	\label{eq:dual-centered-jump}
\end{equation}
Define the total penalty by
\begin{equation}
	M_Y^\Phi(u):=\sum_{n=1}^{N_o}\Delta M_n^{u,\Phi}.
	\label{eq:MY-main}
\end{equation}
Because the sensing action is $\mb F_{t_n^-}^Y$-measurable and the
pre-posterior average uses the correct predictive law,
\begin{equation}
	\E\bigl[\Delta M_n^{u,\Phi}\mid\mb F_{t_n^-}^Y\bigr]=0.
	\label{eq:dual-centered-jump-zero}
\end{equation}
Consequently, $\E[M_Y^\Phi(u)]=0$ for every original admissible policy.

\subsection{Observation-seed relaxation and lower bound}
\label{subsec:dual-predictive-seeds}
\label{subsec:DY-main}

We now specify the future randomness revealed to the relaxed controller.
When sensing changes the observation law, candidate actions must be
evaluated using a specified common-randomness construction.

\subsubsection{A common random scenario}

Consider first a scalar Gaussian forecast and observation model:
\[
X\mid\mb F_{t_n^-}^Y\sim\mathcal N(m^-,s^-),\qquad
Y=X+\beta\varepsilon,\qquad \varepsilon\sim\mathcal N(0,1),
\]
where $\beta$ is chosen from the pre-observation information and
$\varepsilon$ is independent of the forecast state and that information.
The predictive observation law is $\mathcal N(m^-,s^-+\beta^2)$, which
can be generated by
\begin{equation}
	Y^\beta=m^-+\sqrt{s^-+\beta^2}\,\xi_n,
	\qquad \xi_n\sim\mathcal N(0,1).
	\label{eq:dual-gaussian-predictive-seed}
\end{equation}
Here $\xi_n$ randomizes the complete predictive observation, accounting
for both forecast uncertainty and sensor noise. It is distinct from the
physical sensor noise $\varepsilon$. We hold $\xi_n$ fixed while comparing
candidate sensing actions, allowing each action to generate its own
observation.\footnote{Fixing the observation is generally different from
	fixing the random input. If $m^-=s^-=0$, then $Y^\beta=\beta\xi_n$.
	The same seed $\xi_n=1$ gives $Y^1=1$ and $Y^2=2$. Conversely, imposing
	$Y^\beta=1$ requires different seeds for $\beta=1$ and $\beta=2$.}

For a general model, let $\xi=(\xi_1,\ldots,\xi_{N_o})$ be independent
predictive-observation seeds with prescribed laws $\lambda_n$.
Their joint law is fixed independently of the policy. Under every
original admissible policy, $\xi_n$ is independent of $\mb F_{t_n^-}^Y$.
We generate observations using a measurable map
\begin{equation}
	Y_{t_n}=G_n\bigl(\mu_{t_n^-},\Yhist^{[n-1]},\beta_{t_n},\xi_n\bigr).
	\label{eq:dual-observation-seed-representation}
\end{equation}
For every fixed forecast belief, history, and sensing action, the map
$G_n$ must transform $\xi_n\sim\lambda_n$ into an observation with density
$L_n(\cdot;\yhist^{[n-1]},\beta,\mu^-)$. Thus the simulator reproduces
the original predictive transition laws under admissible policies.
A fixed seed sequence specifies a common scenario for all candidate
policies. This coupling is fixed throughout the dual calculation.

\subsubsection{Policy notation and additional information}

For a generic policy $u=(\alpha,\beta)$ and scenario $\xi$, write its
realized actions as $\alpha_t(\xi)$ and $\beta_{t_n}(\xi)$, its belief
as $\mu_t^u(\xi)$, and its observations as $Y_{t_n}^u(\xi)$. Set
\[
\Yhist^{u,[n]}(\xi)
:=\bigl(Y_{t_1}^u(\xi),\ldots,Y_{t_n}^u(\xi)\bigr),
\qquad \Yhist^{u,[0]}(\xi):=\emptyset.
\]
The argument $\xi$ identifies the scenario. Admissibility determines
which information the policy may use. Superscripts $U$ remain reserved
for potential-induced policies, whose realized actions are
$\alpha_t^U(\xi)$ and $\beta_{t_n}^U(\xi)$.

An original policy uses only available observations. The relaxed class
$\mcal U^{Y,\mathrm{rel}}$ additionally permits decisions to depend on
the complete seed sequence. An original policy can ignore this extra
information, so
\begin{equation}
	\mcal U^{\mathrm{PO}}[0,T]\subseteq\mcal U^{Y,\mathrm{rel}}.
	\label{eq:admissible-contained-relaxed}
\end{equation}
Only the information constraint is relaxed: the action constraints and
feasibility requirements for the separated-state dynamics are retained.

For each candidate policy and scenario, we use the original prediction
equation and the same observation-generation and update maps:
\begin{align}
	Y_{t_n}^u(\xi)
	&=G_n\bigl(\mu_{t_n^-}^u(\xi),\Yhist^{u,[n-1]}(\xi),
	\beta_{t_n}(\xi),\xi_n\bigr),
	\label{eq:relaxed-generated-observation}\\
	\mu_{t_n}^u(\xi)
	&=K_{\beta_{t_n}(\xi),Y_{t_n}^u(\xi)}
	\bigl(\mu_{t_n^-}^u(\xi)\bigr).
	\label{eq:relaxed-belief-update}
\end{align}
The history argument of $K$ is suppressed. Under anticipative policies,
$\mu_t^u(\xi)$ denotes the separated state propagated by these maps.
It is not defined by conditioning the hidden state on the revealed seeds.

\subsubsection{Pathwise cost and penalty}

Define
\begin{equation}
	\begin{aligned}
		C(u;\xi):={}&\int_0^T
		\ip{f_t(\cdot,\alpha_t(\xi))}{\mu_t^u(\xi)}\,\dd t\\
		&+\sum_{n=1}^{N_o}
		\ip{c_{t_n}(\cdot,\beta_{t_n}(\xi))}{\mu_{t_n^-}^u(\xi)}
		+\ip{g}{\mu_T^u(\xi)}.
	\end{aligned}
	\label{eq:dual-pathwise-cost}
\end{equation}
Evaluate the penalty along the same candidate trajectory:
\begin{equation}
	\begin{aligned}
		M_Y^\Phi(u;\xi):=\sum_{n=1}^{N_o}\Big[&
		\Phi_n\bigl(t_n,\mu_{t_n}^u(\xi);\Yhist^{u,[n]}(\xi)\bigr)\\
		&-(\mcal E_n^{\beta_{t_n}(\xi)}\Phi_n)
		\bigl(\mu_{t_n^-}^u(\xi);\Yhist^{u,[n-1]}(\xi)\bigr)\Big].
	\end{aligned}
	\label{eq:dual-pathwise-penalty}
\end{equation}
The pre-posterior average always uses the original predictive kernel,
holding the candidate forecast belief, history, and sensing action fixed.
It is not conditioned on the revealed future seeds.

For every original admissible policy, transition-law consistency and
the centering identity give
\begin{equation}
	\E_\xi[C(u;\xi)]=J(u),\qquad
	\E_\xi[M_Y^\Phi(u;\xi)]=0.
	\label{eq:penalty-preserves-admissible-cost}
\end{equation}
The zero-mean property is required only for original admissible policies.

\subsubsection{The information relaxation lower bound}

For each scenario, minimize the penalized cost over relaxed policies and
then average the resulting infima:
\begin{equation}
	D_Y^\Phi(0,\mu_0):=\E_\xi\left[
	\inf_{u\in\mcal U^{Y,\mathrm{rel}}}
	\{C(u;\xi)-M_Y^\Phi(u;\xi)\}\right].
	\label{eq:DY-main}
\end{equation}

\begin{theorem}[Observation-relaxation lower bound]
	\label{thm:DY-certificate-main}
	Suppose that:
	\begin{enumerate}
		\item $\mcal U^{\mathrm{PO}}[0,T]\subseteq\mcal U^{Y,\mathrm{rel}}$.
		\item for every original admissible policy, $C(u;\xi)$ and
		$M_Y^\Phi(u;\xi)$ are integrable and satisfy
		\eqref{eq:penalty-preserves-admissible-cost}.
		\item the pathwise infimum in \eqref{eq:DY-main} is measurable and
		integrable.
	\end{enumerate}
	Then
	\begin{equation}
		D_Y^\Phi(0,\mu_0)\leq V_0(\mu_0).
		\label{eq:DY-leq-V-main}
	\end{equation}
\end{theorem}

\begin{proof}
	Fix $u\in\mcal U^{\mathrm{PO}}[0,T]$. Since this policy is available
	to the relaxed optimizer, almost surely,
	\[
	\inf_{v\in\mcal U^{Y,\mathrm{rel}}}
	\{C(v;\xi)-M_Y^\Phi(v;\xi)\}
	\leq C(u;\xi)-M_Y^\Phi(u;\xi).
	\]
	Taking expectations gives $D_Y^\Phi(0,\mu_0)\leq J(u)$.
	Taking the infimum over original admissible policies proves the claim.
\end{proof}

Combining weak duality with the admissible policy induced by $U$ yields
\begin{equation}
	\boxed{D_Y^\Phi(0,\mu_0)\leq V_0(\mu_0)
		\leq J(\alpha^U,\beta^U)=\cJ(U).}
	\label{eq:duality-sandwich-main}
\end{equation}
Exactness of the test functional is not needed for this inequality.
Its quality determines how informative the lower bound is.

\subsection{Exactness of the observation relaxation}

The exact value functions provide a zero-gap penalty when they satisfy
the following pathwise dynamic-programming conditions.

\begin{proposition}[Zero-gap property for the exact value]
	\label{prop:dual-exactness-at-value}
	Let $V_n(t,\mu;\yhist^{[n]})$, $n=0,\ldots,N_o$, be the belief-space
	value functions, with $t_0=0$ and $t_{N_o+1}=T$. Assume the hypotheses
	of \cref{thm:DY-certificate-main} hold with $\Phi=V$. Suppose also that:
	\begin{enumerate}
		\item $V_{N_o}(T,\mu;\yhist^{[N_o]})=\ip{g}{\mu}$.
		\item along every relaxed prediction trajectory with the scenario
		fixed, the integrated verification inequality holds:
		\begin{equation}
			\begin{aligned}
				&V_n\bigl(t,\mu_t^u(\xi);\Yhist^{u,[n]}(\xi)\bigr)\\
				&\quad\leq\int_t^s
				\ip{f_r(\cdot,\alpha_r(\xi))}{\mu_r^u(\xi)}\,\dd r
				+V_n\bigl(s,\mu_s^u(\xi);\Yhist^{u,[n]}(\xi)\bigr),
			\end{aligned}
			\label{eq:dual-value-continuous-verification}
		\end{equation}
		for $t_n\leq t<s<t_{n+1}$, including the pre-observation endpoint
		with $\mu_{t_{n+1}^-}^u(\xi)$ and the terminal endpoint with $\mu_T^u(\xi)$.
		\item at observation times,
		\begin{equation}
			\begin{aligned}
				V_{n-1}(t_n,\mu^-;\yhist^{[n-1]})
				=\inf_{\beta\in\cB}\Big\{&\ip{c_{t_n}(\cdot,\beta)}{\mu^-}\\
				&+(\mcal E_n^\beta V_n)(\mu^-;\yhist^{[n-1]})\Big\}.
			\end{aligned}
			\label{eq:dual-value-jump-bellman}
		\end{equation}
	\end{enumerate}
	The pathwise conditions hold simultaneously for all relaxed controls
	on a common set of scenarios of probability one. Then
	\begin{equation}
		D_Y^V(0,\mu_0)=V_0(\mu_0).
		\label{eq:dual-zero-gap-at-value}
	\end{equation}
\end{proposition}

\begin{proof}
	Fix a scenario in the common set and a relaxed policy, and suppress
	their labels in the trajectory notation. On each prediction interval,
	\[
	\begin{aligned}
		\int_{t_n}^{t_{n+1}}\ip{f_t(\cdot,\alpha_t)}{\mu_t}\,\dd t
		\geq{}&V_n(t_n,\mu_{t_n};\Yhist^{[n]})\\
		&-V_n(t_{n+1},\mu_{t_{n+1}^-};\Yhist^{[n]}),
	\end{aligned}
	\]
	where at the terminal endpoint $\mu_{t_{N_o+1}^-}$ denotes $\mu_T$.
	At observation time $t_n$, write
	\[
	\Delta M_n^V=V_n(t_n,\mu_{t_n};\Yhist^{[n]})
	-(\mcal E_n^{\beta_{t_n}}V_n)(\mu_{t_n^-};\Yhist^{[n-1]}).
	\]
	The Bellman identity implies, for the candidate action,
	\[
	\begin{aligned}
		\ip{c_{t_n}(\cdot,\beta_{t_n})}{\mu_{t_n^-}}-\Delta M_n^V
		\geq{}&V_{n-1}(t_n,\mu_{t_n^-};\Yhist^{[n-1]})\\
		&-V_n(t_n,\mu_{t_n};\Yhist^{[n]}).
	\end{aligned}
	\]
	Summing cancels all intermediate value terms. The terminal condition
	gives $C(u;\xi)-M_Y^V(u;\xi)\geq V_0(\mu_0)$ for every relaxed policy
	on the same set of scenarios. Taking the pathwise infimum and
	expectation gives $D_Y^V(0,\mu_0)\geq V_0(\mu_0)$.
	Weak duality gives the reverse inequality.
\end{proof}

The integrated condition expresses that following a feasible control
segment and then acting optimally cannot improve on the initial value.
It avoids differentiability assumptions on $V$. The proposition is
conditional on verification of these properties for the chosen model.

\subsection{Potential-based test functionals and an exact LQG certificate}

\subsubsection{Lifting a pathwise potential to belief space}

A convenient test family is the linear belief lift
\begin{equation}
	\Phi_n^U(t,\mu;\yhist^{[n]})
	:=\int_{\R^{d_x}}U(t,x,\yhist^{[n]})\,\mu(\dd x).
	\label{eq:phi-from-potential}
\end{equation}
Substituting the Bayes operator in the pre-posterior average gives
\begin{equation}
	\begin{aligned}
		(\mcal E_n^\beta\Phi_n^U)(\mu^-;\yhist^{[n-1]})
		=\int_{\R^{d_y}}\int_{\R^{d_x}}&
		U(t_n,x,(\yhist^{[n-1]},y))\\
		&\times\pi_n(y\mid x,\yhist^{[n-1]},\beta)
		\,\mu^-(\dd x)\,\dd y,
	\end{aligned}
	\label{eq:phi-from-potential-preposterior}
\end{equation}
because the Bayes normalizing denominator cancels with the predictive
density. More generally, for a belief represented by $\beliefstate$,
one may choose an independent reduced test family
\begin{equation}
	\Phi_n^\eta(t,\mu^\beliefstate;\yhist^{[n]})
	=\phi_n^\eta(t,\beliefstate;\yhist^{[n]}).
	\label{eq:general-reduced-dual-family}
\end{equation}

The full value need not be linear in the belief for a linear lift to
produce an exact penalty. The penalty uses centered observation-update
increments. Components that cancel in this centering do not need to be
represented. The next result illustrates this distinction.

\subsubsection{Scalar LQG with controlled sensing}

Consider
\[
\dd X_t=(aX_t+b\alpha_t)\,\dd t+\sigma\,\dd W_t,
\qquad Y_{t_n}=cX_{t_n}+\beta_{t_n}\varepsilon_n,
\]
with $X_0\sim\mathcal N(m_0,s_0)$ and independent standard normal physical
sensor noises $\varepsilon_n$, independent of the initial state and
Brownian motion. Let $r>0$, $q,q_T\geq0$, and take running cost
$\tfrac12(qX_t^2+r\alpha_t^2)$, terminal cost $\tfrac12q_TX_T^2$, and
sensing cost $\sum_n\kappa_n/\beta_{t_n}$ with deterministic
$\kappa_n\geq0$. The continuous action is unrestricted,
$\alpha_t\in\R$. The sensing-action set $\cB\subset(0,\infty)$ is fixed,
and the initial variance is deterministic.

The Gaussian belief $\mathcal N(m_t,s_t)$ satisfies
\[
\dot m_t=am_t+b\alpha_t,\qquad
\dot s_t=2as_t+\sigma^2
\]
between observations. For the common predictive seed $\xi_n$, its update is
\[
s_n^+=\frac{s_n^-\beta_{t_n}^2}{c^2s_n^-+\beta_{t_n}^2},\qquad
m_n^+=m_n^-+\rho_n\xi_n,\qquad
\rho_n:=\frac{cs_n^-}{\sqrt{c^2s_n^-+\beta_{t_n}^2}},
\]
so $\rho_n^2=s_n^--s_n^+$. Here $m_n^\pm=m_{t_n^\pm}$ and
$s_n^\pm=s_{t_n^\pm}$. In particular, the variance is a deterministic
function of a fixed sensing schedule, independent of the mean and the
continuous control.

\begin{proposition}[Exact certificate for scalar LQG with controlled sensing]
	\label{prop:exact-certificate-lqg}
	In the preceding model, assume that the hypotheses of
	\cref{thm:DY-certificate-main} hold for the test functional below,
	that the Riccati solution is finite on $[0,T]$, and that every
	deterministic schedule in $\cB^{N_o}$ combined with the
	certainty-equivalent feedback is admissible with finite cost.
	Let
	\[
	-\dot A_t=\frac q2+2aA_t-\frac{2b^2}{r}A_t^2,
	\qquad A_T=\frac{q_T}{2}.
	\]
	At every observation date take $U(t_n,x)=A_{t_n}x^2$ and
	$\Phi_n^U(t_n,\mu)=A_{t_n}\int x^2\,\mu(\dd x)$.
	Then
	\[
	D_Y^{\Phi^U}(0,\mu_0)=V_0(\mu_0).
	\]
	The conclusion holds for any finite number of fixed observation dates
	and any nonnegative sensing-cost coefficients satisfying the stated
	assumptions. It does not require attainment of the optimal sensing schedule.
\end{proposition}

\begin{proof}
	The lifted test value is $A_{t_n}(m^2+s)$. For a fixed forecast and
	sensing action, the posterior variance is deterministic and
	$\E[(m_n^+)^2]=(m_n^-)^2+\rho_n^2$. Thus
	\[
	M_Y^{\Phi^U}(u;\xi)
	=\sum_{n=1}^{N_o}A_{t_n}
	\bigl[(m_n^+)^2-(m_n^-)^2-\rho_n^2\bigr].
	\]
	All quantities here are evaluated along the candidate trajectory.
	The pre-posterior averages retain their original-kernel meaning.
	Between observations, the Riccati equation gives
	\[
	\frac12(qm_t^2+r\alpha_t^2)
	+\frac{\dd}{\dd t}(A_tm_t^2)
	=\frac r2\left(\alpha_t+\frac{2bA_t}{r}m_t\right)^2.
	\]
	Integrating over all prediction intervals and accounting for mean jumps,
	\[
	\begin{aligned}
		&\int_0^T\frac12(qm_t^2+r\alpha_t^2)\,\dd t+\frac{q_T}{2}m_T^2\\
		&\quad=A_0m_0^2
		+\sum_nA_{t_n}\bigl[(m_n^+)^2-(m_n^-)^2\bigr]
		+\int_0^T\frac r2\left(\alpha_t+\frac{2bA_t}{r}m_t\right)^2\dd t.
	\end{aligned}
	\]
	For a deterministic sensing schedule
	$\boldsymbol\beta=(\beta_1,\ldots,\beta_{N_o})$, let
	$s_t^{\boldsymbol\beta}$ denote its variance trajectory and define
	\[
	\begin{aligned}
		H(\boldsymbol\beta):={}&A_0m_0^2
		+\sum_nA_{t_n}\rho_n(\boldsymbol\beta)^2
		+\int_0^T\frac q2s_t^{\boldsymbol\beta}\,\dd t\\
		&+\frac{q_T}{2}s_T^{\boldsymbol\beta}
		+\sum_n\frac{\kappa_n}{\beta_n}.
	\end{aligned}
	\]
	Subtracting the penalty and restoring the scenario notation yields
	\[
	\begin{aligned}
		C(u;\xi)-M_Y^{\Phi^U}(u;\xi)
		={}&H\bigl(\boldsymbol\beta(\xi)\bigr)\\
		&+\int_0^T\frac r2
		\left(\alpha_t(\xi)+\frac{2bA_t}{r}m_t^u(\xi)\right)^2\dd t,
	\end{aligned}
	\]
	where $\boldsymbol\beta(\xi)$ is the schedule selected on that scenario.
	The function $H$ has no explicit seed dependence, although
	$H(\boldsymbol\beta(\xi))$ can depend on the seeds through the schedule.
	Therefore every relaxed policy has penalized cost at least
	$\inf_{\boldsymbol\beta\in\cB^{N_o}}H(\boldsymbol\beta)$.
	
	Conversely, fix any deterministic sensing schedule and use
	$\alpha_t=-2bA_tm_t/r$. This is an original admissible policy and
	makes the residual square vanish. Its expected original cost equals
	$H(\boldsymbol\beta)$ because the penalty has zero mean. Consequently,
	\[
	\inf_{\boldsymbol\beta\in\cB^{N_o}}H(\boldsymbol\beta)
	\leq D_Y^{\Phi^U}(0,\mu_0)
	\leq V_0(\mu_0)
	\leq\inf_{\boldsymbol\beta\in\cB^{N_o}}H(\boldsymbol\beta),
	\]
	proving equality throughout.
\end{proof}

The exact coefficients depend on the Riccati equation, not on the
sensing-cost coefficients. Those costs still affect the optimal sensing
schedule and the value. The posterior variance contributes no centered
increment at a fixed sensing action, which explains why the simple lift
can provide an exact penalty without representing the full value.
The proof uses unrestricted continuous actions. It does not establish
the same certificate for bounded controls or nonquadratic state costs.

For three observation dates, the vector
$(A_{t_1},A_{t_2},A_{t_3})$ is therefore an exact maximizer in a dual
family containing independent quadratic coefficients at those dates.
This supplies a quadratic benchmark for the neural-family dual fit in
Experiment~2. It does not by itself guarantee uniqueness or numerical
recovery of the Riccati coefficients.

\subsection{Optimization of the dual test functional}
\label{subsec:dual-maximisation}

Primal training minimizes the cost of the induced policy. Dual training
maximizes the penalized relaxed value. Good primal performance therefore
does not establish that the associated penalty is tight. We use
$U^\theta$ for the primal potential and $\Phi^\eta$ for the dual test
family. The latter may be a lift of a separately parametrized potential
or a direct function of the belief variables.

\subsubsection{Affine penalties and concavity}

Consider a potential linear in its dual coefficients,
\[
U^\eta=U_{\mathrm{base}}+\sum_{j=1}^{p}\eta_j\psi_j,
\qquad \eta\in\mathcal D_{\mathrm{dual}}\subseteq\R^p,
\]
where the base potential and basis functions are fixed. For the LQ
parametrization, the fixed part may enforce the quadratic terminal cost.
Let $\mathcal D_{\mathrm{dual}}$ be a convex parameter set on which the test families
and dual expectations satisfy the required feasibility and integrability
conditions. Linearity of the lift and centering gives
\begin{equation}
	M_Y^{\Phi^\eta}(u;\xi)
	=M_{\mathrm{base}}(u;\xi)+\eta^\top\mathbf p(u;\xi),
	\label{eq:penalty-affine}
\end{equation}
where $\mathbf p=(p_1,\ldots,p_p)$ is the vector of centered penalty
features. It is distinct from the conditional mean $m$.
Write
\begin{equation}
	\begin{aligned}
		\mathcal R(\eta;\xi)
		&:=\inf_{u\in\mcal U^{Y,\mathrm{rel}}}
		\{C(u;\xi)-M_{\mathrm{base}}(u;\xi)-\eta^\top\mathbf p(u;\xi)\},\\
		D(\eta)&:=D_Y^{\Phi^\eta}(0,\mu_0)
		=\E_\xi[\mathcal R(\eta;\xi)].
	\end{aligned}
	\label{eq:dual-G}
\end{equation}
The relaxed control class is fixed independently of $\eta$.
For each scenario, $\mathcal R$ is an infimum of affine functions of
$\eta$ and is therefore concave. Its expectation is concave as well,
and weak duality gives $D(\eta)\leq V_0(\mu_0)$.
These statements concern the exact dual objective. They do not establish
existence of a maximizing coefficient vector or convergence of a
particular numerical optimizer. When a maximizer exists, denote it by
$\eta^\star\in\operatorname*{argmax}_{\eta\in\mathcal D_{\mathrm{dual}}}D(\eta)$.
Write $\widehat\eta$ for the numerically selected coefficients.

\subsubsection{Updating the dual parameters}

The dual objective $D$ is concave in $\eta$. To maximize it, we
first solve the inner problem for each scenario $\xi$ and then
evaluate the penalty features at the resulting optimizer.

More precisely, let $u^\star(\xi,\eta)$ be an exact inner minimizer.
Keeping this control fixed and changing $\eta$ to $\eta'$ gives
\[
\mathcal R(\eta';\xi)
\leq
\mathcal R(\eta;\xi)
-
(\eta'-\eta)^\top
\mathbf p(u^\star(\xi,\eta);\xi).
\]
Consequently, if the minimizer is selected measurably, its penalty
features are integrable, and $D$ is finite on
$\mathcal D_{\mathrm{dual}}$, then
\begin{equation}
	g(\eta)
	:=
	-\E_\xi[
	\mathbf p(u^\star(\xi,\eta);\xi)
	]
	\in \partial^+D(\eta).
	\label{eq:dual-supergradient}
\end{equation}
Here a supergradient is the concave analogue of a gradient:
it provides an affine upper bound,
\[
D(\eta')
\leq D(\eta)+g(\eta)^\top(\eta'-\eta).
\]
The formula only requires evaluating the penalty features at the
inner optimizer. It does not require differentiating that optimizer.

For a fixed parameter $\eta$, independent scenarios
$\xi^{(1)},\ldots,\xi^{(N)}$ give the unbiased estimate
\[
\widehat g(\eta)
=
-\frac1N\sum_{j=1}^N
\mathbf p(u^\star(\xi^{(j)},\eta);\xi^{(j)}).
\]
This estimate can be used in a supergradient-ascent method.
The guarantee requires exact inner minimizers. When only local
solutions are available, the same formula supplies a heuristic
update direction.
\subsubsection{Dual coefficients and induced policies}

If $\Phi^\eta$ is a lift, $U^{\widehat\eta}$ may also be used to generate
a policy through the reduced Hamiltonians. Dual maximization does not
in general guarantee that this policy has a small primal cost.
The continuous LQG selector uses a spatial-gradient moment, while
sensing uses the discrete Hamiltonian. The dual objective instead uses
centered test-value increments. Some changes to a potential affect these
objects differently, whereas others, such as deterministic additive
constants at observation dates, cancel from the penalty altogether.
Any policy generated from dual coefficients must therefore be checked
for admissibility and evaluated separately. The relevant cost is
$J(\alpha^{U^{\widehat\eta}},\beta^{U^{\widehat\eta}})$.

\subsection{Numerical evaluation and primal--dual gaps}
\label{subsec:duality-gap-main}

\subsubsection{Evaluating the penalty and solving the inner problem}

Gaussian separated states give a finite-dimensional representation by
means and covariances. Gaussian propagation does not, however, imply
that the relaxed minimization is easy or globally solvable.
For linear-Gaussian observation models and potentials polynomial in
the relevant state and observation variables, the posterior terms and
pre-posterior averages can be evaluated from joint Gaussian moments.
More general potentials may require quadrature or sampling inside each
evaluation of the inner objective. This additional numerical integration
has its own approximation error.

The accuracy of penalty evaluation and the accuracy of inner minimization
are separate issues. For a fixed scenario, a feasible candidate control
$\widetilde u$ satisfies
\[
C(\widetilde u;\xi)-M_Y^{\Phi^\eta}(\widetilde u;\xi)
\geq\mathcal R(\eta;\xi).
\]
Thus local minimization supplies an upper estimate of the pathwise
infimum, not a certified lower estimate. Its average need not be below
$V_0(\mu_0)$, even with an exact centered penalty and independent
evaluation scenarios. A valid lower estimate requires an exact global
inner solution or a certified lower bound on the inner infimum, together
with appropriate treatment of integration and statistical errors.

In the unrestricted scalar LQG example, continuous controls can be
eliminated analytically. A numerical sensing search still needs its own
global-accuracy assessment: a finite grid with local refinement does
not by itself certify the infimum over a continuous sensing-action set.
Likewise, a quadratic inner objective requires verification of convexity
on the coefficient domain before a convex-solver guarantee can be used.
Bounded-control and nonquadratic examples evaluated by local multistart
are reported as dual diagnostics unless a certified lower estimate of
the inner infimum is available.

Restricting the test family preserves weak duality when each resulting
penalty is feasible. Restricting the inner control class is different:
it may remove original admissible policies and invalidate the bound for
the unrestricted POSOC problem.

\subsubsection{Fitting, selecting, and evaluating dual coefficients}

The true objective is concave, but an inexact inner solver need not
preserve this structure in the computed objective or ascent directions.
For fixed candidate inner controls, the evaluated penalized cost is
affine in $\eta$. On an unrestricted parameter space it can grow without
bound if its aggregate feature slope is nonzero. The true dual objective
remains bounded above by $V_0(\mu_0)$ because it minimizes over all
relaxed controls. An insufficiently accurate inner solver may fail to
find controls that prevent this artificial growth.

Increasing inner iterations and restarts is a useful diagnostic but is
not a proof of global optimality. Independent final evaluation controls
fitting bias from the scenario sample. It does not remove inner-solver
or penalty-integration bias. If coefficient selection uses a validation
sample, final reporting uses a separate sample. If the final iterate
is used without validation selection, final evaluation is still kept
separate from fitting.

The experiment-specific optimization methods, sample sizes, and
integration orders are given in Section~\ref{sec:numerics} and
Appendix~\ref{app:protocol}.

\subsubsection{The reported primal--dual gap}

For an admissible potential and a dual-feasible test family, the exact gap is
\begin{equation}
	\operatorname{Gap}_Y(U,\Phi)
	:=\cJ(U)-D_Y^\Phi(0,\mu_0)\geq0.
	\label{eq:duality-gap-main}
\end{equation}
It bounds the suboptimality of the induced policy:
\[
0\leq\cJ(U)-V_0(\mu_0)\leq\operatorname{Gap}_Y(U,\Phi).
\]
For numerical reporting, first fit and select $\widehat\eta$, and then
evaluate the selected test functional independently of that fitting.
Write
\begin{equation}
	\widehat J_{\mathrm{pr}}(\theta)
	:=\widehat J(\alpha^{U^\theta},\beta^{U^\theta}),
	\qquad
	\widehat D_{\mathrm{dual}}
	:=\widehat D_Y^{\Phi^{\widehat\eta}}(0,\mu_0).
	\label{eq:training-upper-main}
\end{equation}
The hats denote numerical estimates, not deterministic bounds.
Even with exact global inner minimization, a finite Monte Carlo average
or quadrature rule need not lie below the optimal value. A computed
lower certificate requires one-sided error control for penalty evaluation,
inner minimization and outer integration or sampling. A Monte Carlo
standard error quantifies sampling uncertainty only.
One separately optimized dual estimate can be compared with multiple saved
primal checkpoints because the original optimal value is the same.
The reported diagnostic gap is
\begin{equation}
	\widehat{\operatorname{Gap}}_Y(\theta)
	:=\widehat J_{\mathrm{pr}}(\theta)-\widehat D_{\mathrm{dual}}.
	\label{eq:training-gap-main}
\end{equation}
It need not be nonnegative when estimation and optimization errors are
present. A dual estimate above a numerical HJB reference is a warning
that must be assessed against the errors in both calculations.
Certification for the continuous-time problem additionally requires
control of time-discretization and belief-approximation errors.

\begin{remark}[Other information relaxations]
	The construction above reveals predictive seeds in the separated
	belief-state model. Revealing physical sensor noise or Brownian
	increments gives a different information relaxation. Its state
	representation, continuation values, and penalty must be defined
	consistently with that choice, and the weak-duality conditions must
	be verified for the resulting construction.
\end{remark}

\section{From reduced Hamiltonians to a differentiable training problem}
\label{sec:training-workflow}

This section explains how the continuous reduced formulation becomes the finite-dimensional problem solved in the numerical experiments. The goal is to make explicit what information flows through the simulator and what is differentiated during training.

\subsection{Why optimize a potential?}

A partially observed policy is a rule that maps the available observation
history to a decision. Directly parametrizing such rules is possible, but it
can obscure the Hamiltonian structure of the control problem. The present
formulation instead parametrizes a potential $\widehat U^\theta$. Through
the reduced Hamiltonians, this potential determines the continuous
prediction rule and, when sensing is controlled, the observation channel.

At time $\tau_k$, let $\beliefstate_k$ denote the reduced belief state used by the
numerical simulator. In the scalar experiments, $\beliefstate_k=(m_k,\Sigma_k)$.
The realized observation history, or its chosen finite-dimensional feature,
is kept as a separate argument of the potential.
The current potential slice $\widehat U_k^\theta$ and its required spatial
derivatives are inserted into the reduced continuous Hamiltonian. This
determines the continuous control, the prediction of the conditional mean
and covariance, and the running-cost contribution. At an observation time,
the post-observation potential also enters the discrete pre-posterior
Hamiltonian and determines the sensing action when sensing is controlled.
The realized observation is then assimilated through the corresponding
Bayesian or Kalman update.

Thus the trainable variable is $\theta$, while the computational object
induced by $\theta$ is a complete closed-loop Gaussian belief simulator. The
training objective is the expected cost returned by this simulator.

\subsection{Data flow through one simulated trajectory}

For a fixed parameter vector $\theta$, a single simulated trajectory has the following structure:
\begin{equation}
  \beliefstate_0
  \longmapsto
  \beliefstate_1
  \longmapsto
  \cdots
  \longmapsto
  \beliefstate_K,
  \label{eq:simulator-flow}
\end{equation}
where each transition depends on $\theta$ through the potential. More precisely, the transition over one time step can be written as
\begin{equation}
  \beliefstate_{k+1}=\widehat\Phi_k^\theta(\beliefstate_k,\chi_{k+1}),
  \label{eq:trajectory-transition-workflow}
\end{equation}
where $\chi_{k+1}$ contains the random variables needed for state propagation, observation generation, or numerical quadrature. The stage loss is denoted by $\ell_k(\beliefstate_k,\widehat U_k^\theta,\chi_{k+1})$, and the terminal loss by $\mathfrak g(\beliefstate_K)$. The pathwise loss is
\begin{equation}
  \mathfrak L_\Delta(\theta;\xi)
  =\sum_{k=0}^{K-1}\ell_k(\beliefstate_k,\widehat U_k^\theta,\chi_{k+1})+\mathfrak g(\beliefstate_K).
  \label{eq:pathwise-loss-workflow}
\end{equation}
The training objective is the expectation
\begin{equation}
  \cJ_\Delta(\theta)=\E[\mathfrak L_\Delta(\theta;\xi)].
  \label{eq:expected-loss-workflow}
\end{equation}
In practice this expectation is replaced by a mini-batch average, so the plotted minibatch objective is a Monte Carlo estimate on a fresh random batch at each iteration. Its fluctuations around the trained cost are sampling noise and not an instability of the optimizer.

The important point is that the map $\theta\mapsto \mathfrak L_\Delta(\theta;\chi)$ is a deterministic computational graph once the random seed or the random variables $\xi$ are fixed. This is why reverse-mode automatic differentiation can be used. When all components are smooth, the gradient returned by automatic differentiation coincides with the discrete adjoint formula stated and proved in \cref{prop:theta-reduced-discrete-adjoint}. When the reduced Hamiltonian contains a nonsmooth operation such as a sign or an absolute value, the computational graph either has undefined derivatives at the switching set or nearly useless derivatives away from it. The smoothing section addresses precisely this issue.

\subsection{Training objective in reduced notation}

Let $k_n$ denote the time-grid index corresponding to the observation time $t_n$. In the notation of the numerical simulator, the objective has the form
\begin{equation}
  \cJ_\Delta(\theta)
  =\E\bigg[
  \sum_{k=0}^{K-1}\mathfrak f_k(\beliefstate_k,\widehat U_k^\theta)\Delta t_k
  +\sum_{n=1}^{N_o}\mathfrak c_n(\beliefstate_{k_n^-},\widehat U_{k_n}^\theta)
  +\mathfrak g(\beliefstate_K)
  \bigg].
  \label{eq:theta-loss-expanded}
\end{equation}
The terms in this expression have distinct roles.
\begin{itemize}[leftmargin=2em]
  \item $\mathfrak f_k$ is the running-cost slice obtained from the reduced continuous Hamiltonian.
  \item $\mathfrak c_n$ is the discrete observation or sensing-cost slice. It is zero in the fixed-noise regime A examples and nonzero in controlled-sensing regime B.
  \item $\mathfrak g$ is the terminal cost evaluated under the terminal belief.
  \item $\beliefstate_{k_n^-}$ denotes the forecast reduced belief state immediately before the $n$-th observation update.
\end{itemize}
The use of forecast states in the discrete term is important: the sensing decision must be determined before the next observation is known. After the observation is generated or received, the posterior state is used for all subsequent continuous prediction steps.

\begin{remark}[Training and final evaluation]
	In the smooth LQG examples, the same reduced Hamiltonian is used for training
	and evaluation. In the nonsmooth bounded-control problem, the potential is
	trained with a smoothed Hamiltonian but is also evaluated under the exact
	hard switching rule. We therefore report both the smoothed-deployment cost
	and the exact-projection cost, since smoothing is only a numerical device and does not change the original control problem. In Experiment~5 the primal bound is the cost of the potential trained with the smoothed Hamiltonian and deployed with the exact hard switch.
\end{remark}
\subsection{Algorithmic summary}
The following pseudocode summarizes the training loop used by the five numerical experiments. The individual formulas inserted into the Hamiltonian are model-dependent,
but the organization of prediction, update, cost accumulation, and differentiation is
the same.

\begin{algorithm}[h]
\small
\caption{Hamiltonian potential optimization}
\label{alg:training}
\SetKwInOut{Input}{Inputs}\SetKwInOut{Output}{Output}
\Input{Potential family $\widehat U^\theta$, initial parameter $\theta_0$, initial belief,\newline
 time grid and observation schedule, batch size, optimizer, optional smoothing parameter $\delta$.}
\Output{Trained parameter $\theta$, checkpoint policy-cost estimates, and one numerical dual estimate.}
\BlankLine
$\theta\leftarrow \theta_0$\;
\For{optimization iterations}{
  draw a mini-batch of process and observation noise scenarios\;
  \For{each scenario}{
    initialize the numerical belief representation\;
    \For{$k=0,\ldots,K-1$}{
      evaluate $\widehat U_k^\theta$ and the required spatial derivatives\;
      insert these quantities into the reduced continuous Hamiltonian\;
      propagate the belief over $[\tau_k,\tau_{k+1}]$ by the induced prediction map\;
      accumulate the continuous stage-cost contribution\;
      \If{$\tau_{k+1}$ is an observation time}{
        evaluate the reduced discrete Hamiltonian, when sensing is controlled\;
        generate or read the observation\;
        update the belief by the chosen differentiable update map\;
        accumulate the discrete sensing or observation-cost contribution\;
      }
    }
    add the terminal cost\;
  }
  compute the mini-batch estimate of $\cJ_\Delta(\theta)$\;
  compute a parameter update by automatic differentiation or the discrete adjoint, with exact gradients only for differentiable numerical maps\;
  update $\theta$ with the chosen optimizer\;
  \If{the iteration is a checkpoint}{
    estimate the policy cost of the current potential on fresh paths, using a separate random generator\;
  }
}
maximize the dual test functional over $\eta$ on a fit sample, select it on a validation sample, and evaluate it on a fresh sample\;
\Return $\theta$, the policy-cost estimates, and the numerical dual estimate\;
\end{algorithm}
The algorithm is deliberately written at a level that covers all five numerical experiments. The concrete formulas needed by the scalar Gaussian experiments are collected in Appendix~\ref{sec:gaussian-formulas}, while the experiment-specific objectives and parameter values are given in \cref{sec:numerics}.

\section{Numerical discretization and gradient computation}
\label{sec:numerical-method}

The continuous formulation is implemented on a time grid
\[
0=\tau_0<\tau_1<\cdots<\tau_K=T,
\]
which contains all observation times. We write
$\Delta t_k=\tau_{k+1}-\tau_k$.

The discretization has two types of steps. Between observations, the
conditional mean and covariance are propagated by the prediction map
induced by the continuous reduced Hamiltonian. At an observation time, the
forecast mean and covariance are transformed into posterior quantities by
the Bayesian or Kalman update and, when observation noise is controlled, by
the reduced discrete Hamiltonian.

Keeping these two operations separate is important. The sensing action must
be chosen using the forecast information available before the observation,
whereas all subsequent prediction steps use the posterior belief.

All numerical experiments in this paper use Gaussian beliefs, represented
by their conditional means and covariances. The discrete adjoint is written
using an abstract finite-dimensional state $\beliefstate$ so that the scalar and
multidimensional Gaussian cases can be treated with the same notation.
Particle and ensemble approximations are not implemented in the present
work.
\subsection{Gaussian belief representation}
\label{subsec:gaussian-belief-representation}

The numerical state is the mean--covariance representation
\begin{equation}
	\beliefstate_k=(m_k,\Sigma_k,h_k),
	\qquad
	\mu^{\beliefstate_k}=\mathcal N(m_k,\Sigma_k),
\end{equation}

All state and observation models in the numerical experiments are
linear-Gaussian. Conditional on the full observation history, the control
between observations is a known input, so Gaussian prediction and
conditioning remain exact. The obstacle running cost is nonquadratic,
but it does not alter the conditional-law argument. Time discretization
still introduces numerical error. The Gaussian belief variables remain
$(m_k,\Sigma_k)$. $h_k$ records any retained observation-history features
used by the potential. It is included in the numerical state solely to
close the differentiation of the simulator. The associated probability
measure depends only on $(m_k,\Sigma_k)$. A fixed-size padded history with
deterministic availability masks may be used. Between observations the
history is unchanged. At an observation it is updated with the newly
generated observation. If no history is used, $h_k$ is omitted.
The initial state, including retained initial data, is independent of
$\theta$. The discrete sensing decision uses the pre-observation history.
Subsequent prediction steps use the updated history.

This representation has two advantages for the present study. First, the
prediction and observation updates are finite-dimensional differentiable
maps. Second, expectations of polynomial potential derivatives and of the
cost functions used in the experiments can be evaluated analytically from
Gaussian moments or Gaussian probabilities.

\begin{remark}[Scope of the numerical representation]
	The Hamiltonian-potential formulation itself is not restricted to Gaussian
	beliefs. A particle or ensemble representation could be introduced by
	replacing $(m_k,\Sigma_k)$ with an appropriate finite-dimensional numerical
	state and by defining corresponding prediction and observation-update maps.
	Such an extension raises additional questions concerning resampling,
	quadrature, and gradient estimation. These issues are outside the scope of
	the present paper.
\end{remark}

\subsection{Discrete adjoint identity for parametrized potentials}
\label{subsec:dual-gradients}

We now rewrite the discrete adjoint construction entirely at the level of the
\emph{reduced closed-loop maps induced by the parametrized potential}
$\widehat U^\theta$, with $\theta\in\Theta\subset\mathbb R^p$.
This is the discrete counterpart of Section~\ref{sec:method}: once $\widehat U^\theta$ is fixed,
it induces a closed-loop continuous propagation mechanism between observations
and a closed-loop Bayesian update mechanism at the observation times.
The optimization variable is therefore $\theta$ alone.
No independent control variable appears in the formulation below.

Let
\[
0=\tau_0<\tau_1<\dots<\tau_K=T,
\qquad
\Delta t_k:=\tau_{k+1}-\tau_k,
\qquad k=0,\dots,K-1,
\]
be a numerical time grid, and let
\[
k_{\bullet}:\{1,\dots,N_o\}\to\{1,\dots,K\}
\]
denote the observation-index map, so that the $n$-th observation time satisfies
$\tau_{k_n}=t_n$.
Here the observation schedule is that of the implemented discrete problem.
When nominal dates are rounded to the simulation grid, the implemented
dates are distinguished from the nominal ones, as recorded for
Experiment~2 in Appendix~\ref{app:protocol}.
We also write
\[
\nu(\tau_k):=\max\{j\in\{0,\dots,N_o\}:t_j\le \tau_k\}
\]
for the index of the most recent observation available at time $\tau_k$.

\paragraph{Reduced closed-loop maps induced by $\widehat U^\theta$.}
For each time step $k=0,\dots,K-1$, let
\[
\widehat\Phi_k^\theta(\beliefstate,\widehat U_k^\theta,\chi_{k+1}^c;\Delta t_k)
\]
denote the numerical scheme map on $[\tau_k,\tau_{k+1}]$ obtained after inserting the
candidate potential $\widehat U^\theta$ into the continuous Hamiltonian minimisation
and then applying the chosen numerical propagation scheme.
Likewise, for each observation time $t_n$, let
\[
\widehat\Psi_n^\theta(\beliefstate^-,\widehat U_{k_n}^\theta,y,\chi_n^d)
\]
denote the numerical update map obtained after inserting $\widehat U^\theta$ into
the discrete Hamiltonian minimisation and then applying the chosen normalized
observation update scheme.
Thus we keep the dependence on the potential explicit in the argument list:
the maps are evaluated at the current reduced state \emph{and} at the current
potential slice $\widehat U_k^\theta$ or $\widehat U_{k_n}^\theta$. The superscript $\theta$ then records that these slices themselves come from the parametrized family $\widehat U^\theta$.
The state arguments $\beliefstate$ and $\beliefstate^-$ are dummy inputs of the maps, whereas
$\widehat\beliefstate_k^\theta$ and $\widehat\beliefstate_{k_n^-}^\theta$ below are the actual
trajectory values produced by those maps.

Starting from $\widehat\beliefstate_0^\theta=\beliefstate_0$, we define the reduced closed-loop recursion by
\begin{equation}
	\widehat\beliefstate_{k+1^-}^\theta
	=
	\widehat\Phi_k^\theta\bigl(\widehat\beliefstate_k^\theta,\widehat U_k^\theta,\chi_{k+1}^c;\Delta t_k\bigr),
	\qquad k=0,\dots,K-1,
	\label{eq:theta-reduced-scheme}
\end{equation}

\label{eq:theta-forward-phi-hamiltonian}
and, at the observation times,
\begin{equation}
	\widehat\beliefstate_{k_n}^\theta
	=
	\widehat\Psi_n^\theta\bigl(\widehat\beliefstate_{k_n^-}^\theta,\widehat U_{k_n}^\theta,Y_{t_n}^\theta,\chi_n^d\bigr),
	\qquad n=1,\dots,N_o.
	\label{eq:theta-reduced-update}
\end{equation}

\label{eq:theta-forward-psi-hamiltonian}
At non-observation steps we set
\begin{equation}
	\widehat\beliefstate_{k+1}^\theta=\widehat\beliefstate_{k+1^-}^\theta
	\qquad\text{whenever }k+1\neq k_n\text{ for all }n.
	\label{eq:theta-no-observation}
\end{equation}

\label{eq:theta-forward-no-update-hamiltonian}

The role of \eqref{eq:theta-reduced-scheme} is to mimic, at the numerical level,
the continuous belief dynamics under $\widehat U$ between observation times,
while \eqref{eq:theta-reduced-update} mimics the normalized discrete update
at the observation times.
The important point is that these are \emph{full numerical maps}.
All dependence on $\Delta t_k$ is carried by the actual scheme.

\label{eq:theta-particle-em-weights}

\paragraph{Reduced stage costs.}
We now absorb the closed-loop dependence on $\widehat U^\theta$ into reduced cost
slices.
For $k=0,\dots,K-1$, let
\begin{equation}
	\mathfrak{f}_k^{\Delta}(\beliefstate,\widehat U_k^\theta; \Delta t_k ) = \mathfrak{f}_k(\beliefstate,\widehat U_k^\theta) \Delta t_k  
	\label{eq:theta-reduced-running-cost}
\end{equation}

be the reduced running-cost density associated with the
closed-loop dynamics induced by $\widehat U^\theta$.
At the observation times, let
\begin{equation}
	\mathfrak c_n^\theta(\beliefstate^-,\widehat U_{k_n}^\theta),
	\label{eq:theta-reduced-observation-cost}
\end{equation}

be the reduced observation cost induced by
$\widehat U^\theta$.
Finally, let
\begin{equation}
	\mathfrak g(\beliefstate)
	\label{eq:theta-reduced-terminal-cost}
\end{equation}

be the terminal cost.

With this notation, the discretized objective reads
\begin{equation}
	\mcal J_\Delta(\widehat U^\theta)
	:=
	\mathbb E\Bigg[
	\sum_{k=0}^{K-1} \mathfrak f_k(\widehat\beliefstate_k^\theta,\widehat U_k^\theta)\,\Delta t_k
	+
	\sum_{n=1}^{N_o} \mathfrak c_n^\theta(\widehat\beliefstate_{k_n^-}^\theta,\widehat U_{k_n}^\theta)
	+
	\mathfrak g(\widehat\beliefstate_K^\theta)
	\Bigg].
	\label{eq:theta-discrete-objective}
\end{equation}

\paragraph{Parameter-dependent observations.}
All expectations below are taken over base random inputs
$\chi^c,\chi^d$ whose joint law is independent of $\theta$.
These inputs are held fixed during differentiation. At an observation time,
write the predictive simulation of the observation as
\begin{equation}
	Y_{t_n}^\theta
	=G_n^\theta\bigl(\widehat\beliefstate_{k_n^-}^\theta,
	\widehat U_{k_n}^\theta,\chi_n^d\bigr).
	\label{eq:theta-generated-observation}
\end{equation}

Here the potential-induced sensing action is selected from forecast
information, before the current observation is generated. Even with fixed
observation noise, the simulated observation can depend on $\theta$
through the forecast distribution.

Define the composed observation-update map
\begin{equation}
	\overline\Psi_n^\theta(\beliefstate,U,\chi_n^d)
	:=\widehat\Psi_n^\theta
	\bigl(\beliefstate,U,G_n^\theta(\beliefstate,U,\chi_n^d),\chi_n^d\bigr).
	\label{eq:theta-composed-update}
\end{equation}

The map includes the update of any retained observation-history features.
Thus \eqref{eq:theta-reduced-update} is evaluated at $Y_{t_n}^\theta$,
and is equivalently written using $\overline\Psi_n^\theta$.

\paragraph{Derivative notation.}
For a reduced object $F^\theta(\beliefstate,U,y)$,
$\nabla_\beliefstate F^\theta$ and $\partial_y F^\theta$ denote
derivatives in the indicated arguments, with the other arguments fixed.
The symbol $D_\theta$ denotes the direct parameter derivative at fixed
reduced state and, when $y$ is an independent argument, at fixed $y$:
\[
D_\theta F^\theta(\beliefstate,\widehat U^\theta,y)
:=\left.\frac{\dd}{\dd\vartheta}
F^\vartheta(\beliefstate,\widehat U^\vartheta,y)
\right|_{\vartheta=\theta}.
\]
It includes the dependence on the parametrized potential. For a composed
map such as $\overline\Psi_n^\theta$, it also differentiates through the
generated observation at fixed reduced state. In particular,
\begin{align}
	\partial_\beliefstate\overline\Psi_n^\theta
	&=\partial_\beliefstate\widehat\Psi_n^\theta
	+\partial_y\widehat\Psi_n^\theta\,
	\partial_\beliefstate G_n^\theta,
	\label{eq:theta-composed-state-derivative}\\
	D_\theta\overline\Psi_n^\theta
	&=D_\theta\widehat\Psi_n^\theta
	+\partial_y\widehat\Psi_n^\theta\,D_\theta G_n^\theta.
	\label{eq:theta-composed-parameter-derivative}
\end{align}

The right-hand sides are evaluated at
$y=G_n^\theta(\beliefstate,\widehat U_{k_n}^\theta,\chi_n^d)$.
Holding the potential argument fixed means holding its function fixed.
Its evaluation at state or history features is still differentiated with
respect to those features.

For a scalar reduced object without an independent observation argument,
the total chain rule is
\begin{equation}
	\frac{\dd}{\dd\theta}
	F^\theta(\widehat\beliefstate_k^\theta,\widehat U_k^\theta)
	=(Q_k^\theta)^\top\nabla_\beliefstate
	F^\theta(\widehat\beliefstate_k^\theta,\widehat U_k^\theta)
	+D_\theta F^\theta(\widehat\beliefstate_k^\theta,\widehat U_k^\theta),
	\qquad
	Q_k^\theta:=\frac{\dd\widehat\beliefstate_k^\theta}{\dd\theta}.
	\label{eq:theta-master-chain-rule}
\end{equation}

We use column gradients for scalar functions and output-by-input
Jacobians for vector maps.

\begin{proposition}[Reduced discrete adjoint representation in Hamiltonian form]
	\label{prop:theta-reduced-discrete-adjoint}
	\label{prop:theta-reduced-discrete-adjoint-hamiltonian}
	Assume that the reduced prediction and observation-update maps are the
	chosen numerical maps induced by the continuous and discrete augmented
	Hamiltonians at $\gamma=1$. Define the composed discrete Hamiltonian by
	\begin{equation}
		\overline{\mcal H}_n^{d,\theta}
		(\beliefstate,U,\gamma;\chi_n^d)
		:=\widehat{\mcal H}_n^{d,\theta}
		\bigl(\beliefstate,U,\gamma;
		G_n^\theta(\beliefstate,U,\chi_n^d),\chi_n^d\bigr).
		\label{eq:theta-composed-hamiltonian}
	\end{equation}
	
	The observation generator in this definition is the one selected at
	$\gamma=1$ and is held fixed when taking $\partial_\gamma$.
	Suppose the reduced cost slices are identified with
	\begin{align}
		\mathfrak f_k(\beliefstate,U)
		&=\partial_\gamma\widehat{\mcal H}_k^{c,\theta}
		(\beliefstate,U,1;\chi_{k+1}^c),
		\label{eq:theta-running-hamiltonian-slice}\\
		\mathfrak c_n^\theta(\beliefstate,U)
		&=\partial_\gamma\overline{\mcal H}_n^{d,\theta}
		(\beliefstate,U,1;\chi_n^d).
		\label{eq:theta-observation-hamiltonian-slice}
	\end{align}
	
	The fixed seed arguments in the cost-slice notation are suppressed.
	If the observation cost is independent of the realized observation,
	the second identity simply has no observation-argument contribution.
	Assume that these slices, the composed maps, and the terminal cost are
	pathwise differentiable in the relevant arguments, and that differentiation
	can be exchanged with expectation.
	
	Starting from the fixed initial state, use
	\eqref{eq:theta-reduced-scheme}--\eqref{eq:theta-no-observation}
	with $Y_{t_n}=Y_{t_n}^\theta$. Then
	\begin{align}
		\mcal J_\Delta(\widehat U^\theta)
		={}&\mathbb E\Bigg[
		\sum_{k=0}^{K-1}\partial_\gamma\widehat{\mcal H}_k^{c,\theta}
		(\widehat\beliefstate_k^\theta,\widehat U_k^\theta,1;\chi_{k+1}^c)
		\Delta t_k+\sum_{n=1}^{N_o}\partial_\gamma\overline{\mcal H}_n^{d,\theta}
		(\widehat\beliefstate_{k_n^-}^\theta,\widehat U_{k_n}^\theta,1;\chi_n^d)
		+\mathfrak g(\widehat\beliefstate_K^\theta)\Bigg].
		\label{eq:theta-discrete-objective-hamiltonian}
	\end{align}
	
	Define the terminal adjoint by
	\begin{equation}
		p_K^\theta:=\nabla_\beliefstate\mathfrak g
		(\widehat\beliefstate_K^\theta).
		\label{eq:theta-terminal-adjoint}
	\end{equation}
	\label{eq:theta-terminal-adjoint-hamiltonian}
	At each observation time, define the pre-update adjoint by
	\begin{align}
		p_{k_n^-}^\theta:={}&
		\nabla_\beliefstate\partial_\gamma\overline{\mcal H}_n^{d,\theta}
		(\widehat\beliefstate_{k_n^-}^\theta,\widehat U_{k_n}^\theta,1;\chi_n^d)
		+\Bigl(\partial_\beliefstate\overline\Psi_n^\theta
		(\widehat\beliefstate_{k_n^-}^\theta,\widehat U_{k_n}^\theta,\chi_n^d)
		\Bigr)^\top p_{k_n}^\theta.
		\label{eq:theta-preupdate-adjoint}
	\end{align}
	\label{eq:theta-preupdate-adjoint-hamiltonian}
	At non-observation indices set
	\begin{equation}
		p_{k+1^-}^\theta:=p_{k+1}^\theta
		\qquad\text{whenever }k+1\ne k_n\text{ for all }n.
		\label{eq:theta-no-update-adjoint}
	\end{equation}
	\label{eq:theta-no-update-adjoint-hamiltonian}
	For $k=K-1,\ldots,0$, define
	\begin{align}
		p_k^\theta:={}&
		\nabla_\beliefstate\partial_\gamma\widehat{\mcal H}_k^{c,\theta}
		(\widehat\beliefstate_k^\theta,\widehat U_k^\theta,1;\chi_{k+1}^c)
		\Delta t_k+\Bigl(\partial_\beliefstate\widehat\Phi_k^\theta
		(\widehat\beliefstate_k^\theta,\widehat U_k^\theta,
		\chi_{k+1}^c;\Delta t_k)\Bigr)^\top p_{k+1^-}^\theta.
		\label{eq:theta-backward-adjoint}
	\end{align}
	\label{eq:theta-backward-adjoint-hamiltonian}
	Then
	\begin{align}
		\nabla_\theta\mcal J_\Delta(\widehat U^\theta)
		={}&\mathbb E\Bigg[
		\sum_{k=0}^{K-1}\Bigl(
		D_\theta\partial_\gamma\widehat{\mcal H}_k^{c,\theta}
		(\widehat\beliefstate_k^\theta,\widehat U_k^\theta,1;\chi_{k+1}^c)
		\Delta t_k
		+\bigl(D_\theta\widehat\Phi_k^\theta
		(\widehat\beliefstate_k^\theta,\widehat U_k^\theta,
		\chi_{k+1}^c;\Delta t_k)\bigr)^\top p_{k+1^-}^\theta\Bigr)
		\nonumber\\
		&+\sum_{n=1}^{N_o}\Bigl(
		D_\theta\partial_\gamma\overline{\mcal H}_n^{d,\theta}
		(\widehat\beliefstate_{k_n^-}^\theta,\widehat U_{k_n}^\theta,1;\chi_n^d)
		+\bigl(D_\theta\overline\Psi_n^\theta
		(\widehat\beliefstate_{k_n^-}^\theta,\widehat U_{k_n}^\theta,
		\chi_n^d)\bigr)^\top p_{k_n}^\theta\Bigr)
		+D_\theta\mathfrak g(\widehat\beliefstate_K^\theta)\Bigg].
		\label{eq:theta-dual-gradient}
	\end{align}
	\label{eq:theta-dual-gradient-hamiltonian}
\end{proposition}
The proof is given in \cref{app:-first-adjoint-proof}.

The representation \eqref{eq:theta-dual-gradient} expresses the gradient
entirely through the numerical maps and cost slices induced by
$\widehat U^\theta$, without introducing independent continuous or
discrete control variables. The composed observation-update map includes
the dependence of the generated observations on $\theta$.
The backward-gradient algorithm therefore provides a Monte Carlo
estimator of $\nabla_\theta\cJ_\Delta(\widehat U^\theta)$ by evaluating
this representation along simulated trajectories.

This interpretation requires the complete numerical simulator, including
observation generation and any retained history, to be pathwise
differentiable, with sufficient integrability to exchange differentiation
and expectation. For the Euler prediction and Kalman update maps used
in the LQG experiments, differentiability holds when the innovation
covariance remains nonsingular and the potential-induced coefficients
are differentiable. The same requirement applies to the Gaussian
prediction, observation-generation, and update maps used in the obstacle
experiment.

A discontinuous minimizing selector can violate these assumptions even
when the belief prediction and observation update are otherwise smooth.
This occurs in the bounded-control example, where changing the potential
moves the boundary between minimizing Hamiltonian branches.
We therefore introduce smooth Hamiltonians whose derivatives define
differentiable admissible feedback rules and the corresponding cost
slices. The discrete adjoint then computes the gradient of the resulting
smoothed discrete objective. In the affine-control example below, these
cost slices coincide with the original running cost evaluated under the
smoothed feedback.

\subsection{Smoothing nonsmooth reduced Hamiltonians}
\label{subsec:smoothing}

The bounded-control example in \cref{subsec:nonsmooth} isolates the
difficulty caused by branch switching while retaining an exact Gaussian
belief representation, $\beliefstate_k=(m_k,\Sigma_k)$.
We first identify the nonsmooth term in the reduced Hamiltonian and
then construct smooth feedback rules for training.

\subsubsection{Reduced bounded-control Hamiltonian.}

Consider the scalar controlled generator and running cost
\begin{align}
	(\cG^\alpha U)(x)
	&=(ax+b\alpha)\partial_xU(x)
	+\frac{\sigma^2}{2}\partial_{xx}U(x),
	\qquad |\alpha|\le\alpha_{\max},
	\label{eq:ns-generator}\\
	f(x,\alpha)
	&=\frac q2x^2+\kappa\alpha.
	\label{eq:ns-running-cost}
\end{align}
For $\mu_k=\mcal N(m_k,\Sigma_k)$, the augmented reduced Hamiltonian is
\begin{equation}
	\tH^c_k(\mu_k,U,\gamma)
	:=
	\inf_{|\alpha|\le\alpha_{\max}}
	\left\langle
	\cG^\alpha U+\gamma f(\cdot,\alpha),\mu_k
	\right\rangle.
	\label{eq:ns-augmented-hamiltonian-before-reduction}
\end{equation}
The variable $\gamma$ is a running-cost multiplier. The original
Hamiltonian is recovered at $\gamma=1$. Where the minimizing branch
is unique, differentiation with respect to $\gamma$ extracts its
running-cost contribution.

For the quadratic potential
\begin{equation}
	\widehat U_k^\theta(x)
	=\frac12P_kx^2+r_kx+c_k^U,
	\label{eq:ns-quadratic-potential}
\end{equation}
define the gradient moment and switching score by
\begin{equation}
	\eta_k^\theta
	:=\left\langle
	\partial_x\widehat U_k^\theta,\mu_k
	\right\rangle
	=P_km_k+r_k,
	\qquad
	\zeta_k^\gamma:=b\eta_k^\theta+\gamma\kappa.
	\label{eq:ns-gradient-moment}
\end{equation}
Since
\begin{equation}
	\inf_{|\alpha|\le\alpha_{\max}}\alpha\zeta_k^\gamma
	=-\alpha_{\max}|\zeta_k^\gamma|,
	\label{eq:ns-kink-explicit}
\end{equation}
the Hamiltonian becomes
\begin{equation}
	\begin{aligned}
		\tH^c_k(m_k,\Sigma_k,P_k,r_k,\gamma)
		={}&
		a\bigl(P_k(m_k^2+\Sigma_k)+r_km_k\bigr)
		+\frac{\sigma^2}{2}P_k\\
		&+\gamma\frac q2(m_k^2+\Sigma_k)
		-\alpha_{\max}|\zeta_k^\gamma|.
	\end{aligned}
	\label{eq:ns-exact-ham}
\end{equation}
The constant $c_k^U$ does not enter this continuous Hamiltonian.

At $\gamma=1$, a minimizing selector is
$-\alpha_{\max}\operatorname{sign}(\zeta_k^1)$.
We set $\operatorname{sign}(0)=0$. When the switching score vanishes,
every admissible action is a minimizer. The resulting Euler
belief prediction is
\begin{align}
	m_{k+1}^{-}
	&=m_k+
	\left(
	am_k-b\alpha_{\max}
	\operatorname{sign}(\zeta_k^1)
	\right)\Delta t,
	\label{eq:ns-exact-mean-prediction}\\
	\Sigma_{k+1}^{-}
	&=\Sigma_k+(2a\Sigma_k+\sigma^2)\Delta t.
	\label{eq:ns-variance-prediction}
\end{align}

\subsubsection{Branch switching and pathwise differentiation.}

The action-dependent Hamiltonian is the minimum of two endpoint
branches, corresponding to $\alpha=\pm\alpha_{\max}$. These branches
have equal values on the tie set
\begin{equation}
	\zeta_k^\gamma=b(P_km_k+r_k)+\gamma\kappa=0.
	\label{eq:ns-switching-hypersurface-gamma}
\end{equation}
At $\gamma=1$, this is the switching surface
\begin{equation}
	b(P_km_k+r_k)+\kappa=0.
	\label{eq:ns-switching-hypersurface}
\end{equation}
The reduced Hamiltonian is continuous but generally nonsmooth
on this set, while its minimizing selector switches between
the two endpoints.

Away from the tie set, the hard selector is locally constant,
so its derivatives with respect to the potential coefficients
vanish. Direct differentiation through
\eqref{eq:ns-exact-mean-prediction} therefore misses the contribution
from shifting the switching boundary, although this shift can change
the expected cost. Even when a sampled trajectory almost surely
avoids the boundary, differentiating individual trajectories need
not recover the derivative of the expected objective.
Hamiltonian smoothing replaces this discontinuous branch selection
by a differentiable transition.

\subsubsection{Smooth Hamiltonians and induced controls.}

Since the reduced Hamiltonian is available explicitly, we smooth
its absolute-value term directly. We replace
$\alpha_{\max}|\zeta|$ by either
\begin{align}
	\psi_\delta^{\mathrm{soft}}(\zeta)
	&=\delta\log\left(
	2\cosh\left(\frac{\alpha_{\max}\zeta}{\delta}\right)
	\right),
	\label{eq:ns-soft-smoothing}\\
	\psi_\delta^{\mathrm{sqrt}}(\zeta)
	&=\alpha_{\max}\sqrt{\zeta^2+\delta^2},
	\qquad \delta>0.
	\label{eq:ns-sqrt-smoothing}
\end{align}
The first choice gives the log-sum-exp soft-min of the two
endpoint branches. The second is an analytic square-root
regularization. For either choice,
\begin{equation}
	\begin{aligned}
		\tH^{c,\psi}_{\delta,k}
		(m_k,\Sigma_k,P_k,r_k,\gamma)
		={}&
		a\bigl(P_k(m_k^2+\Sigma_k)+r_km_k\bigr)
		+\frac{\sigma^2}{2}P_k\\
		&+\gamma\frac q2(m_k^2+\Sigma_k)
		-\psi_\delta(\zeta_k^\gamma).
	\end{aligned}
	\label{eq:ns-smoothed-ham}
\end{equation}
At $\gamma=1$, differentiation of the smoothed action term
defines the effective smoothed control
$\alpha_k^{\theta,\delta}=-\psi_\delta'(\zeta_k^1)$.
Explicitly,
\begin{align}
	\alpha_k^{\theta,\mathrm{soft}}
	&=-\alpha_{\max}
	\tanh\left(
	\frac{\alpha_{\max}\zeta_k^1}{\delta}
	\right),
	\label{eq:ns-soft-control}\\
	\alpha_k^{\theta,\mathrm{sqrt}}
	&=-\alpha_{\max}
	\frac{\zeta_k^1}
	{\sqrt{(\zeta_k^1)^2+\delta^2}}.
	\label{eq:ns-sqrt-control}
\end{align}
The associated mean prediction is
\begin{equation}
	m_{k+1}^{-}
	=m_k+\left(
	am_k+b\alpha_k^{\theta,\delta}
	\right)\Delta t,
	\label{eq:ns-smoothed-mean-prediction}
\end{equation}
while the variance prediction remains
\eqref{eq:ns-variance-prediction}. The running-cost slice satisfies
\begin{equation}
	\left.
	\partial_\gamma\tH^{c,\psi}_{\delta,k}
	\right|_{\gamma=1}
	=
	\frac q2(m_k^2+\Sigma_k)
	+\kappa\alpha_k^{\theta,\delta}.
	\label{eq:ns-smoothed-cost-slice}
\end{equation}
Thus, in this affine-control example, the cost slice evaluates
the original running cost under the effective smoothed control.
Using this slice together with
\eqref{eq:ns-smoothed-mean-prediction} defines a consistent differentiable simulator.
Under the stated regularity and integrability assumptions,
Proposition~\ref{prop:theta-reduced-discrete-adjoint} applies.

The parameter $\delta$ controls both the Hamiltonian approximation
and the conditioning of the induced simulator. Larger values broaden
the transition between active branches. As $\delta$ decreases,
the effective smoothed control approaches the hard selector away
from the tie set, but its derivatives become concentrated near
that set, potentially making gradient-based training difficult.
Convergence of the Hamiltonian values alone does not imply
convergence of the corresponding objective gradients.

We train by minimizing the discretized objective
$\cJ_{\Delta,\delta}(\theta)$ using the effective smoothed control
and its associated belief dynamics. After training, we evaluate
both this control and the exact minimizing selector reconstructed
from the learned potential, always using the original control cost.
Both controls remain within
$[-\alpha_{\max},\alpha_{\max}]$ and respect the observation
information structure. Their expected costs therefore provide
primal upper bounds for the same discretized control problem.
Evaluation with the exact selector additionally assesses the
switching boundary learned during training.

\section{Numerical experiments}
\label{sec:numerics}

\paragraph{Potential parametrization and implementation details.}
The potential parametrization should provide sufficient expressiveness, stable spatial derivatives, and inexpensive expectations under the belief. The experiments use polynomial potentials, a structured quadratic potential with neural time coefficients, and a history-dependent polynomial--trigonometric family. Observation features supplied to the potential are kept distinct from the posterior used in the Hamiltonian expectations: compressing the potential's history input does not replace the full belief. Detailed potential families and basis definitions are given in Appendix~\ref{sec:parametrization}. The scalar Gaussian moment, prediction, observation-update, and cost formulas are collected with their derivations in Appendix~\ref{sec:gaussian-formulas}. Experimental settings are given in Appendix~\ref{app:protocol}.

This section presents five numerical experiments of increasing complexity.
Experiment~1 is the LQG problem with fixed observation noise (regime~A), in
dimensions $d=1$ to $d=100$. Experiment~2 introduces controlled sensing
(regime~B) with three observations. Experiment~3 adds bounds on the control, so
that the optimal sensing can depend on the observation history. Experiment~4 is a
nonlinear obstacle problem, and Experiment~5 is a bounded-control LQG problem that
requires Hamiltonian smoothing.

\paragraph{Common protocol.}
In all experiments the plotted curves have the same meaning. The potential is
trained by a first-order method, and the training objective is shown as a grey
curve. At selected iterations the current potential is used as an admissible
policy, giving its numerical cost estimate (orange). The exact admissible-policy
cost is an upper bound on the optimal value. A Monte Carlo estimate is not a
deterministic upper bound. The green horizontal line is a separately optimized
numerical dual estimate. Weak duality applies to the exact dual functional for
every feasible test functional, whereas its computed estimate also depends on
inner minimization, penalty evaluation and outer integration or sampling.
Exact inner minimization alone does not certify that estimate. Certification
requires appropriate one-sided control of these errors, as discussed in
Section~\ref{subsec:duality-gap-main}. Monte Carlo standard errors describe
sampling uncertainty only. Quadrature and solver errors require separate
assessment. The plotted primal--dual gaps are numerical diagnostics.
In Experiments~1, 4 and~5 the potentials are trained on random Monte Carlo
minibatches, so the grey curve fluctuates around the trained cost because of the
finite batch. In Experiment~2 the structured potential is trained with exact Gaussian
expectations over the observation innovations, so the training curve is not noisy and
the reported primal costs have no Monte Carlo error. The dual estimate is evaluated separately on fresh Monte Carlo scenarios. Experiment~3 is also trained with a tensor Gauss--Hermite rule over the observation innovations, so its training curve is smooth. Its policy costs are Monte Carlo estimates. In Experiment~3 the dual estimate is
obtained once, by maximizing the dual from a random start without validation selection
(Section~\ref{subsec:dual-maximisation}). It is a numerical estimate and is not evaluated along the training.
 Numerical settings are collected in Appendix~\ref{app:protocol}.

\subsection{Experiment 1: LQG regime A with fixed observation noise and dimensional scaling}
\label{subsec:lqg-regime-a}
We consider the LQG under regime A, in which the observation times and observation-noise level are fixed. Its purpose is twofold. We first use the scalar problem to verify the basic mechanism of the method in a setting where the belief recursion is exact, low-dimensional, and analytically understood. We then consider  coupled multidimensional problems to investigate the approximation quality and computational scaling of the same method.

The hidden state satisfies
\begin{equation}
	\dd X_t
	=
	(aX_t+b\alpha_t)\,\dd t
	+\sigma\,\dd W_t,
	\qquad
	X_0\sim\mcal N(m_0,\Sigma_0),
	\label{eq:num-regA-state}
\end{equation}
and the observations are
\begin{equation}
	Y_{t_n}
	=
	cX_{t_n}+\varepsilon\xi_n,
	\qquad
	\xi_n\stackrel{\mathrm{iid}}{\sim}\mcal N(0,1),
	\qquad
	n=1,\ldots,N_o.
	\label{eq:num-regA-observation}
\end{equation}
The observation-noise level $\varepsilon$ is fixed and  the objective is
\begin{equation}
	J(\alpha)
	=
	\E\bigg[
	\int_0^T
	\left(
	\frac12 qX_t^2
	+\frac12 r\alpha_t^2
	\right)\dd t
	+\frac12 q_TX_T^2
	\bigg].
	\label{eq:num-regA-cost}
\end{equation}
Since the model is linear Gaussian, the belief remains Gaussian and the reduced state is $\beliefstate_k=(m_k,\Sigma_k)$. In both the scalar and multidimensional experiments, the potential is quadratic in $x$ and has no explicit observation-history input. In the scalar case, we use
\begin{equation}
	U_k^\theta(x)=c_k+\ell_k x+\frac12 h_kx^2,
	\qquad
	\eta_k^\theta
	=\int_{\mathbb R}\partial_xU_k^\theta(x)\,\mu^{\beliefstate_k}(\dd x)
	=\ell_k+h_km_k.
\end{equation}
The coefficients are parametrized separately at each grid time in the scalar experiment. Although the potential has no explicit observation-history input, the induced control uses all past observations through the Kalman mean $m_k$. The resulting discrete training loss is
\begin{equation}
	\cJ_\Delta(\theta)
	=
	\E\bigg[
	\sum_{k=0}^{K-1}
	\left(
	\frac12 q(m_k^2+\Sigma_k)
	+\frac{b^2}{2r}(\eta_k^\theta)^2
	\right)\Delta t_k
	+\frac12 q_T(m_K^2+\Sigma_K)
	\bigg].
	\label{eq:num-regA-loss}
\end{equation}

\begin{figure}[h]
	\centering
	\includegraphics[width=0.6\textwidth]
	{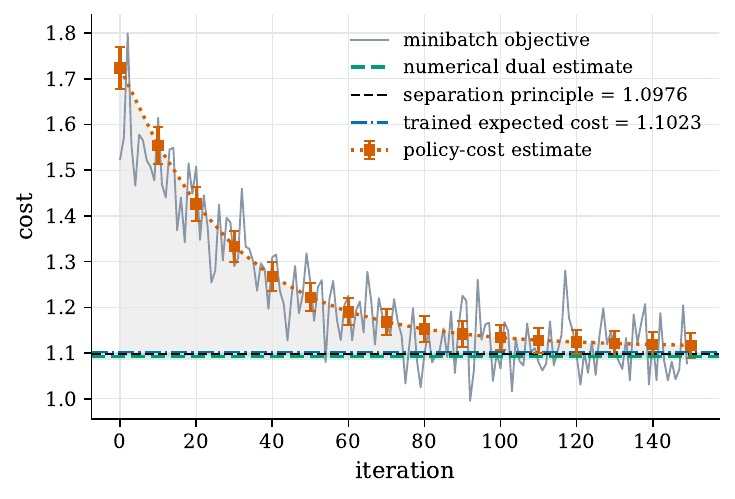}
	\caption{Training for the scalar LQG problem under regime A, $d=1$. Grey: minibatch objective. Orange: policy-cost estimate on fresh paths (mean $\pm$ $1.96$ standard errors). Green dashed: numerical dual estimate, maximized once after training, with a band of $\pm1.96$ standard errors. Black: separation-principle value. Blue: final trained expected cost, computed by Gaussian second-moment propagation.}
	\label{fig:lqg-regimeA-convergence}
\end{figure}
Figure~\ref{fig:lqg-regimeA-convergence} illustrates the training in the scalar setting. The minibatch objective decreases rapidly during the first iterations and then fluctuates around the final trained cost, as expected for a finite Monte Carlo batch. The final trained expected cost is $1.1023$, within $0.42\%$ of the separation-principle value $1.0976$. The trained expected cost is computed by Gaussian second-moment propagation for the discretized belief dynamics. The checkpoint estimates retain their Monte Carlo error bars. The numerical dual estimate, maximized once after training, equals $1.0924\pm0.0039$ (standard error), so the estimated primal--dual gap is $0.0099$, about $0.9\%$ of the value. The separation-principle value $1.0976$ lies between the two estimates.

We next retain the same fixed-noise regime-A structure and increase the state, control, and observation dimensions. For every tested dimension, the problem is dense and fully coupled. Thus, the dimensional study does not consist of independent copies of the scalar model.

Let
$
X_t\in\mathbb R^d,
\qquad
\alpha_t\in\mathbb R^d,
\qquad
Y_{t_n}\in\mathbb R^d.
$
The state dynamics are
\begin{equation}
	\dd X_t
	=
	(AX_t+B\alpha_t)\dd t
	+
	\sigma_X\,\dd W_t,
	\qquad
	X_0\sim\mathcal N(m_0,\Sigma_0).
	\label{eq:num-regA-multidim-state}
\end{equation}
and the observations are
\begin{equation}
	Y_{t_n}
	=
	CX_{t_n}+\Gamma\xi_n,
	\qquad
	\xi_n\stackrel{\mathrm{iid}}{\sim}\mcal N(0,I_d).
	\label{eq:num-regA-multidim-observation}
\end{equation}
The objective is
\begin{equation}
	J(\alpha)
	=
	\E\left[
	\int_0^T
	\left(
	\frac12 X_t^\top QX_t
	+\frac12\alpha_t^\top R\alpha_t
	\right)\dd t
	+\frac12 X_T^\top Q_TX_T
	\right],
	\label{eq:num-regA-multidim-cost}
\end{equation}
where $Q,Q_T\succeq0$ and $R\succ0$.

For each dimension $d$, the matrices
\begin{equation}
	A,\quad B,\quad C,\quad
	\sigma_X\sigma_X^\top,\quad
	\Gamma\Gamma^\top,\quad
	Q,\quad R,\quad Q_T,\quad
	\Sigma_0,
\end{equation}
are dense. They are generated using distinct orthogonal bases and controlled spectra. In particular, the matrices are not constructed from one common rotation of a diagonal problem and are not simultaneously diagonalizable. Consequently, the resulting system cannot be reduced to $d$ independent scalar LQG problems by a common change of coordinates. The spectra are chosen so that $A$ is stable and the covariance and cost matrices remain moderately conditioned as $d$ increases.

The belief remains Gaussian,
\begin{equation}
	\beliefstate_k=(m_k,\Sigma_k),
	\qquad
	\mu^{\beliefstate_k}=\mcal N(m_k,\Sigma_k).
\end{equation}
Since the observation schedule and observation-noise covariance are fixed, the covariance sequence $(\Sigma_k)_{k=0}^K$ and the Kalman gains are deterministic and can be precomputed before training. Between observation times, the belief prediction is
\begin{align}
	m_{k+1}^{-}
	&=
	m_k
	+
	\left(
	Am_k+B\alpha_k
	\right)\Delta t_k,
	\label{eq:num-regA-multidim-pred-mean}
	\\
	\Sigma_{k+1}^{-}
	&=
	\Sigma_k
	+
	\left(
	A\Sigma_k+\Sigma_kA^\top
	+\sigma_X\sigma_X^\top
	\right)\Delta t_k.
	\label{eq:num-regA-multidim-pred-cov}
\end{align}
At an observation time, the Kalman gain is
\begin{equation}
	K_{k+1}
	=
	\Sigma_{k+1}^{-}C^\top
	\left(
	C\Sigma_{k+1}^{-}C^\top+\Gamma\Gamma^\top
	\right)^{-1},
	\label{eq:num-regA-multidim-kalman-gain}
\end{equation}
and the posterior mean is updated according to
\begin{equation}
	m_{k+1}
	=
	m_{k+1}^{-}
	+
	K_{k+1}
	\left(
	Y_{\tau_{k+1}}-Cm_{k+1}^{-}
	\right).
	\label{eq:num-regA-multidim-update-mean}
\end{equation}
The posterior covariance is evaluated using the Joseph form
\begin{equation}
	\Sigma_{k+1}
	=
	(I-K_{k+1}C)
	\Sigma_{k+1}^{-}
	(I-K_{k+1}C)^\top
	+
	K_{k+1}\Gamma\Gamma^\top K_{k+1}^\top.
	\label{eq:num-regA-multidim-update-cov}
\end{equation}

The reduced running cost is
\begin{equation}
	\widetilde f(\mu^{\beliefstate_k},\alpha_k)
	=
	\frac12
	\left(
	m_k^\top Qm_k+\operatorname{Tr}(Q\Sigma_k)
	\right)
	+
	\frac12\alpha_k^\top R\alpha_k.
	\label{eq:num-regA-multidim-running-cost}
\end{equation}
The control induced by the potential is
\begin{equation}
	\alpha_k^\theta
	=
	-R^{-1}B^\top\eta_k^\theta,
	\qquad
	\eta_k^\theta
	=
	\int_{\mathbb R^d}
	\nabla_x U_k^\theta(x)
	\mu^{\beliefstate_k}(\dd x).
	\label{eq:num-regA-multidim-control}
\end{equation}
In the implementation, the multiplication by $R^{-1}$ is carried out using a linear solve rather than by forming the inverse explicitly.

The multidimensional training objective is therefore
\begin{equation}
	\cJ_\Delta^{(d)}(\theta)
	=
	\E\bigg[
	\sum_{k=0}^{K-1}
	\bigg(
	\frac12
	\left(
	m_k^\top Qm_k+\operatorname{Tr}(Q\Sigma_k)
	\right)
	+
	\frac12
	(\alpha_k^\theta)^\top
	R\alpha_k^\theta
	\bigg)\Delta t_k
	+
	\frac12
	\left(
	m_K^\top Q_Tm_K+\operatorname{Tr}(Q_T\Sigma_K)
	\right)
	\bigg].
	\label{eq:num-regA-multidim-loss}
\end{equation}
\paragraph{Quadratic potential and time parametrization.}
In the multidimensional experiment, we use a quadratic potential in the state variable $x\in\mathbb R^d$,
\begin{equation}
	U_k^\theta(x)
	=
	c_k
	+
	\ell_{x,k}^\top x
	+
	\frac12 x^\top H_{xx,k}x,
	\qquad H_{xx,k}=H_{xx,k}^\top.
	\label{eq:num-regA-multidim-potential}
\end{equation}
Its state gradient is
\begin{equation}
	\nabla_x U_k^\theta(x)
	=
	\ell_{x,k}
	+
	H_{xx,k}x.
\end{equation}
Since $X_k\mid\F^Y_{\tau_k}\sim\mcal N(m_k,\Sigma_k)$, the corresponding
gradient moment is
\begin{equation}
	\eta_k^\theta
	=
	\ell_{x,k}
	+
	H_{xx,k}m_k.
	\label{eq:num-regA-multidim-gradient}
\end{equation}
It can therefore be evaluated analytically.
The potential has no explicit observation-history input. The induced
control nevertheless uses the full observation history through the Kalman
mean $m_k$. The conditional quadratic running and terminal costs are
also evaluated analytically using the Gaussian moments. The remaining
outer expectation in the training objective is estimated by Monte Carlo
over Gaussian observation innovations.

The potential coefficients are shared over $K_\theta$ time blocks. Since the
potential input $x$ has dimension $d$, a complete quadratic polynomial
contains
$\binom{d+2}{2}$
coefficients per time block. Hence, the total number of trainable parameters
is
\begin{equation}
	N_\theta
	=
	K_\theta\binom{d+2}{2}.
	\label{eq:num-regA-total-parameter-count}
\end{equation}
In the experiments, we use $K_\theta=8$. This reduces the number of
time-dependent parameters but does not change the cost of propagating the
dense Kalman covariance matrices.

\begin{figure}[h]
	\centering
	\includegraphics[width=0.6\textwidth]
	{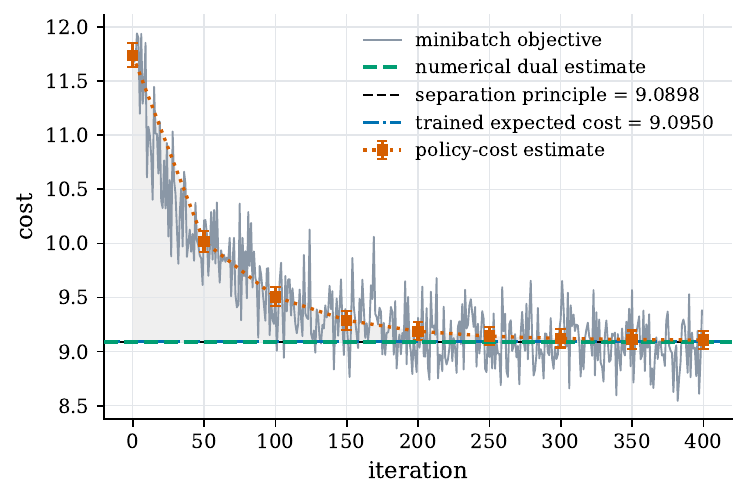}
	\caption{Training for the LQG problem under regime A, $d=10$. Same conventions as in Figure~\ref{fig:lqg-regimeA-convergence}, including the band of $\pm1.96$ standard errors around the numerical dual estimate. The black line is the separation-principle value.}
	\label{fig:lqg-regimeA-multidimensional}
\end{figure}

The trained policy is compared with the finite-horizon separation-principle controller.
The detailed $d=10$ experiment is reported in Figure~\ref{fig:lqg-regimeA-multidimensional} . For the additional
scalability study, we consider
$
d\in\{20,50,100\}.
$
For every value of $d$, a dense coupled
problem is constructed. Each problem is discretized using $40$ time steps
and three observation times. The potential uses eight block-constant time
slices and is trained for $400$ Adam iterations with batch size $64$.

The scalability study records the trained and reference costs, the relative
gap, the number of trainable coefficients, the wall-clock time per
optimization iteration, the total training time.
The number of coefficients in the complete quadratic potential grows
quadratically with the state dimension. The batched propagation of
conditional means and controls requires dense matrix products, while the
covariance propagation, innovation solves, and Riccati recursion involve
dense matrix factorizations and matrix--matrix products. Since the covariance
sequence and Kalman gains are deterministic in regime A, these quantities are
precomputed once and reused throughout training.

For $d=10$ the trained expected cost is $9.0950$, compared with the separation-principle value $9.0898$ (a relative gap of $0.057\%$). The numerical dual estimate, maximized once, is $9.0868\pm0.0029$, so the estimated primal--dual gap is $0.0082$, about $0.09\%$ of the value. The trained expected cost is computed by Gaussian second-moment propagation. The checkpoint estimates retain their Monte Carlo error bars.

\begin{table}[h]
\centering
\small
\caption{Accuracy and computational scaling for state-only quadratic potentials in dense regime-A LQG problems, trained for $400$ Adam updates. Costs are evaluated deterministically by Gaussian second-moment propagation for the discretized belief dynamics. The relative gap is $(J^\theta-J^{\mathrm{sep}})/J^{\mathrm{sep}}$. Timings are measured in the new runs with one CPU thread.}
\label{tab:lqg-regimeA-multidimensional}
\begin{tabular}{rrrrrrr}
\toprule
$d$ & parameters & $J^{\mathrm{sep}}$ & $J^\theta$ & relative gap & time/iteration & training time\\
 & & & & & $(\mathrm{s})$ & $(\mathrm{s})$\\
\midrule
  $20$ & $1\,848$ & $18.1582$ & $18.1679$ & $0.053\%$ & $0.0065$ & $2.62$ \\
  $50$ & $10\,608$ & $45.0251$ & $45.0514$ & $0.058\%$ & $0.0084$ & $3.34$ \\
  $100$ & $41\,208$ & $89.9387$ & $90.0045$ & $0.073\%$ & $0.0156$ & $6.26$ \\
\bottomrule
\end{tabular}
\end{table}

Figure~\ref{fig:lqg-regimeA-multidimensional} and Table~\ref{tab:lqg-regimeA-multidimensional} summarize the dimensional study. As expected, both the separation-principle cost and the trained cost increase with the dimension. The expected separation-principle values are $18.1582$, $45.0251$, and $89.9387$ for $d=20$, $50$, and $100$, while the corresponding trained expected costs are $18.1679$, $45.0514$, and $90.0045$. Both costs are computed by Gaussian second-moment propagation for the discretized belief dynamics. The relative gaps are $0.053\%$, $0.058\%$, and $0.073\%$, respectively. This indicates that the quality of the learned feedback does not deteriorate substantially over the tested dimensions.

The computational cost increases with the dimension but remains manageable. The number of trainable coefficients grows from $1\,848$ at $d=20$ to $41\,208$ at $d=100$. Over the same range, the measured time per optimization iteration increases from approximately $0.0065$ seconds to $0.0156$ seconds, while the total training time increases from $2.62$ seconds to $6.26$ seconds for $400$ updates. These timings were measured anew using one CPU thread. Thus, although the complete quadratic potential induces a substantial growth in the number of parameters, the measured runtime increases moderately over the tested range.

At $d=100$, the trained policy achieves an expected cost of $90.0045$, compared with $89.9387$ for the separation-principle controller, a relative gap of $0.073\%$. This result confirms that the method remains accurate for a dense, genuinely coupled $100$-dimensional problem rather than only for a collection of independent scalar systems.

The scalar example provides a primal--dual optimality diagnostic and verifies the basic optimization mechanism. The dense multidimensional study evaluates accuracy, parameter growth, runtime, and memory relative to the exact
separation-principle benchmark.

\subsection{Experiment 2: LQG regime B with controlled observation noise}
\label{subsec:lqg-regime-b}

The second experiment keeps the scalar linear-Gaussian state equation but makes the observation noise controlled. There are in total three controlled observations, at $t=0.25,0.5,0.75$,
\begin{equation}
	Y_{t_n}=cX_{t_n}+\beta_{t_n}\chi_n^d,
	\qquad \chi_n^d\sim\mcal N(0,1),
	\label{eq:num-regB-observation}
\end{equation}
with $\beta_{t_n}\in[\beta_{\min},\beta_{\max}]\subset(0,\infty)$. The cost is
\begin{equation}
	J(\alpha,\beta)=\E\bigg[
	\int_0^T\left(\frac12 qX_t^2+\frac12 r\alpha_t^2\right)\dd t
	+\sum_{n=1}^{3}\frac{\kappa_n}{\beta_{t_n}}
	+\frac12 q_T X_T^2
	\bigg].
	\label{eq:num-regB-cost}
\end{equation}
The continuous reduced formulas remain the same as in regime A. The key
difference is the discrete pre-posterior Hamiltonian, which balances the
sensing cost against the expected potential after the observation.
Thus regime B tests whether potential optimization can jointly produce
continuous feedback and controlled information acquisition. We compare
the resulting policy with the HJB and Pontryagin approaches
in~\cite{bayer2025cont,bayer2025pontryagin}.

 In this unconstrained LQG problem, an optimal sensing schedule can be
 chosen independently of the realized observations. The value separates
into a quadratic contribution in the conditional mean and a covariance
contribution, and the covariance recursion depends on the sensing
controls but not directly on the observation realizations. The sensing
decision is made before the current observation. Although dependence on
 previous observations is admissible, it is unnecessary for an optimal
schedule in this example.

Experiment 2 therefore does not test the learning of history-dependent
sensing. In Experiment 3, we impose bounds on the continuous control
$\alpha$. Control saturation breaks the quadratic separation, allowing
the value of information, and hence the optimal sensing decision, to
depend on the conditional mean and previous observations.

The potential is the structured quadratic
potential~\eqref{eq:structured-potential}, of the form
\[
U^\theta(t,x,z)
=a_\theta(t)x^2+b_\theta(t)xz+d_\theta(t)z^2.
\]
Its coefficients are produced by a neural network taking time as its
only input, with two hidden layers of width $32$ and Tanh activations.
The parametrization enforces
$U^\theta(T,x,z)=q_Tx^2/2$ and $d_\theta(t)>0$ for $t<T$.
The variable $z$ contains the latest observation, whereas the Kalman
mean incorporates the full observation history.

Unlike regime A, the potential must also provide a nontrivial
continuation contribution to the sensing decision. For the present
family, the pre-posterior expectation satisfies
\[
\E\!\left[
U^\theta(t_n,X_{t_n},Y_{t_n})
\mid\F^Y_{t_n^-}
\right]
=
\bigl(a_\theta(t_n)+cb_\theta(t_n)+c^2d_\theta(t_n)\bigr)
\bigl(m_{t_n^-}^2+\Sigma_{t_n^-}\bigr)
+d_\theta(t_n)\beta_{t_n}^2.
\]
Hence, its sensing-dependent part is
\[
d_\theta(t_n)\beta_{t_n}^2+\frac{\kappa_n}{\beta_{t_n}},
\]
yielding the closed-form minimizer
\[
\beta_{t_n}^\theta
=
\operatorname{clip}_{[\beta_{\min},\beta_{\max}]}
\left[
\left(\frac{\kappa_n}{2d_\theta(t_n)}\right)^{1/3}
\right].
\]
The observation-dependent part of the potential thus supplies the
continuation surrogate needed to balance sensing quality and cost.
If the potential depended only on $(t,x)$, its pre-posterior expectation
would be independent of $\beta$ by the tower property, and this
construction would select $\beta_{\max}$.

The problem data are $T=1$ (discretized with $350$ time steps, with observation dates rounded as recorded in Appendix~\ref{app:protocol}), $\kappa=(0.05,0.01,0.001)$, $\beta_{t_n}\in[0.02,1.5]$, and the observation memory has length one. The network is trained with Adam (learning rate $0.003$, $500$ updates, seed $1$, $\ell_2$ weight $10^{-5}$).

Conditional Gaussian moments are evaluated analytically. The outer
training expectation over the three observation innovations is
evaluated by an order-$2$ tensor Gauss--Hermite rule. This rule is exact for the quadratic innovation cost of
the present linear controller, up to floating-point error.
During training, sensing is selected by hard minimization over a discrete
grid, with a straight-through softmin gradient.
This straight-through update is a surrogate gradient and is not covered
by the exact discrete-adjoint identity of
Proposition~\ref{prop:theta-reduced-discrete-adjoint} for the hard-selection
forward map. At evaluation, sensing is selected by the closed-form rule of
Appendix~\ref{app:regime-b-beta-proof}.

 We compare saved checkpoints from two training runs, using respectively
 the grid selector and the continuous closed-form selector during
 training. Selection minimizes each checkpoint's own expected policy
 cost under the continuous selector. The selected checkpoint is update
 $500$ of the grid-selector run. No reference value, control, or Riccati
 coefficient is used to train or select the potential.
 Reported policy costs are evaluated by an order-$3$ tensor
 Gauss--Hermite rule, which is also exact for this quadratic innovation
 cost. They therefore have no Monte Carlo error, although the simulation
 retains time-discretization error.
 
The green line is obtained by a separate dual optimization over the same
neural potential family as the primal, initialized from the saved primal
weights. All $1\,219$ parameters remain trainable. The selected dual
checkpoint is evaluated on fresh innovation scenarios.
\begin{figure}[h]
	\centering
	\includegraphics[width=0.6\textwidth]{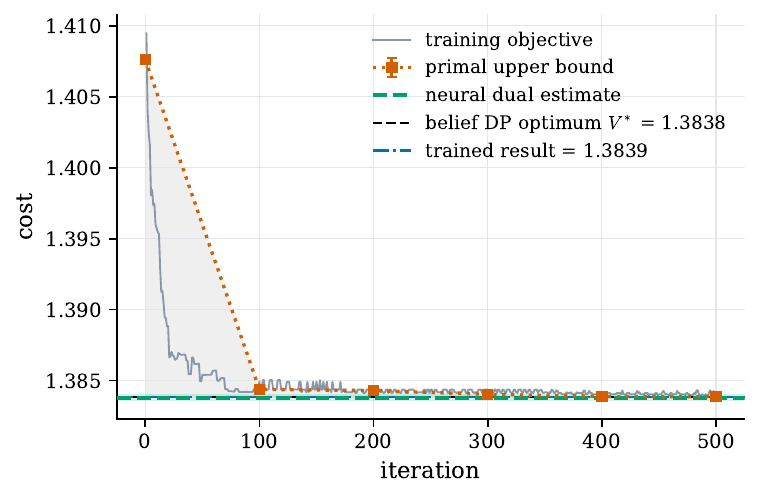}
	\caption{Training for the scalar LQG problem under regime B.
		Grey: training objective. Orange: expected checkpoint cost
		with the continuous closed-form sensing rule. Both expectations are
		integrated exactly for the present quadratic innovation costs, so
		no sampling error bars are shown for the primal costs. Green dashed:
        independently fitted neural dual estimate, with a band of
        $\pm1.96$ Monte Carlo standard errors. Black: numerical optimal
        value $V^*$ of the discretized belief problem.}
	\label{fig:lqg-regimeB-convergence}
\end{figure}

Figure~\ref{fig:lqg-regimeB-convergence} shows the training under
regime B. The expected policy cost decreases from $1.4076$ at
initialization to $1.38386$ after $500$ updates, which is $0.0023\%$
above the numerical optimum $V^*=1.38383$ of the discretized belief
problem. The reference is evaluated using a discrete Riccati
recursion and numerical optimization over the three deterministic
sensing levels.
The independently optimized neural dual estimate is
$1.38380\pm0.00012$, evaluated on $5\,000$ fresh innovation scenarios
(Monte Carlo standard error). The estimated primal--dual difference,
$6.1\times10^{-5}$, is smaller than this standard error and therefore
cannot be distinguished from sampling fluctuations.
Independent quadrature checks give consistent values.
These numerical evaluations do not provide a certified lower bound:
the Monte Carlo standard error measures sampling uncertainty only,
and does not account for integration errors or incomplete inner
minimization.

\begin{figure}[h]
	\begin{minipage}{.45\textwidth}
		\centering
		\includegraphics[width=0.9\textwidth]{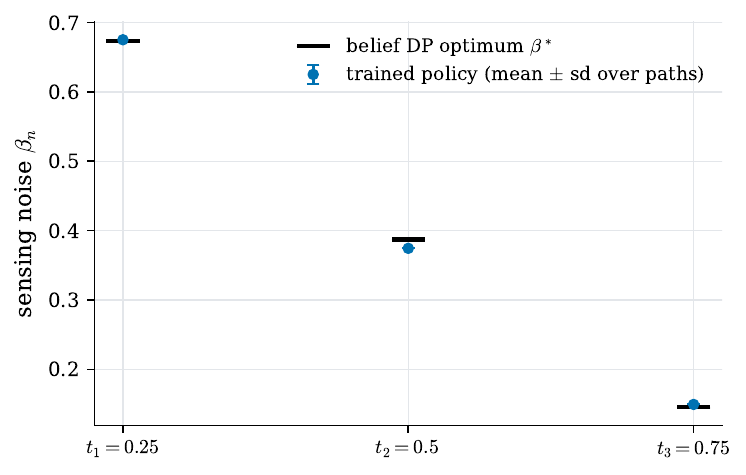}
		\captionof{figure}{Noise levels $\beta_{t_n}$ at each observation time for LQG regime B: learned policy (mean $\pm$ standard deviation over paths) and belief dynamic programming optimum.}
		\label{fig:regimeB-noise-control}
	\end{minipage}\hfill
	\begin{minipage}{.55\textwidth}
		\centering
		\begin{tabular}{lrrrr}
			\toprule
			Method & $\bar\beta_{t_1}$ & $\mathbb{E}[\bar\beta_{t_2}]$ & $\mathbb{E}[\bar\beta_{t_3}]$ & Cost\\
			\midrule
			HJB (Optimal)~\cite{bayer2025cont} & 0.677 & $0.389$ & $0.147$ & 1.383\\
			Pontryagin~\cite{bayer2025pontryagin} & 0.765 & 0.352 & 0.145 & 1.385\\
			Optimizer & $0.675$ & $0.375$ & $0.150$ & 1.384\\
			\bottomrule
		\end{tabular}
		\captionof{table}{Comparison of noise-control methods for LQG regime B.}
		\label{tab:regimeB-noise-control}
	\end{minipage}
\end{figure}

The comparison in Figure~\ref{fig:regimeB-noise-control} and Table~\ref{tab:regimeB-noise-control} shows that the learned optimizer produces noise levels close to those obtained by the HJB and Pontryagin approaches. The learned values are $0.675$, $0.375$ and $0.150$, compared with the belief dynamic programming optimum $0.674$, $0.387$ and $0.145$. The cost of the learned policy, $1.38386$, agrees with the optimal value $1.38383$ to $0.0023\%$, and is consistent with the HJB cost $1.383$ and the Pontryagin cost $1.385$ at the precision of the table.

The dispersion of the learned noise levels across trajectories is zero
up to rounding errors (below $10^{-16}$). This follows directly from
the chosen potential family: the closed-form sensing rule depends
only on $d_\theta(t_n)$, so every trajectory receives the same
$\beta_{t_n}$. The experiment therefore evaluates the joint
optimization of continuous feedback and a deterministic sensing
schedule, consistently with the structure of the unconstrained
LQG optimum.

\subsection{Experiment 3: LQG regime B with bounded control and controlled observation noise}
\label{subsec:lqg-regime-b-bounded}

The third experiment keeps the same scalar linear-Gaussian state equation and controlled observation model as in Experiment 2, but additionally imposes bounds on the continuous control. The observations are again given by
\begin{equation}
	Y_{t_n}=cX_{t_n}+\beta_{t_n}\chi_n^d,
	\qquad \chi_n^d\sim\mcal N(0,1),
\end{equation}
with $\beta_{t_n}\in[\beta_{\min},\beta_{\max}]\subset(0,\infty)$, while the continuous control is now constrained by
\begin{equation}
	\alpha_t\in[\underline\alpha,\overline\alpha].
\end{equation}
The cost remains
\begin{equation}
	J(\alpha,\beta)=\E\bigg[
	\int_0^T\left(\frac12 qX_t^2+\frac12 r\alpha_t^2\right)\dd t
	+\sum_{n=1}^{3}\frac{\kappa_n}{\beta_{t_n}}
	+\frac12 q_T X_T^2
	\bigg].
\end{equation}

The bound on the continuous control changes the structure of the problem in an essential way. In Experiment 2, the unconstrained LQG feedback preserves the quadratic structure of the value function, and the optimization of the sensing control $\beta_{t_n}$ depends only on the covariance. Consequently, the optimal sensing strategy is independent of the conditional mean and hence independent of the realized observation history.

When the continuous control is bounded, however, the optimal feedback saturates whenever the unconstrained control would lie outside the admissible interval. The resulting value function is no longer globally quadratic in the conditional mean. Consequently, the value of information can depend on both the conditional covariance and the conditional mean.

Since the conditional mean depends on previous observations, bounded
control allows the optimal sensing strategy to depend on the observation
history. For $n\geq2$, the sensing decision can therefore be represented
as
\begin{equation}
	\beta_{t_n}^\star
	=
	\beta_{t_n}^\star
	\left(
	\yhist^{[n-1]}
	\right),
\end{equation}

The main difference from Experiment 2 is that different observation
histories can now lead to different sensing decisions. Dispersion in
the learned noise levels can therefore reflect feedback on the
accumulated information, rather than merely numerical effects.
Experiment 3 tests the joint learning of bounded continuous feedback
and history-dependent sensing.

The numerical parameters are $\alpha_{\max}=0.6$,
$\beta_n\in[0.02,1.5]$, and
$(\kappa_1,\kappa_2,\kappa_3)=(0.05,0.01,0.001)$, with observations at
$(t_1,t_2,t_3)=(0.25,0.50,0.75)$, initial belief
$\mathcal N(0,1)$, and $200$ time steps.
The optimizer parametrizes one potential
$U_{\theta_n}(t,x,\yhist^{[n]})$ on each observation interval, namely
the polynomial--trigonometric potential~\eqref{eq:poly-trig-potential}
of Section~\ref{subsec:poly-trig-potential}, with \AnsParams{} trained
coefficients in total. Each slice receives the full available
observation history and no observation-memory truncation is used.
Training minimizes the primal objective without HJB or DP targets.

Expectations are approximated by Gauss--Hermite quadrature, yielding
a deterministic training objective. The induced controls satisfy
the prescribed bounds.
\begin{figure}[h]
\centering
\includegraphics[width=.75\textwidth]{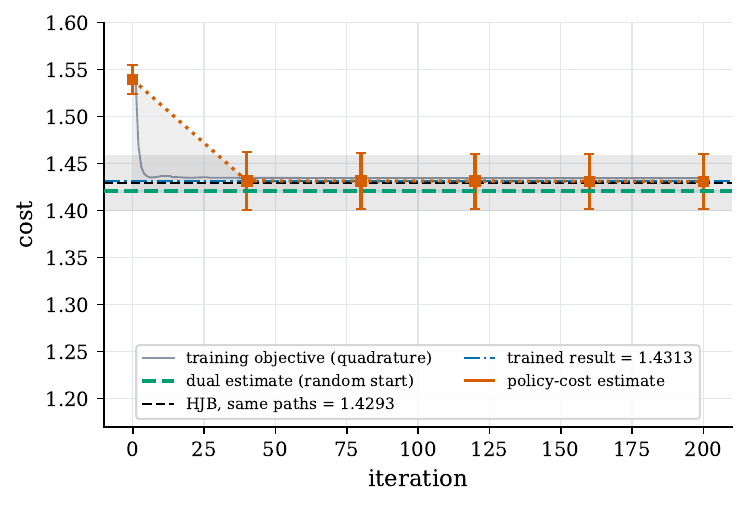}
\caption{Training of the polynomial--trigonometric potential
	(\AnsParams{} parameters) for bounded LQG regime B.
	Grey: deterministic training objective.
	Orange: checkpoint policy costs.
	Green dashed: independently optimized numerical dual estimate.
	Black dashed: HJB policy cost, evaluated using the same innovations
	as the checkpoints.
	Blue dash-dotted: cost of the validation-selected potential.
	Bands show $\pm1.96$ Monte Carlo standard errors.
	The shaded primal--dual gap is a numerical diagnostic,
	not a certified optimality gap.}
\label{fig:regimeB-bounded-convergence}
\end{figure}
Figure~\ref{fig:regimeB-bounded-convergence} shows the training. The primal cost drops within the first $40$ iterations and then stays within the error bars of the HJB cost. The dual estimate is a single value obtained after the training, by maximizing the dual separately from a random start, and it is therefore drawn as a horizontal line.

The potential~\eqref{eq:poly-trig-potential} is trained by minimizing
the primal cost for \AnsTrainIters{} AdamW iterations at learning
rate \AnsTrainLR{}, without HJB or DP targets.
The validation-selected policy costs
$J(\theta_{\mathrm{primal}})=\AnsJ\pm\AnsJse$ on
\AnsEvalPaths{} innovation scenarios.
Using the same standard-normal innovations, the HJB policy costs
$\AnsHJBPaired$, giving the paired difference
$\AnsPairedDiff\pm\AnsPairedDiffse$.
The sample does not resolve a cost difference.
The potential and HJB policy use $200$ and $800$ time steps,
respectively, so the comparison also includes discretization effects.

A separate coefficient vector in the same polynomial--trigonometric
family is optimized for the dual by approximate envelope-gradient
ascent (Section~\ref{subsec:dual-maximisation}).
Starting from independent Gaussian coefficients with standard
deviation \AnsRandScale{}, we perform \AnsSteps{} AdamW updates
at learning rate \AnsDualLR{}, using $256$ fresh innovation
scenarios per update.
Evaluation on \AnsDualMCPaths{} fresh scenarios gives
$D=\AnsD\pm\AnsDualSE$, with diagnostic gap
$J-D=\AnsGap$, about \AnsGapPct\% of the primal cost.
The dual estimate and the associated gap are numerical diagnostics,
not a certified lower bound or optimality gap. The Monte Carlo standard
error measures sampling uncertainty only. It does not account for
errors from local inner minimization or numerical integration.

\begin{figure}[h]
	\centering
	\includegraphics[width=\textwidth]{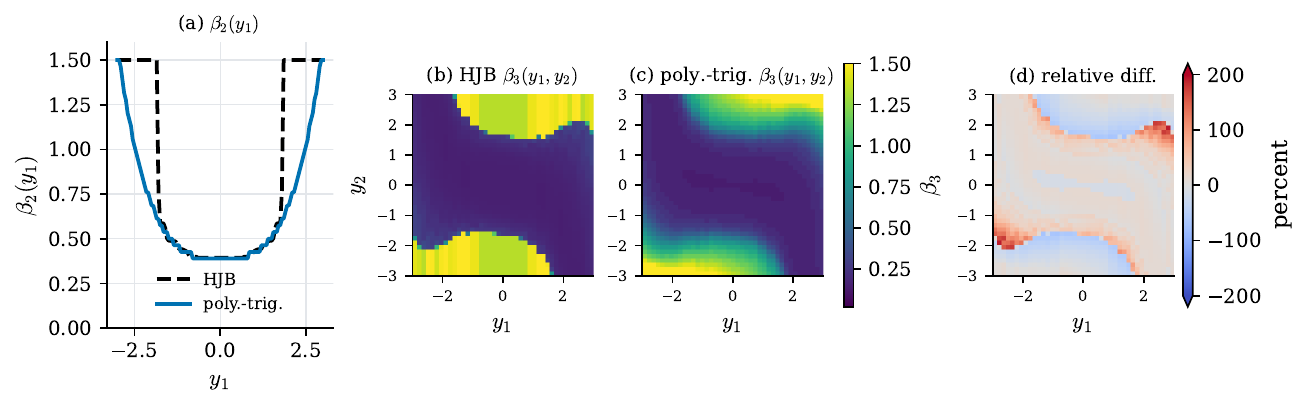}
	\caption{Sensing policies in Experiment~3.
		(a)~$\beta_2(y_1)$.
		(b), (c)~$\beta_3(y_1,y_2)$ for the HJB policy and the
		polynomial--trigonometric potential. Both rules use the same
		forecast beliefs generated by the potential-induced policy.
		(d)~Relative difference
		$(\beta_3-\beta_3^{\mathrm{HJB}})/\beta_3^{\mathrm{HJB}}$
		in percent.}
	\label{fig:regimeB-bounded-polytrig}
\end{figure}

Figure~\ref{fig:regimeB-bounded-polytrig} compares the sensing policies.
The first sensing control is deterministic:
$\beta_1=\AnsBetaOneAns$ for the potential and
$\beta_1^{\mathrm{HJB}}=\AnsBetaOneHJB$ for HJB.
For $\beta_2(y_1)$, both policies choose larger noise levels for
large $|y_1|$, but the potential approaches the upper bound more
gradually. The mean absolute difference over $y_1\in[-3,3]$ is
$0.2039$.
Both sensing rules are evaluated at the same forecast beliefs,
generated by the preceding potential-induced policy.
The $\beta_3$ surfaces have similar qualitative patterns but differ
near the transition regions, with mean absolute difference
$0.1375$ over the displayed observation grid.
The relative differences range from $-69.1\%$ to $191.2\%$.
Similar expected costs do not establish equality of the sensing
policies.

\subsection{Experiment 4: NLQG regime A with an obstacle penalty}
\label{subsec:nlqg-obstacle}

The fourth experiment keeps fixed observation noise but introduces a nonlinear, nonquadratic running penalty. The state and observation model are
\begin{equation}
  \dd X_t=(aX_t+b\alpha_t)\,\dd t+\sigma\,\dd W_t,
  \qquad
  Y_{t_n}=cX_{t_n}+\varepsilon\xi_n.
  \label{eq:num-nlqg-state-observation}
\end{equation}
The cost is
\begin{equation}
  J(\alpha)=\E\bigg[
  \int_0^T\left(\frac12 qX_t^2+\frac12 r\alpha_t^2+P(t,X_t)\right)\dd t
  +\frac12 q_T(X_T-x_\star)^2
  \bigg],
  \label{eq:num-nlqg-cost}
\end{equation}
where the obstacle penalty is
\begin{equation}
  P(t,x)=\Lambda\,\mathbbm 1_{[t_{\min},t_{\max}]}(t)
  \mathbbm 1_{\{x_{\min}\le |x|\le x_{\max}\}}.
  \label{eq:num-obstacle-penalty}
\end{equation}
For $X\sim\mcal N(m,\Sigma)$, the obstacle contribution is computed exactly through the Gaussian band probability (with standard normal CDF $\Phi_{\mathrm N}$),
\begin{align}
  \E[P(t,X)]
  &=\Lambda\,\mathbbm 1_{[t_{\min},t_{\max}]}(t)
  \Big[
  \Phi_{\mathrm N}\left(\frac{x_{\max}-m}{\sqrt \Sigma}\right)
  -\Phi_{\mathrm N}\left(\frac{x_{\min}-m}{\sqrt \Sigma}\right)
  \nonumber\\
  &\hspace{8em}
  +\Phi_{\mathrm N}\left(\frac{-x_{\min}-m}{\sqrt \Sigma}\right)
  -\Phi_{\mathrm N}\left(\frac{-x_{\max}-m}{\sqrt \Sigma   }\right)
  \Big].
  \label{eq:num-obstacle-expectation}
\end{align}
The continuous quadratic control term is eliminated as before:
\begin{equation}
  \widehat\alpha_k^\theta=-\frac{b}{r}\eta_k^\theta,
  \qquad
  \frac12 r(\widehat\alpha_k^\theta)^2=\frac12\rho(\eta_k^\theta)^2.
  \label{eq:num-nlqg-control-elim}
\end{equation}
The reduced stage cost is therefore
\begin{equation}
  \mathfrak f_k(\beliefstate_k,\widehat U_k^\theta)
  =\frac12 q(m_k^2+\Sigma_k)+\E[P(\tau_k,X_k)]+\frac12\rho(\eta_k^\theta)^2.
  \label{eq:num-nlqg-stage-cost}
\end{equation}
The potential is affine in $x$ before the observation and quadratic
in $(x,Y_{t_1})$ afterward, with separate coefficients at each grid
time. Conditional polynomial expectations are evaluated analytically
using Gaussian moments. Since there is only one observation, the
retained memory contains the full available observation history.
The experiment uses the following parameters:
\begin{equation}
\begin{gathered}
  T=1,
  \qquad K=300,
  \qquad N_o=1,
  \qquad t_1=0.5,
  \qquad \Delta t=1/300,\\
  a=0,
  \quad b=1,
  \quad \sigma=0.5,
  \quad c=1,
  \quad \varepsilon=0.1,
  \quad q=0,
  \quad r=2,
  \quad q_T=20,
  \quad x_\star=0,\\
  m_0=0,       
  \quad \Sigma_0=0.01,
  \quad \Lambda=1000,
  \quad [t_{\min},t_{\max}]=[0.3,0.6],
  \quad [x_{\min},x_{\max}]=[0.1,2.0].
\end{gathered}
\label{eq:num-nlqg-params}
\end{equation}

Although the obstacle penalty is discontinuous in $x$, its Gaussian
expectation is smooth in the conditional mean for positive variance.
The resulting reduced objective is nevertheless nonconvex and can
have unfavorable local minima.

We initialize the potential by simulating $200$ state trajectories
under a provisional feedback, computing empirical cost-to-go values
backward, and regressing them onto the polynomial basis.
One damped regression update produces the initial coefficients.
These simulated trajectories are used only for initialization.
Subsequent training and cost evaluation use the exact Gaussian
belief representation.

The initialized potential has cost approximately $305$, which is
then reduced by gradient-based optimization.

\begin{figure}[h]
	\centering
	\includegraphics[width=0.6\textwidth]{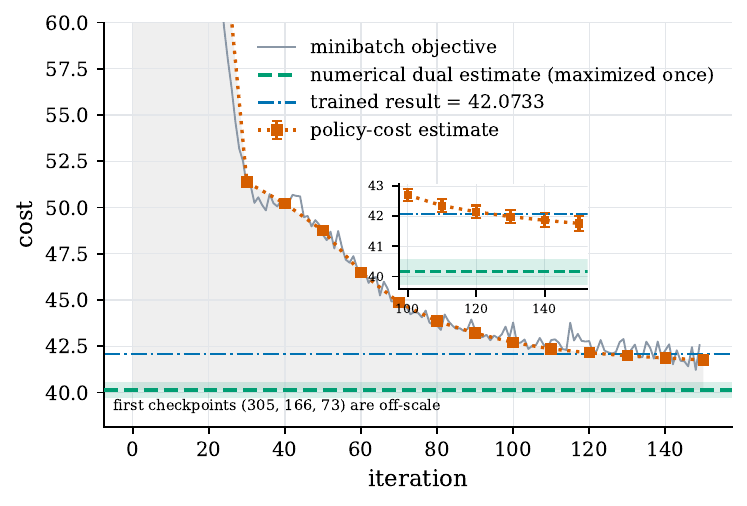}
\caption{Training for the obstacle problem.
	Grey: minibatch objective.
	Orange: checkpoint policy costs.
	Green dashed: independently maximized numerical dual estimate,
	evaluated on $2048$ fresh innovation scenarios.
	Bands and error bars show $\pm1.96$ Monte Carlo standard errors.
	The inset shows iterations $100$--$150$.
	The initial cost, approximately $305$, is outside the displayed range.
	Local minimization of the nonconvex relaxed problem makes the
	primal--dual gap a numerical diagnostic, not a certificate.}
	\label{fig:nlqg-1d-certificate}
\end{figure}

\begin{figure}[h]
	\centering
	\includegraphics[width=0.6\textwidth]{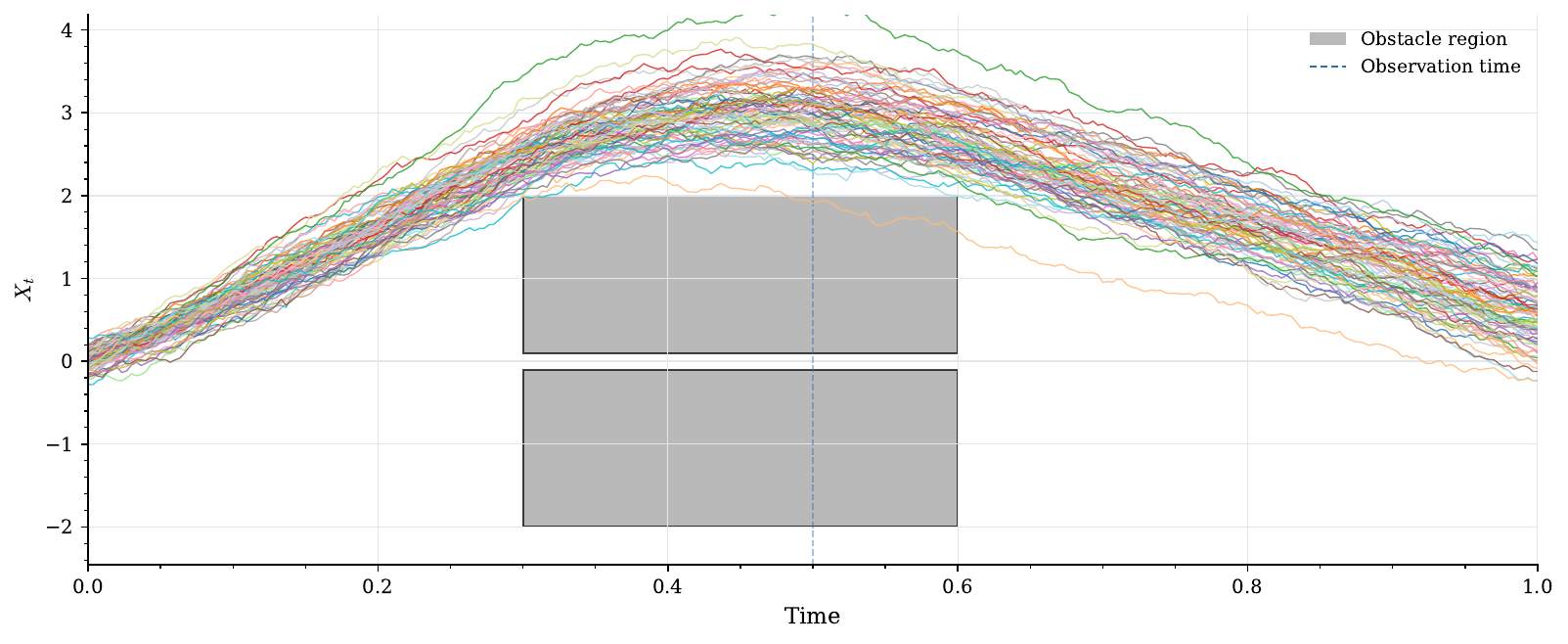}
	\caption{Controlled trajectories of the obstacle problem under the trained potential.}
	\label{fig:nlqg-1d-paths}
\end{figure}

The results show that the restricted polynomial potential produces an
effective avoidance policy for this nonquadratic problem.
Starting from the regression initialization, whose cost is approximately
$305$, optimization reduces the expected cost to $42.07\pm0.06$.
The separately maximized numerical dual estimate, evaluated on
$2048$ fresh innovation scenarios, is $40.15\pm0.21$.
The resulting diagnostic gap is $1.93$, approximately $4.6\%$ of
the primal cost.
Since the nonconvex relaxed problems are solved by local multistart
minimization, this estimate is not a certified lower bound.
The reported standard errors measure sampling uncertainty, while
inner-optimization and numerical-integration errors remain outside
these error bars.
The gap therefore provides an empirical diagnostic rather than a
proof of global optimality. The polynomial parametrization restricts
the policies that training can represent, and the numerical results
also retain time-discretization error. Gaussian filtering itself is
exact for the present dynamics and observation model.

The controlled trajectories illustrate how the policy balances
obstacle avoidance, control effort and the terminal penalty.
The state distribution is moved toward the positive side before the
penalty becomes active, so that most paths lie above $x=2$ during
$[0.3,0.6]$. After the penalty window, the control directs the state
back toward the terminal target $x_\star=0$, without requiring every
trajectory to reach it.
Because the first observation arrives at $t_1=0.5$, the initial
avoidance maneuver is necessarily planned without observations.
Feedback from the measurement can only adjust the subsequent motion.
The positive direction is one of two symmetry-related avoidance
strategies permitted by the symmetric obstacle and terminal cost.

This example demonstrates that potential optimization can produce
useful controls beyond quadratic running costs, while retaining
analytical Gaussian expectations and a low-degree polynomial
parametrization.
The regression initialization is part of the tested procedure.
The experiment does not establish that the same performance would
be obtained from an arbitrary initialization.
We additionally compare potential parametrizations in dimensions
$d\in\{2,4,10\}$. For $i=1,\ldots,d$, the state dynamics are
\[
\dd X_t^i
=
\left[
0.6\left(\frac1d\sum_{j=1}^d X_t^j-X_t^i\right)
+\alpha_t^i
\right]\dd t
+0.5\,\dd W_t^i,
\qquad
X_0\sim\mathcal N(0,0.01I_d).
\]
The drift couples the coordinates by pulling each toward their average.
Ten observations are available at $t_n=n/11$, with
\[
Y_{t_n}=X_{t_n}+0.1\,\xi_n,
\qquad
\xi_n\stackrel{\mathrm{iid}}{\sim}\mathcal N(0,I_d).
\]
The objective is
\[
J(\alpha)
=
\E\left[
\int_0^1
\left(
\frac r2\|\alpha_t\|^2+\sum_{i=1}^d P(t,X_t^i)
\right)\dd t
+\frac{q_T}{2}\|X_1\|^2
\right],
\]
where $r=2$, $q_T=20$, and $P$ is the scalar obstacle
penalty~\eqref{eq:num-obstacle-penalty} with the parameters
in~\eqref{eq:num-nlqg-params}.

Table~\ref{tab:nlqg-representation-comparison} compares four potential
parametrizations. The observation-free controller uses no measurements.
The full-history and moment-based controllers use all observations
received so far, represented respectively by the raw history and the
Gaussian mean.

The potential families are chosen to provide practical parametrizations
for each information representation, using neural networks of comparable
width and depth. Since the architectures differ, the comparison
assesses these combinations of representation and parametrization
rather than isolating their individual effects.
\begin{table}[htbp]
	\centering
	\begin{tabular}{lcrrr}
		\hline
		Potential & Observations used
		& $d=2$ & $d=4$ & $d=10$ \\
		\hline
		$U(t,x)$
		& None
		& 84.26 & 168.89 & 541.62 \\
				$U(t,x,Y_{n-1},Y_{n})$
		& Two, delayed
		& 72.88 & 144.87 & 360.79 \\
		$U(t,x,\mf{Y}^{[n]})$
		& All ten
		& 70.83 & 149.99 & 519.58 \\
		$U(t,x,m)$
		& All ten
		& 65.11 & 132.66 & 336.08 \\
		\hline
		Full-history numerical dual estimate
		& All ten
		& 58.33 & 116.41 & 290.65 \\
		\hline
	\end{tabular}
\caption{Total expected costs for the multidimensional obstacle problem
	under different potential parametrizations and observation inputs.
	Results use one training seed and $8192$ evaluation scenarios.
	The final row gives a common full-history numerical dual estimate
	for the discretized problem. It is not a certified lower bound.}
	\label{tab:nlqg-representation-comparison}
\end{table}

\subsection{Experiment 5: nonsmooth bounded-control LQG and Hamiltonian smoothing}
\label{subsec:nonsmooth}

This experiment illustrates the role of Hamiltonian smoothing in a simple
bounded-control LQG problem.  The purpose is to isolate a case where the reduced
Hamiltonian is continuous but not differentiable with respect to the trainable
potential coefficients.  We compare the training of the same quadratic potential
using either the original nonsmooth Hamiltonian or a smooth Hamiltonian
surrogate.  The quality of the trained potential is then assessed through the
corresponding policy-cost estimate, numerical dual estimate, and diagnostic gap.

The hidden state is scalar and evolves according to
\begin{equation}
	\dd X_t=(aX_t+b\alpha_t)\dd t+\sigma\dd W_t,
	\qquad X_0\sim\mcal N(m_0,\Sigma_0),
	\label{eq:num-nonsmooth-state}
\end{equation}
with bounded control
\begin{equation}
	\alpha_t\in[-\alpha_{\max},\alpha_{\max}].
\end{equation}
The state is observed at fixed observation times through
\begin{equation}
	Y_{t_n}=cX_{t_n}+\varepsilon \xi_n,
	\qquad \xi_n\sim\mcal N(0,1),
	\label{eq:num-nonsmooth-observation}
\end{equation}
where the observation noises are independent.  The conditional law is represented
by the Gaussian reduced state \(\beliefstate_k=(m_k,\Sigma_k)\), whose prediction and update
steps are given by the standard scalar Kalman formulas.

The cost functional is
\begin{equation}
	J(\alpha)=
	\E\left[
	\int_0^T
	\left(
	\frac q2(m_t^2+\Sigma_t)+\kappa\alpha_t
	\right)\dd t
	+
	\frac{q_T}{2}\left((m_T-x_\star)^2+\Sigma_T\right)
	\right].
	\label{eq:num-nonsmooth-cost}
\end{equation}
The affine term \(\kappa\alpha_t\), combined with the bounded action set, is what
creates the nonsmooth reduced Hamiltonian.

\begin{figure}[h]
	\centering
	\includegraphics[width=0.6\textwidth]
	{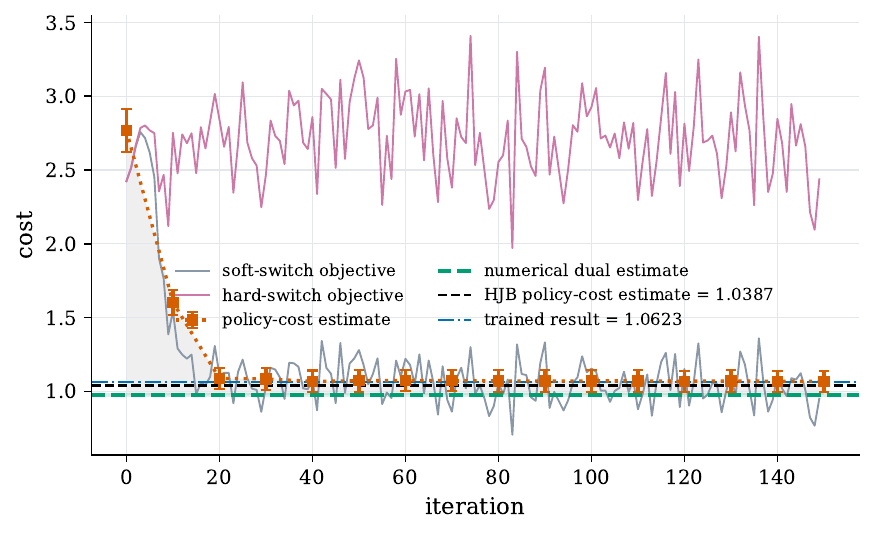}
	\caption{Hard- and soft-switch training for the bounded-control
		problem, seed $0$.
		Pink and grey: minibatch costs under the hard and smooth controllers.
		Orange: soft-trained potentials deployed with the hard switch.
		Green dashed: separately optimized numerical dual estimate.
		Black dash-dotted: HJB policy reference. Blue dash-dotted: cost of the final soft-trained potential,
		evaluated with the exact hard-switch controller.
		Bands and error bars show $\pm1.96$ Monte Carlo standard errors.}
	\label{fig:lqg-1d-nonsmooth-soft-certificate}
\end{figure}

We use the time-dependent quadratic potential ansatz
\begin{equation}
	\widehat U_k^\theta(x)=\frac12 P_kx^2+r_kx+c_k^U.
	\label{eq:num-nonsmooth-quadratic-potential}
\end{equation}
Only the coefficients \((P_k,r_k)\) affect the induced control, since
\(c_k^U\) disappears after differentiating with respect to \(x\).  In terms of
the reduced mean \(m_k\), the nonsmooth Hamiltonian contribution is
\begin{equation}
	-\alpha_{\max}
	\left|
	b(P_km_k+r_k)+\kappa
	\right|.
	\label{eq:num-nonsmooth-hamiltonian}
\end{equation}
This expression has a kink at the switching surface
\begin{equation}
	b(P_km_k+r_k)+\kappa=0,
\end{equation}
which makes the training objective difficult to optimize using first-order
methods.  The smoothed training replaces this term by the soft-Hamiltonian
surrogate
\begin{equation}
	-\delta
	\log\left(
	2\cosh\left(
	\frac{\alpha_{\max}\{b(P_km_k+r_k)+\kappa\}}{\delta}
	\right)
	\right),
	\label{eq:num-soft-hamiltonian}
\end{equation}
where \(\delta>0\) controls the smoothing level.  Both training curves in
Figure~\ref{fig:lqg-1d-nonsmooth-soft-certificate} therefore optimize the same
potential coefficients \((P_k,r_k)\). The only difference is the Hamiltonian used
during training.

The numerical parameters are
\begin{equation}
	\begin{gathered}
		T=1,
		\qquad K=60,
		\qquad N_o=3,
		\qquad t_n\in\{0.25,0.5,0.75\},\\
		a=0,
		\quad b=1,
		\quad \sigma=0.5,
		\quad c=1,
		\quad \varepsilon=0.1,
		\quad q=2,
		\quad q_T=2,
		\quad \kappa=1,
		\quad \alpha_{\max}=1,\\
		m_0=0,
		\quad \Sigma_0=1,
		\quad x_\star=0.
	\end{gathered}
	\label{eq:num-nonsmooth-params}
\end{equation}

Figure~\ref{fig:lqg-1d-nonsmooth-soft-certificate} compares hard-
and soft-switch training using the same potential family,
initialization and learning rate.
With the hard switch, automatic differentiation through the sign
function gives zero gradients, so the potential coefficients remain
unchanged. The oscillations of the pink curve therefore reflect
minibatch sampling rather than changes in the policy.
Across the three seeds, the hard-switch policies have an average
evaluated cost of approximately $2.71$.

The smooth controller provides nonzero gradients and allows the
training cost to decrease substantially.
When the resulting potentials are deployed with the exact hard
switch, their costs are $1.062$, $1.017$ and $1.029$, with Monte Carlo
standard errors between $0.016$ and $0.020$.
These values are close to the HJB grid value $1.037$ and the
Monte Carlo cost of the HJB policy, $1.039\pm0.006$.
The comparison therefore shows that smoothing can enable effective
gradient-based training while retaining good performance under the
original hard-switch controller.

For seed $0$, separate optimization over the six quadratic-potential
coefficients entering the observation penalties gives a numerical
dual estimate of $0.9756\pm0.0082$ on $4\,096$ fresh Monte Carlo
scenarios. The corresponding diagnostic gap is approximately
$0.0867$.
The pathwise inner problems are strictly convex box-constrained
quadratic programs, which permits reliable numerical minimization.
Global solution of these convex inner problems controls only the
minimization step. Solver tolerances, penalty evaluation and approximation
of the outer expectation remain sources of error. Without one-sided
control of these errors, the reported estimate is not a certified
finite-computation lower bound. Its Monte Carlo standard error measures
sampling uncertainty only.
The main conclusion is that smoothing supplies useful training
gradients in this example, and that the learned potential remains
effective when deployed with the exact hard switch.

\section{Conclusion}
\label{sec:conclusion}

We presented a Hamiltonian-based optimization framework for partially observed stochastic control. The method replaces direct optimization over policies by optimization over a pathwise potential. After the reduced Hamiltonians are evaluated, this potential determines the numerical belief recursion and the on-policy cost. After discretization and parametrization, the method becomes a stochastic optimization problem that can be trained by automatic differentiation or an equivalent discrete adjoint.

The central computational question is how to evaluate and differentiate the reduced Hamiltonians. In structured scalar LQG examples, exact Gaussian filtering gives explicit formulas for the continuous reduction and the Kalman update. In controlled-sensing regime B, an additional quadratic structure yields an explicit cubic-root formula for the sensing channel. In the obstacle example, the state and observation model remains linear-Gaussian, and the nonquadratic penalty is evaluated using Gaussian probabilities. In nonsmooth bounded-control problems, smoothing provides useful gradients and enables first-order optimization where the exact hard switch has zero or unstable gradients.

The numerical experiments support the following picture. In the LQG problems the trained potential reaches a cost close to the optimal one: the relative gap is $0.42\%$ in the scalar fixed-noise problem and stays below one percent from $d=10$ up to $d=100$, and in the controlled-sensing problem the structured potential matches the optimal value of the belief dynamic program. The separately optimized numerical dual estimates are close to the reference values. The controlled-sensing comparison uses an independently optimized dual in the same neural family as the primal. In the bounded-control problem trained with smoothing, the convex inner problems can be solved globally to numerical tolerance, but the outer expectation and solver errors still require assessment. The trained potential, deployed with the exact hard switch, is within about two standard errors of the HJB value. These numerical agreements do not certify the computed dual estimates: the theoretical lower bound applies to the exact dual functional, and a computed certificate also requires one-sided control of integration, sampling and solver errors. Monte Carlo standard errors account only for sampling uncertainty. In the obstacle problem the primal--dual gap of about $4.6\%$ is a diagnostic rather than a certificate. The same holds for the bounded controlled-sensing problem, where the gap is about \AnsGapPct\%.

The limitations of the present study include the restricted potential families, experiment-dependent restrictions on observation inputs, time discretization, and numerical errors in dual evaluation. Length-one observation memory is used in some parametrizations, whereas the bounded controlled-sensing example uses all available observations. Particle and ensemble belief representations, and a systematic study of history-dependent sensing, are left for future work.

\appendix

\section{Derivations used by the numerical formulation}
\label{app:derivations}

This appendix derives the belief updates, Hamiltonian cost slices,
Gaussian expectation formulas, controlled-sensing reduction,
Hamiltonian smoothing error bounds, and discrete adjoint
underlying the numerical setup to reproduce the results.

\subsection{Belief update and the pre-posterior tower property}
\label{app:bayes-proof}

Let $\mu^-$ be the forecast law of $X$ immediately before an observation and suppose that, conditionally on $X=x$, the observation $Y$ has density $\pi(y\mid x,\beta)$. For a bounded measurable test function $\varphi$, Bayes' rule gives
\begin{equation}
  \E[\varphi(X)\mid Y=y]
  =\frac{\int \varphi(x)\pi(y\mid x,\beta)\mu^-(\dd x)}
  {\int \pi(y\mid x,\beta)\mu^-(\dd x)}.
  \label{eq:app-bayes-test}
\end{equation}
Therefore the posterior measure is
\begin{equation}
  \mu^+(\dd x)
  =\frac{\pi(y\mid x,\beta)\mu^-(\dd x)}
  {\int \pi(y\mid z,\beta)\mu^-(\dd z)}.
  \label{eq:app-bayes-measure}
\end{equation}
The predictive density of the observation is
\begin{equation}
  p_Y(y)=\int \pi(y\mid x,\beta)\mu^-(\dd x).
  \label{eq:app-predictive-density}
\end{equation}
Using \eqref{eq:app-bayes-test} and \eqref{eq:app-predictive-density}, we obtain
\begin{equation}
  \E\big[\ip{\varphi}{\mu^+}\big]
  =\int \left(\frac{\int \varphi(x)\pi(y\mid x,\beta)\mu^-(\dd x)}{p_Y(y)}\right)p_Y(y)\dd y
  =\ip{\varphi}{\mu^-}.
  \label{eq:app-tower-posterior}
\end{equation}
This is the pre-posterior tower property used in the discrete Hamiltonian. Before observing $Y$, the continuation value is an expectation over all possible posterior beliefs generated by the future observation.

\subsection{Cost slices from augmented Hamiltonians}
\label{app:cost-slice-proof}

We prove the continuous identity under the assumptions of
Proposition~\ref{prop:cost-slice-main}. The discrete case is identical. Write
\begin{equation}
  F(a,\gamma)=\ip{\cG_aU+\gamma f_t(\cdot,a)}{\mu},
  \qquad
  \tH^c_t(\mu,U,\gamma)=\inf_{a\in\cA}F(a,\gamma).
\end{equation}
Let $a_\gamma$ be a minimizer and set
$B(a):=\ip{f_t(\cdot,a)}{\mu}$. Since $F$ is affine in $\gamma$,
optimality at $1+h$ and at $1$ gives, for $h\ne0$ sufficiently small,
\begin{equation}
  hB(a_{1+h})
  \le \tH^c_t(\mu,U,1+h)-\tH^c_t(\mu,U,1)
  \le hB(a_1).
\end{equation}
Compactness and continuity imply that every limit point of $a_{1+h}$
minimizes $F(\cdot,1)$. Uniqueness, or agreement of $B$ on all such
minimizers, therefore yields
\begin{equation}
  B(a_{1+h})\longrightarrow B(a_1)
  \qquad\text{as }h\longrightarrow0.
\end{equation}
Dividing the bounds by $h$, reversing them when $h<0$, proves
\begin{equation}
  \partial_\gamma\tH^c_t(\mu,U,1)=\ip{f_t(\cdot,a_1)}{\mu}.
\end{equation}
This is the envelope identity: it differentiates the optimized value
without requiring a derivative of $a_\gamma$. The same argument applies
to the discrete Hamiltonian, whose dependence on $\gamma$ is also affine.
If minimizing actions have different cost contributions, a classical
$\gamma$-derivative need not exist and this cost slice cannot be identified
without an additional qualification.

\subsection{Relation between \texorpdfstring{$\cJ(U)$}{J(U)}
	and the original cost}
\label{app:cJ-versus-J-proof}

Let $U$ be an admissible potential, with induced controls
$(\alpha^U,\beta^U)$ and exact conditional state law $\mu_t^U$.
Assume the costs are integrable and the cost slice identifications
\eqref{eq:cJ-J-running-identification}--
\eqref{eq:cJ-J-discrete-identification} hold.

Since $\alpha_t^U$ is observation-adapted and $\beta_{t_n}^U$
is measurable before the observation at $t_n$, conditional expectation
gives
\begin{align}
	\E[f_t(X_t,\alpha_t^U)\mid\mathbb F_t^Y]
	&=\ip{f_t(\cdot,\alpha_t^U)}{\mu_t^U},\\
	\E[c_{t_n}(X_{t_n},\beta_{t_n}^U)
	\mid\mathbb F_{t_n^-}^Y]
	&=\ip{c_{t_n}(\cdot,\beta_{t_n}^U)}{\mu_{t_n^-}^U},\\
	\E[g(X_T)\mid\mathbb F_T^Y]
	&=\ip{g}{\mu_T^U}.
\end{align}
Here continuity of the state process identifies the conditional law
of $X_{t_n}$ before assimilation with $\mu_{t_n^-}^U$.
The tower property and the cost slice identities therefore yield
\begin{align}
	J(\alpha^U,\beta^U)
	&=\E\bigg[
	\int_0^T\mathfrak f_t(\mu_t^U,U_t)\,\dd t
	+\sum_{n=1}^{N_o}
	\mathfrak c_n(\mu_{t_n^-}^U,U_{t_n})
	+\ip{g}{\mu_T^U}
	\bigg]\\
	&=\cJ(U).
\end{align}
This proves \eqref{eq:cJ-equals-J-induced}.
For fixed sensing parameters, the same argument applies with
$\beta$ fixed. The discrete cost term is omitted when no
observation-stage cost is present.

\subsection{Proof of the discrete adjoint representation}
\label{app:-first-adjoint-proof}

\begin{proof}
	First relate the Hamiltonian notation to the reduced objects. By
	\eqref{eq:theta-running-hamiltonian-slice} and
	\eqref{eq:theta-observation-hamiltonian-slice},
	\begin{align*}
		\nabla_\beliefstate\mathfrak f_k
		&=\nabla_\beliefstate\partial_\gamma\widehat{\mcal H}_k^{c,\theta},
		&D_\theta\mathfrak f_k
		&=D_\theta\partial_\gamma\widehat{\mcal H}_k^{c,\theta},\\
		\nabla_\beliefstate\mathfrak c_n^\theta
		&=\nabla_\beliefstate\partial_\gamma\overline{\mcal H}_n^{d,\theta},
		&D_\theta\mathfrak c_n^\theta
		&=D_\theta\partial_\gamma\overline{\mcal H}_n^{d,\theta},
	\end{align*}
	with the arguments and $\gamma=1$ understood. In particular,
	\begin{align}
		\nabla_\beliefstate\partial_\gamma\overline{\mcal H}_n^{d,\theta}
		&=\nabla_\beliefstate\partial_\gamma\widehat{\mcal H}_n^{d,\theta}
		+(\partial_\beliefstate G_n^\theta)^\top
		\nabla_y\partial_\gamma\widehat{\mcal H}_n^{d,\theta},
		\label{eq:theta-hamiltonian-observation-state}\\
		D_\theta\partial_\gamma\overline{\mcal H}_n^{d,\theta}
		&=D_\theta\partial_\gamma\widehat{\mcal H}_n^{d,\theta}
		+(D_\theta G_n^\theta)^\top
		\nabla_y\partial_\gamma\widehat{\mcal H}_n^{d,\theta}.
		\label{eq:theta-hamiltonian-observation-parameter}
	\end{align}
	
	The right-hand sides are evaluated at the generated observation, with
	the derivatives of the uncomposed Hamiltonian taken at fixed $y$.
	The forward maps are the full numerical maps constructed from the
	Hamiltonian forward derivatives: prediction includes the chosen time
	scheme, and the observation map includes normalization and retained-history
	updates. They are not identified with an unnormalized Hamiltonian
	derivative. Consequently their full Jacobians must be used.
	
	Define the forward sensitivities
	\[
	Q_k^\theta:=\frac{\dd\widehat\beliefstate_k^\theta}{\dd\theta},
	\qquad
	Q_{k+1^-}^\theta:=\frac{\dd\widehat\beliefstate_{k+1^-}^\theta}{\dd\theta},
	\qquad
	R_n^\theta:=\frac{\dd Y_{t_n}^\theta}{\dd\theta}.
	\]
	Since the initial condition is fixed, $Q_0^\theta=0$.
	Differentiating the prediction step gives
	\begin{equation}
		Q_{k+1^-}^\theta
		=\partial_\beliefstate\widehat\Phi_k^\theta\,Q_k^\theta
		+D_\theta\widehat\Phi_k^\theta.
		\label{eq:theta-forward-sensitivity-scheme}
	\end{equation}
	
	At an observation time, the chain rule gives
	\begin{align}
		R_n^\theta
		&=\partial_\beliefstate G_n^\theta\,Q_{k_n^-}^\theta
		+D_\theta G_n^\theta,
		\label{eq:theta-observation-sensitivity}\\
		Q_{k_n}^\theta
		&=\partial_\beliefstate\widehat\Psi_n^\theta\,Q_{k_n^-}^\theta
		+D_\theta\widehat\Psi_n^\theta
		+\partial_y\widehat\Psi_n^\theta\,R_n^\theta
		\nonumber\\
		&=\partial_\beliefstate\overline\Psi_n^\theta\,Q_{k_n^-}^\theta
		+D_\theta\overline\Psi_n^\theta.
		\label{eq:theta-forward-sensitivity-update}
	\end{align}
	
	All quantities are evaluated along the simulated trajectory. The first
	line of \eqref{eq:theta-forward-sensitivity-update} makes explicit the
	additional derivative through $Y_{t_n}^\theta$. Its second line incorporates
	that derivative into the composed map. At non-observation steps,
	$Q_{k+1}^\theta=Q_{k+1^-}^\theta$.
	
	Differentiating the objective and using the cost slice identifications
	above yields
	\begin{align}
		\nabla_\theta\mcal J_\Delta(\widehat U^\theta)
		={}&\mathbb E\Bigg[
		\sum_{k=0}^{K-1}
		\Bigl((Q_k^\theta)^\top\nabla_\beliefstate\mathfrak f_k
		+D_\theta\mathfrak f_k\Bigr)\Delta t_k
		\nonumber\\
		&+\sum_{n=1}^{N_o}
		\Bigl((Q_{k_n^-}^\theta)^\top\nabla_\beliefstate\mathfrak c_n^\theta
		+D_\theta\mathfrak c_n^\theta\Bigr)
		+(Q_K^\theta)^\top\nabla_\beliefstate\mathfrak g
		+D_\theta\mathfrak g\Bigg].
		\label{eq:theta-diff-objective}
	\end{align}
	
	The observation derivatives in this expression are already included in
	$\nabla_\beliefstate\mathfrak c_n^\theta$ and
	$D_\theta\mathfrak c_n^\theta$ through
	\eqref{eq:theta-hamiltonian-observation-state}--
	\eqref{eq:theta-hamiltonian-observation-parameter}.
	
	We now eliminate the forward sensitivities using the adjoints. More
	explicitly, \eqref{eq:theta-backward-adjoint} and
	\eqref{eq:theta-forward-sensitivity-scheme} imply
	\begin{align*}
		(Q_k^\theta)^\top\nabla_\beliefstate\mathfrak f_k\,\Delta t_k
		={}&(Q_k^\theta)^\top p_k^\theta
		-(Q_{k+1^-}^\theta)^\top p_{k+1^-}^\theta
		+(D_\theta\widehat\Phi_k^\theta)^\top p_{k+1^-}^\theta.
	\end{align*}
	Likewise, \eqref{eq:theta-preupdate-adjoint} and
	\eqref{eq:theta-forward-sensitivity-update} give
	\[
	(Q_{k_n^-}^\theta)^\top\nabla_\beliefstate\mathfrak c_n^\theta
	=(Q_{k_n^-}^\theta)^\top p_{k_n^-}^\theta
	-(Q_{k_n}^\theta)^\top p_{k_n}^\theta
	+(D_\theta\overline\Psi_n^\theta)^\top p_{k_n}^\theta.
	\]
	Substituting these identities into \eqref{eq:theta-diff-objective}, all
	interior sensitivity terms telescope. The terminal term cancels because
	$p_K^\theta=\nabla_\beliefstate\mathfrak g$, and the initial boundary
	term vanishes because $Q_0^\theta=0$. Thus
	\begin{align*}
		\nabla_\theta\mcal J_\Delta(\widehat U^\theta)
		=\mathbb E\Bigg[&
		\sum_{k=0}^{K-1}
		\Bigl(D_\theta\mathfrak f_k\,\Delta t_k
		+(D_\theta\widehat\Phi_k^\theta)^\top p_{k+1^-}^\theta\Bigr)\\
		&+\sum_{n=1}^{N_o}
		\Bigl(D_\theta\mathfrak c_n^\theta
		+(D_\theta\overline\Psi_n^\theta)^\top p_{k_n}^\theta\Bigr)
		+D_\theta\mathfrak g\Bigg].
	\end{align*}
	Replacing the cost slice derivatives by their Hamiltonian expressions
	gives \eqref{eq:theta-dual-gradient}. This proves the sole proposition
	directly, while also recovering the reduced-map representation.
\end{proof}

\subsection{Gaussian moments and the Kalman update}
\label{app:gaussian-reduction-proof}
\label{sec:gaussian-formulas}

All numerical examples use exact Gaussian filtering for their linear-Gaussian
state and observation models, with time-discretization error retained.
They use the reduced belief state
\[
\beliefstate_k=(m_k,\Sigma_k),
\qquad
\mu^{\beliefstate_k}=\mcal N(m_k,\Sigma_k).
\]
For $X\sim\mcal N(m,\Sigma)$, the raw moment is
\begin{equation}
	M_n(m,\Sigma)
	=
	\E[X^n]
	=
	\sum_{q=0}^{\lfloor n/2\rfloor}
	\binom{n}{2q}(2q-1)!!\,
	m^{n-2q}\Sigma^q,
	\qquad n\geq0.
	\label{eq:gaussian-raw-moments}
\end{equation}
with $(-1)!!=1$.

For the polynomial potential~\eqref{eq:general-poly-potential}, define the gradient moment
\begin{equation}
  \eta_k^\theta
  :=\int_\R \partial_x \widehat U_k^\theta(x,\yhist^{[\nu(\tau_k)]})\,\mu^{\beliefstate_k}(\dd x).
  \label{eq:eta-general}
\end{equation}
Using~\eqref{eq:general-poly-Ux} and~\eqref{eq:gaussian-raw-moments},
\begin{equation}
  \eta_k^\theta
  =\sum_{j=1}^{p_x}\sum_{\ell=0}^{L_{\nu(\tau_k)}}
  j\theta^{(\nu(\tau_k))}_{k,j,\ell}
  M_{j-1}(m_k,\Sigma_k)\varphi^{(\nu(\tau_k))}_\ell(\yhist^{[\nu(\tau_k)]}).
  \label{eq:eta-general-explicit}
\end{equation}
The derivatives
$\eta_{m,k}=\partial_m\eta_k^\theta$ and
$\eta_{\Sigma,k}=\partial_\Sigma\eta_k^\theta$ are obtained by differentiating the raw moments.

In the scalar linear dynamics
\begin{equation}
  \dd X_t=(aX_t+b\alpha_t)\,\dd t+\sigma\,\dd W_t,
  \label{eq:scalar-linear-dynamics}
\end{equation}
With quadratic control cost $\frac12 r\alpha^2$, minimizing the
continuous Hamiltonian over $\alpha$ yields
\begin{equation}
  \widehat\alpha_k^\theta=-\frac{b}{r}\eta_k^\theta,
  \qquad
  \rho:=\frac{b^2}{r}.
  \label{eq:gaussian-eliminated-control}
\end{equation}
The reduced running cost slice is
\begin{equation}
\mathfrak f_k(\beliefstate_k,\widehat U_k^\theta)
=
\frac12 q(m_k^2+\Sigma_k)
+
\frac12\rho(\eta_k^\theta)^2.
  \label{eq:gaussian-running-cost}
\end{equation}
The Euler prediction map is
\begin{align}
  m_{k+1^-}&=m_k+(a m_k-\rho\eta_k^\theta)\Delta t_k,
  \label{eq:gaussian-pred-m}\\
  \Sigma_{k+1^-}&=\Sigma_k+(2a \Sigma_k+\sigma^2)\Delta t_k.
  \label{eq:gaussian-pred-s}
\end{align}
For a scalar observation with fixed noise,
\begin{equation}
  Y_{t_n}=cX_{t_n}+\varepsilon\xi_n,
  \qquad \xi_n\sim\mcal N(0,1),
  \label{eq:fixed-noise-observation}
\end{equation}
the predictive observation variance and Kalman gain are
\begin{equation}
  v_n=c^2\Sigma_{t_n^-}+\varepsilon^2,
  \qquad
  K_n=\frac{c\Sigma_{t_n^-}}{v_n}.
  \label{eq:kalman-v-K}
\end{equation}
The update is
\begin{align}
  m_{t_n}&=m_{t_n^-}+K_n(y_n-cm_{t_n^-}),
  \label{eq:kalman-update-m}\\
  \Sigma_{t_n}&=(1-K_nc)\Sigma_{t_n^-}
  =\frac{\varepsilon^2}{c^2\Sigma_{t_n^-}+\varepsilon^2}\Sigma_{t_n^-}.
  \label{eq:kalman-update-s}
\end{align}
Finally, the terminal cost used in the scalar quadratic examples is
\begin{equation}
  \mathfrak g(\beliefstate_K)=\frac12 q_T\big((m_K-x_\star)^2+\Sigma_K\big).
  \label{eq:gaussian-terminal-cost}
\end{equation}

\paragraph{Derivation of the moment and conditioning formulas.}

Let $X\sim\mcal N(m,s)$ and write $X=m+\sqrt{s}Z$ with $Z\sim\mcal N(0,1)$. Then
\begin{equation}
  \E[X^n]=\sum_{p=0}^n\binom{n}{p}m^{n-p}s^{p/2}\E[Z^p].
\end{equation}
The odd moments of $Z$ vanish and $\E[Z^{2q}]=(2q-1)!!$. Therefore
\begin{equation}
  \E[X^n]
  =\sum_{q=0}^{\lfloor n/2\rfloor}\binom{n}{2q}(2q-1)!!m^{n-2q}s^q,
\end{equation}
which is \eqref{eq:gaussian-raw-moments}.

For the Kalman update, suppose
\begin{equation}
  Y=cX+\varepsilon\xi,
  \qquad \xi\sim\mcal N(0,1),
\end{equation}
with $\xi$ independent of $X$. The joint law of $(X,Y)$ is Gaussian with
\begin{equation}
  \E[Y]=cm,
  \qquad
  \operatorname{Var}(Y)=c^2s+\varepsilon^2,
  \qquad
  \operatorname{Cov}(X,Y)=cs.
\end{equation}
Conditioning a joint Gaussian gives
\begin{equation}
  \E[X\mid Y=y]=m+\frac{cs}{c^2s+\varepsilon^2}(y-cm),
\end{equation}
and
\begin{equation}
  \operatorname{Var}(X\mid Y)=s-\frac{c^2s^2}{c^2s+\varepsilon^2}
  =\frac{\varepsilon^2}{c^2s+\varepsilon^2}s.
\end{equation}
These are \eqref{eq:kalman-update-m}--\eqref{eq:kalman-update-s}.

\subsection{Scalar LQG continuous Hamiltonian reduction}
\label{app:lqg-continuous-proof}

For a smooth potential $U$ and a Gaussian belief, define
\begin{equation}
  \eta=\ip{\partial_xU}{\mu}.
\end{equation}
In the scalar linear-quadratic model, the part of the continuous Hamiltonian that depends on the continuous action is
\begin{equation}
  b a\,\eta+\frac{r}{2}a^2.
\end{equation}
This quadratic expression is minimized at $a=-(b/r)\eta$, and its value is $-b^2\eta^2/(2r)$. The induced drift of the belief mean is therefore
\begin{equation}
  \dot m=a_0m-\frac{b^2}{r}\eta,
\end{equation}
where $a_0$ denotes the linear drift coefficient in the state equation. The variance equation is independent of the action in the additive-noise LQG setting and is
\begin{equation}
  \dot s=2a_0s+\sigma^2.
\end{equation}
Applying forward Euler gives \eqref{eq:gaussian-pred-m}--\eqref{eq:gaussian-pred-s}. The running-cost slice follows from the envelope identity.

\subsection{Explicit controlled-sensing formula in scalar regime B}
\label{app:regime-b-beta-proof}

Regime B uses one controlled observation
\begin{equation}
  Y=cX^-+\beta\xi,
  \qquad \xi\sim\mcal N(0,1),
\end{equation}
with sensing cost $\kappa_1/\beta$. Before the observation, the reduced state is $(m^-,s^-)$. The post-observation potential is quadratic in both $x$ and the realized observation $y$:
\begin{equation}
  U^+(x,y)=\sum_{j=0}^{2}\sum_{\ell=0}^{2}\theta^+_{j,\ell}x^jy^\ell.
\end{equation}
The pre-posterior continuation value is
\begin{equation}
  \mcal C^\theta(\beta;m^-,s^-)
  =\E\left[\E\left[U^+(X^-,Y)\mid Y\right]\right]
  =\E[U^+(X^-,Y)].
\end{equation}
Expanding $Y^\ell=(cX^-+\beta\xi)^\ell$ for $\ell=0,1,2$ shows that every term is independent of $\beta$ or proportional to $\beta^2$. The odd terms in $\xi$ vanish. Therefore
\begin{equation}
  \mcal C^\theta(\beta;m^-,s^-)=K_0^\theta(m^-,s^-)+K_1^\theta(m^-,s^-)\beta^2.
\end{equation}
The reduced sensing objective is
\begin{equation}
  H_1^\theta(\beta;m^-,s^-)
  =\frac{\kappa_1}{\beta}+K_0^\theta(m^-,s^-)+K_1^\theta(m^-,s^-)\beta^2.
\end{equation}
When $K_1^\theta>0$ and the minimizer is interior, the first-order condition is
\begin{equation}
  -\frac{\kappa_1}{\beta^2}+2K_1^\theta\beta=0.
\end{equation}
Thus
\begin{equation}
  \beta_{\mathrm{unc}}^\theta=\left(\frac{\kappa_1}{2K_1^\theta}\right)^{1/3}.
\end{equation}
The constrained sensing rule is the projection of this value onto $[\beta_{\min},\beta_{\max}]$.

With several observations the same computation applies at each observation time $t_n$, with the potential evaluated at $t_n$ and cost $\kappa_n/\beta$. In Experiment~2, the potential~\eqref{eq:structured-potential} is quadratic in $(x,z)$, so $K_1^\theta$ equals $(T-t_n)\operatorname{softplus}(c_\theta(t_n))$ and does not depend on $(m^-,s^-)$ or on the observations. The sensing rule is then deterministic. For potentials that are not of this form, the pre-posterior value is not of the form $K_0+K_1\beta^2$, and the sensing variable is selected on a finite grid of candidate values.

\subsection{Gaussian obstacle expectation}
\label{app:obstacle-proof}

In the obstacle example, the running penalty is supported on the spatial band $x_{\min}\le |x|\le x_{\max}$ during the active time interval. For $X\sim\mcal N(m,s)$,
\begin{align}
  \Pbb(x_{\min}\le |X|\le x_{\max})
  &=\Pbb(x_{\min}\le X\le x_{\max})
    +\Pbb(-x_{\max}\le X\le -x_{\min})\nonumber\\
  &=\Phi_{\mathrm N}\left(\frac{x_{\max}-m}{\sqrt s}\right)
    -\Phi_{\mathrm N}\left(\frac{x_{\min}-m}{\sqrt s}\right)
    +\Phi_{\mathrm N}\left(\frac{-x_{\min}-m}{\sqrt s}\right)
    -\Phi_{\mathrm N}\left(\frac{-x_{\max}-m}{\sqrt s}\right),
\end{align}
where $\Phi_{\mathrm N}$ is the standard normal distribution function. Multiplication by the penalty amplitude and by the time-window indicator gives the expected obstacle penalty in \eqref{eq:num-obstacle-expectation}. This formula avoids particle quadrature in the optimizer while still preserving the nonlinear dependence of the cost on the belief mean and variance.
\subsection{Smoothing bounds and gradients}
\label{app:smoothing-proof}

We quantify the Hamiltonian perturbations introduced in
\cref{subsec:smoothing}.
Write $\zeta=b(P_km_k+r_k)+\gamma\kappa$ and compare the
Hamiltonians at the same arguments.

For the square-root smoothing,
\begin{equation}
	0\le
	\sqrt{\zeta^2+\delta^2}-|\zeta|
	=
	\frac{\delta^2}
	{\sqrt{\zeta^2+\delta^2}+|\zeta|}
	\le\delta.
	\label{eq:app-sqrt-absolute-bias}
\end{equation}
Consequently,
\begin{equation}
	0\le
	\tH^c_k-\tH^{c,\mathrm{sqrt}}_{\delta,k}
	\le\alpha_{\max}\delta.
\end{equation}

For the log-sum-exp smoothing, the identity
\begin{equation}
	\psi_\delta^{\mathrm{soft}}(\zeta)
	=
	\alpha_{\max}|\zeta|
	+\delta\log\left(
	1+e^{-2\alpha_{\max}|\zeta|/\delta}
	\right)
\end{equation}
gives
\begin{equation}
	0\le
	\tH^c_k-\tH^{c,\mathrm{soft}}_{\delta,k}
	\le\delta\log 2.
	\label{eq:app-soft-absolute-bias}
\end{equation}
Both smoothed Hamiltonians therefore converge uniformly to the
hard Hamiltonian as $\delta\downarrow0$.

Differentiating the smoothed action terms gives
\begin{align}
	-(\psi_\delta^{\mathrm{soft}})'(\zeta)
	&=
	-\alpha_{\max}
	\tanh\left(\frac{\alpha_{\max}\zeta}{\delta}\right),
	\label{eq:app-soft-absolute-gradient}\\
	-(\psi_\delta^{\mathrm{sqrt}})'(\zeta)
	&=
	-\alpha_{\max}
	\frac{\zeta}{\sqrt{\zeta^2+\delta^2}}.
	\label{eq:app-sqrt-absolute-gradient}
\end{align}
At $\gamma=1$, these are precisely the controls in
\eqref{eq:ns-soft-control} and \eqref{eq:ns-sqrt-control}.
They remain admissible and converge to the hard selector
away from the switching surface.

Since $\partial_\gamma\zeta=\kappa$, differentiation at fixed
$(m_k,\Sigma_k,P_k,r_k)$ also yields
\begin{equation}
	\left.
	\partial_\gamma\tH^{c,\psi}_{\delta,k}
	\right|_{\gamma=1}
	=
	\frac q2(m_k^2+\Sigma_k)
	+\kappa\alpha_k^{\theta,\delta}.
\end{equation}
Thus the smoothed cost slice evaluates the original running cost
under the smoothed control, consistently with
\eqref{eq:ns-smoothed-cost-slice}.
Training minimizes this cost along the associated smoothed
belief dynamics. It does not minimize the smoothed Hamiltonian
value itself.

The bounds above concern Hamiltonian values at fixed arguments.
They do not directly bound the difference between the costs of
the induced policies or between their objective gradients.
In particular, control derivatives become concentrated near the
switching surface as $\delta$ decreases.
Differentiating the expected smoothed objective therefore still
requires the regularity and integrability assumptions stated
for the discrete adjoint.

\section{Potential parametrizations used in the experiments}
\label{sec:parametrization}

This appendix records the potential parametrizations used in the numerical experiments. A short overview is given at the beginning of Section~\ref{sec:numerics}.

\subsection{History-dependent polynomial potentials}

Let $n=\nu(\tau_k)$ be the number of observations available at time step $k$. We write a general history-dependent polynomial slice as
\begin{equation}
	\widehat U_k^\theta(x,\yhist^{[n]})
	=\sum_{j=0}^{p_x}\sum_{\ell=0}^{L_n}
	\theta^{(n)}_{k,j,\ell}\,x^j\,\varphi^{(n)}_\ell(\yhist^{[n]}),
	\label{eq:general-poly-potential}
\end{equation}
where $p_x$ is the state-polynomial degree and $\varphi^{(n)}_\ell$ are deterministic features of the realized observation history. The state derivative is
\begin{equation}
	\partial_x \widehat U_k^\theta(x,\yhist^{[n]})
	=\sum_{j=1}^{p_x}\sum_{\ell=0}^{L_n}
	j\theta^{(n)}_{k,j,\ell}\,x^{j-1}\,\varphi^{(n)}_\ell(\yhist^{[n]}).
	\label{eq:general-poly-Ux}
\end{equation}
In a single-observation experiment, this reduces to separate pre-observation and post-observation slices, with observation-polynomial degree $p_y$,
\begin{equation}
	\widehat U_k^{\theta,-}(x)=\sum_{j=0}^{p_x}\theta^-_{k,j}x^j,
	\qquad
	\widehat U_k^{\theta,+}(x,y)=\sum_{j=0}^{p_x}\sum_{\ell=0}^{p_y}\theta^+_{k,j,\ell}x^j y^\ell.
	\label{eq:single-observation-slices}
\end{equation}

\subsection{Compressed observation memory}

For numerical implementation, it is often convenient to keep only a compressed history feature $Z_{t_n}=H_n(\yhist^{[n]})$, for example the newest observation or a fixed-length queue of the most recent observations. Under the viewpoint used in this paper, $Z_{t_n}$ is not part of the Gaussian belief state. It is a deterministic feature supplied to the potential. Thus one may set
\[
\varphi^{(n)}_\ell(\yhist^{[n]})=\widetilde\varphi^{(n)}_\ell(Z_{t_n}),
\]
while keeping the reduced belief state equal to $\beliefstate_k=(m_k,\Sigma_k)$ in  Gaussian experiments.

Keeping only the newest observation discards information about the earlier ones. A longer memory or a learned summary feature can be used in the same framework, at the price of more coefficients.

\subsection{Structured potential for controlled-sensing LQG}
\label{subsec:structured-potential}

For the scalar LQG problem with controlled sensing, the potential is
\begin{equation}
	\widehat U^\theta(t,x,z)
	=\frac{q_T}{2}x^2+(T-t)\Bigl(a_\theta(t)\,x^2+b_\theta(t)\,x\,z+\operatorname{softplus}\bigl(c_\theta(t)\bigr)z^2\Bigr),
	\label{eq:structured-potential}
\end{equation}
where $z$ is the latest observation and $(a_\theta,b_\theta,c_\theta)$ are the three outputs of a feed-forward network with two tanh layers of width $32$ whose only input is the time. The potential is therefore quadratic in $(x,z)$ for each $t$, so that the Gaussian expectations $\E[U]$ and $\E[\partial_xU]$ are explicit in $(m,\Sigma)$. This is a controller ansatz for the LQG problem. It is not meant to represent the value function for every belief, and no reference solution is used to train or to select it.

\subsection{Polynomial--trigonometric potential for Experiment~3}
\label{subsec:poly-trig-potential}

Experiment~3 has three observations, so the potential has four slices $U_n$, $n=0,\dots,3$, where $n$ is the number of observations already received. We use a potential that is linear in its coefficients,
\begin{equation}
	U_n^\theta(t,x,y_1,\dots,y_n)
	=\frac{q_T}{2}x^2+(T-t)\sum_{c\in\{\mathrm{pol},\,\mathrm{trig}\}}\;\sum_{i}\;\sum_{\ell=0}^{2}
	\theta^{c}_{n,i,\ell}\,\phi^{c}_{n,i}(s)\,P_\ell(\tau),
	\label{eq:poly-trig-potential}
\end{equation}
with $s=(x/2.5,\,y_1/3,\dots,y_n/3)$, $\tau=2t/T-1$, and $P_\ell$ the Legendre polynomials. The factor $T-t$ keeps the terminal condition $U(T,\cdot)=\tfrac{q_T}{2}x^2$. The polynomial features $\phi^{\mathrm{pol}}_{n,i}$ are the monomials $s^a$ of even total degree $2\le|a|\le4$. The trigonometric features are $\phi^{\mathrm{trig}}_{n,i}(s)=\cos(k_i\cdot u)$ with $u=\pi\tanh(s)$, where $k_i$ runs over the nonzero integer vectors with $|k_i|_1\le2$, one from each pair $\pm k$. Both families are even under $(x,y)\mapsto(-x,-y)$, a symmetry of the problem. The feature family includes dependence on every available observation
$y_1,\dots,y_n$, so neither the potential nor its induced controls are
restricted to a truncated observation history.
Individual features may involve only a subset of the variables.
For randomly chosen coefficients, we checked numerically that
$\partial_{y_i}U_n\not\equiv0$ and
$\partial_x\partial_{y_i}U_n\not\equiv0$ for every $i$.
This checks the representational capacity of the family rather than
the dependence of every trained policy. Expectations of the trigonometric features under the Gaussian belief are computed by Gauss--Hermite quadrature, as for the network.

The four slices contain \AnsFeatPoly{} polynomial and \AnsFeatTrig{} trigonometric features, and each feature carries \AnsTimeCoef{} time coefficients. This gives \AnsParams{} trainable coefficients. Since $U^\theta$ is linear in $\theta$, the same family with a separate coefficient vector $\eta$ is used for the dual test family in Section~\ref{subsec:dual-maximisation}.

\subsection{Quadratic potentials for the nonsmooth bounded-control example}

The nonsmooth LQG example uses the quadratic time-dependent ansatz
\begin{equation}
	\widehat U_k^\theta(x)=\frac12 P_k x^2+r_kx+c_k^U.
	\label{eq:quadratic-potential}
\end{equation}
Only $(P_k,r_k)$ affect the induced bounded control and the reduced Hamiltonian. The constant $c_k^U$ is therefore omitted from the trained bundle in this restricted parametrization.

\section{Experimental settings}
\label{app:protocol}

Table~\ref{tab:parameter-counts} summarizes the potential
parametrizations and their numbers of trainable parameters.
\begin{table}[htbp]
	\centering
	\begin{tabular}{llr}
		\toprule
		Experiment & Potential parametrization & Trainable parameters \\
		\midrule
		1, $d=1$ & quadratic in $x$, $251$ time nodes & 753 \\
		1, $d=10$ & quadratic in $x$, $8$ time blocks & 528 \\
		1, $d=20$ & quadratic in $x$, $8$ time blocks & 1\,848 \\
		1, $d=50$ & quadratic in $x$, $8$ time blocks & 10\,608 \\
		1, $d=100$ & quadratic in $x$, $8$ time blocks & 41\,208 \\
		2 & time network, two width-$32$ hidden layers & 1\,219 \\
		3 & polynomial--trigonometric & 348 \\
		4 & polynomial in $(x,z)$ after observation & 2\,408 \\
		5 & $P_k,r_k$ on $60$ time steps & 120 \\
		\bottomrule
	\end{tabular}
	\caption{Numbers of trainable parameters in the potential
		parametrizations, including all allocated coefficients.
		The multidimensional obstacle comparison is excluded.}
	\label{tab:parameter-counts}
\end{table}

Experiment~1 uses Adam with learning rate $10^{-2}$,
$150$ iterations and batch size $512$ in the scalar case.
The multidimensional runs use learning rate $5\cdot10^{-3}$,
with $400$ iterations and batch size $128$ for $d=10$,
and $400$ iterations and batch size $64$ for
$d\in\{20,50,100\}$.
Conditional Gaussian moments are computed analytically.
Final costs use $10\,000$ evaluation scenarios in the scalar case
and $16\,384$ in the multidimensional cases.
The dual is optimized separately, starting from the trained potential.

Experiment~2 uses a uniform mesh with $350$ steps on $[0,1]$ and the
structured potential~\eqref{eq:structured-potential}.
The nominal observation dates $(0.25,0.5,0.75)$ are mapped by
$k_n=\operatorname{round}(350t_n)$ to indices $(88,175,262)$, with
ties rounded to the nearest even integer. Thus the implemented
assimilation dates are
\[
(88/350,175/350,262/350)\approx(0.251429,0.5,0.748571).
\]
The primal simulator, the dual calculation and the reference use these
same indices. The observation-stage neural coefficients are evaluated
at the nominal dates, as in the saved implementation.
Training uses $500$ Adam updates at learning rate $0.003$.
Gaussian quadrature computes the policy cost exactly in this
quadratic setting.
Sensing is selected on a grid of $41$ values during training, using a
straight-through softmin gradient with temperature $0.05$, and by the
closed-form rule at evaluation.
The reported checkpoint minimizes its own expected policy cost,
without reference targets.
The dual uses the same $1\,219$-parameter neural family, initialized from
the saved primal weights. Fitting uses $1\,000$ Adam updates with learning
rate $0.003$ and cosine decay, order-$7$ tensor Gauss--Hermite quadrature,
and order-$9$ validation every $50$ updates. Validation selects update
$950$. Subsequent L-BFGS refinement does not improve it.
Continuous relaxed controls are minimized analytically and the three
sensing actions numerically. Final evaluation uses $5\,000$ fresh
Monte Carlo scenarios and independent order-$9$, $11$ and $13$ quadrature
checks, giving values between $1.383758$ and $1.383768$.
The code archive includes the saved weights, traces and evaluation scenarios,
and a provenance record linking the checkpoint, evaluation and Figure~\ref{fig:lqg-regimeB-convergence}.

Experiment~3 uses $200$ time steps and the
polynomial--trigonometric potential~\eqref{eq:poly-trig-potential}.
Training uses $200$ AdamW updates at learning rate $0.01$,
with deterministic Gauss--Hermite quadrature and selection by
validation cost.
Policy costs use $5000$ innovation scenarios, shared with the
HJB policy for the paired comparison. The HJB simulation uses
$800$ time steps.
Dual optimization starts independently and uses $100$ stochastic
ascent updates, followed by evaluation on $5000$ fresh scenarios.

Experiment~4 uses a particle-regression initialization,
followed by Gaussian belief propagation with analytic conditional
moments.
Training uses $150$ Adam updates at learning rate $0.3$
and batch size $200$.
The primal cost is evaluated on $10\,000$ scenarios.
The dual is optimized separately using local multistart
minimization and evaluated on $2048$ fresh scenarios.

Experiment~5 compares hard- and soft-switch training with
smoothing parameter $\delta=0.25$.
Both use $150$ Adam updates at learning rate $0.05$
and batch size $256$.
Final costs are evaluated on $10\,000$ scenarios, with the
soft-trained potentials deployed using the exact hard switch.
The dual is optimized over the six coefficients entering the
observation penalties and evaluated on $4096$ fresh scenarios.

Reported Monte Carlo uncertainties are standard errors.
Numerical dual estimates also involve integration and solver errors,
which are not included in these standard errors.
\bibliographystyle{plain}
\bibliography{references}

\end{document}